\documentstyle[amssymb,amsfonts,11pt,epsfig]{amsart}

\input{xypic}

\newtheorem{TEO}{Theorem}[section]
\newtheorem{PROP}[TEO]{Proposition}
\newtheorem{LEM}[TEO]{Lemma}
\newtheorem{DEF}[TEO]{Definition}
\newtheorem{COR}[TEO]{Corollary}
\newtheorem{REM}[TEO]{Remark}

\newtheorem{FACT}[TEO]{Fact}

\newcommand\dual{\mathrel{\raise3pt\hbox{$\underline{\mathrm{\thinspace d
\thinspace}}$}}}

\newcommand\proj{\Bbb P}
\newcommand\Z{\Bbb Z}

\newcommand\R{\Bbb R}
\newcommand\Q{\Bbb Q}
\newcommand\Co{\Bbb C}

\def\Z{{\Bbb Z}}
\def\R{{\Bbb R}}
\def\Z{{\Bbb Z}}

\begin{document}
\title{Real Hyperelliptic Surfaces and the Orbifold Fundamental Group}
\author{Fabrizio Catanese\\ Paola Frediani}

\thanks{ Research performed in the realm of the  {\em EU Project
"E.A.G.E.R." }}
\date{}
\maketitle
\pagestyle{myheadings}
\markboth{FABRIZIO CATANESE \hspace{1cm}  PAOLA FREDIANI}{REAL HYPERELLIPTIC SURFACES AND THE ORBIFOLD...}

\begin{abstract}

\noindent In this paper we finish the topological classification of real algebraic surfaces of Kodaira dimension zero and we make a step towards the Enriques classification of real algebraic
surfaces, by describing in detail the structure of the moduli space of real
hyperelliptic surfaces.

Moreover, we point out the relevance in real geometry of the notion of the
orbifold fundamental group of a real variety, and we discuss related 
questions on real varieties $(X, \sigma)$ whose underlying complex 
manifold $X$ is a $K ( \pi, 1)$.

Our first result  is  that if $(S, \sigma)$ is a real 
hyperelliptic surface, then the differentiable type of the pair $(S, 
\sigma)$ is completely determined by the orbifold fundamental group 
exact sequence. \\
This result allows us to determine all the possible topological types of
$(S, \sigma)$, and to prove that they are exactly 78.

It follows also as a corollary  that there are exactly eleven cases 
for the topological type of the real part of S.

Finally,  we show that once we fix the topological type of $(S, \sigma)$
corresponding to a real hyperelliptic surface, the corresponding 
moduli space is irreducible (and connected).

We also give, through a series of
tables, explicit analytic representations of the 78 components of the
moduli space.

\vspace{0.2cm}

\end{abstract}

\section*{Introduction}

The purpose of this paper is twofold: on the one side, we finish the
topological classification  of real algebraic surfaces of Kodaira dimension
zero and we make a step towards the Enriques classification of real algebraic
surfaces, by describing in detail the structure of the moduli space of real
hyperelliptic surfaces; on the other hand, we point out the relevance in real
geometry of the notion of the orbifold fundamental group.

  In order to illustrate the latter concept, let us begin by answering the
reader's question: what is a real variety?

A smooth real variety is a pair  $(X, \sigma)$, consisting of the data of a
smooth complex manifold $X$ of complex dimension
$n$ and of an antiholomorphic involution $\sigma: X
\rightarrow X$ (an involution  $\sigma$ is a map whose square is the
identity).

The quickest explanation of what  antiholomorphic means goes as follows: the
smooth complex manifold $X$ is determined by the differentiable manifold $M$
underlying  $X$ and by a complex structure
$J$ on the complexification of the real tangent bundle of $M$.

If instead we consider the same manifold $M$ together with the complex
structure $-J$, we obtain a complex manifold  which is called the conjugate
of $X$ and denoted by $\bar{X}$.

The involution $\sigma$ is now said to be antiholomorphic if it provides an
isomorphism between the complex manifolds $X$ and $\bar{X}$ (and then $(X,
\sigma)$ and
$(\bar{X}, \sigma)$ are also isomorphic as pairs).

What are the main problems concerning real varieties? (we may restrict
ourselves to the case where X is compact).

\begin{itemize}
\item Describe  the isomorphism classes of such pairs $(X, \sigma)$.
\item Or, at least describe the possible topological or differentiable types
of the pairs $(X, \sigma)$.
\item At least describe the possible topological types for the real parts
$X': = X({\bf R}) = Fix(\sigma)$.
\end{itemize}

\begin{REM} Recall that Hilbert's 16-th problem is a special case of the last
question but for the more general case of a pair of real varieties
$(Z \subset X , \sigma)$.
\end {REM}

For a smooth real variety, we have the quotient double covering
$\pi : X \rightarrow Y = X/<\sigma>$, and the quotient $Y$ is called the Klein
variety of
$(X, \sigma)$.

In dimension $n=1$ the datum of the Klein variety is equivalent to  the datum
of the pair $(X, \sigma)$, but this is no longer true in higher dimension,
where we will need also to specify the covering $\pi$.

The covering  $\pi$ is  ramified on the so-called  real part of
$X$, namely,
$X': = X({\R}) = Fix(\sigma)$, which
  is either empty, or a real submanifold of real dimension $n$.

If $X': = X({\R}) = Fix(\sigma)$ is empty,
  the orbifold fundamental group of
$Y$ is just defined as
  the fundamental group of $Y$.

If instead $ X' \neq \emptyset $, we may take a fixed point $x_0 \in
Fix(\sigma) $ and observe that $\sigma$ acts on the fundamental group
$\pi_1(X, x_0)$: we can therefore define the orbifold fundamental group as
the semidirect product of the normal subgroup $\pi_1(X, x_0)$ with the cyclic
subgroup of order two generated by $\sigma$. It is easy to verify then that
changing the base point does not alter  the isomorphism class of the
following exact sequence, yielding the orbifold fundamental group as an
extension

$ 1 \rightarrow \pi_1(X) \rightarrow \pi_1^{orb}(Y) \rightarrow {\Z}/2
\rightarrow 1 $

(changing  the base point only affects the choice of a splitting of the above
sequence).

We claim that the orbifold fundamental group  exact sequence is a powerful
topological invariant of the pair
$(X, \sigma)$ in the case where $X$ has large fundamental group.

To illustrate a concrete issue of such a statement, let us consider the case
of a $K(\pi,1)$, i.e., the case where the universal cover of $X$ is
contractible: this case includes the case of complex tori, of hyperelliptic
surfaces and their generalizations, as well as the case of quotients of the
complex n-ball  and  of polydisks.

Then the homotopy type of $X$ is determined by the fundamental group
$\pi$ and some interesting quite general open questions are the following
ones (where, by abuse of language, we shall talk about orbifold fundamental
group of a real variety, to refer to the above orbifold fundamental group
exact sequence).

\begin{itemize}
\item If $X$ is a compact complex manifold  which is a $K(\pi,1)$, to what
extent does $\pi$ also determine the differentiable type?

\item Same assumptions as above and fix the group $\pi$: when is then the
moduli space of those manifolds $X$ irreducible or connected?

For instance, is it sufficient to put
  the further hypothesis that we only consider K\"ahler manifolds?
(a vast literature already exists, cf.  \cite{mos},
\cite{yau}, \cite{siu},\cite{j-y1}, cf. also \cite{ca} for more examples and
for further references)

\item Assume that we consider real varieties, $( X, \sigma)$ where $X$ is as
above  a $K(\pi,1)$: does the orbifold fundamental group determine the
differentiable type of the real variety?

\item Fixing the orbifold fundamental group, when do we get  a connected
moduli space?
\end{itemize}

Our purpose here is thus to give an issue where the classification of real
varieties can be given in terms of the orbifold fundamental group
(concerning the previous questions, even  the apparently easier case  of real
and complex tori appears to be unsettled as soon as the complex dimension
becomes at least 3, at least if we drop the K\"ahler hypothesis; however,
as it is going to be shown elsewhere,  the complex
manifolds isomorphic to the product of a curve $C$ of genus at least $3$
and of a complex torus of dimension $2$ form a connected component of
their moduli space, which admits however another component,
corresponding to  the examples of \cite{cal}, providing on the underlying
differentiable manifold a non K\"ahler complex structure with trivial
canonical bundle).

Our first result, concerning the topological type of real hyperelliptic
surfaces, can now be briefly stated as follows:

\begin{TEO} Let $(S, \sigma)$ be a real hyperelliptic surface. Then the
differentiable type of the pair $(S, \sigma)$ is completely determined by the
orbifold fundamental group exact sequence.
\end{TEO}

In a sequel to this paper, we plan to show other issues (e.g., in the Kodaira
classification of non algebraic real surfaces) where the topology of the pair
$(S, \sigma)$ is determined by the orbifold fundamental group exact sequence.

Returning  to the case of real hyperelliptic surfaces, the previous theorem
allows us to easily determine completely the possible topological types of
$(S, \sigma)$, and in particular we have the

\begin{TEO} Real hyperelliptic surfaces fall into exactly $78$ topological
types.
\end{TEO}

and the

\begin{COR}

Let $(S, \sigma)$ be a real hyperelliptic surface. Then the real part
  $S(\Bbb R))$ is either
  \begin{itemize}
  \item
  a disjoint union of $t$ tori, where $ 0 \leq t \leq 4$
  \item
  a disjoint union of $b$ Klein bottles,
  where  $ 1 \leq b \leq 4$.
  \item
  the disjoint union of one torus and one Klein bottle
  \item
  the disjoint union of one torus and  two Klein bottles.
  \end{itemize}
\end{COR}

As the reader may guess, the above results are too complicated to be
  described in detail here in the introduction: therefore we will limit
ourselves to illustrate the underlying philosophy by describing it in the
much simpler case of the real elliptic curves.

Classically (cf. e.g. \cite{a-g} for a modern account) real elliptic curves
have been classified according to the number $\nu$ of connected components
(these are circles) of their real part: $\nu$ can only attain the values
$0,1,2$ and completely determines the differentiable type of the involution.

The orbifold fundamental group explains easily this result: if there is a
fixed point for $\sigma$, the orbifold fundamental group sequence splits and
the action of $\sigma$ on the elliptic curve $C$ is completely determined by
the action $s$ of $\sigma$ on $H_1(C, \Z)$.

This situation gives rise to only two cases: $s$ is diagonalizable, and
$C(\R)$ consists of two circles, or $s$ is not diagonalizable, and $C(\R)$
consists of only one circle.

If instead there are no fixed points, an easy linear algebra argument (cf.
lemma \ref{tra}) shows that $s$ is diagonalizable, and the translation vector
of the affine transformation inducing $\sigma$ can be chosen to be
$1/2$ of the $+1$-eigenvector $e_1$ of $s$.

In fact, $\sigma$ is represented by an affine transformation
$(x,y) \rightarrow s(x,y) + (a,b) $, and $s$ is not diagonalizable if and
only if $s(x,y)= (y,x)$ for a suitable choice of two basis vectors. The
identity map, which equals the square of
$\sigma$, is the transformation $(x,y) \rightarrow (x,y) + (a+b,a+b) $, thus
$a+b$ is an integer, and therefore the points $(x, x-a)$ yield a fixed $S^1$
on the elliptic curve.

The complete description of the moduli space of real hyperelliptic surfaces is
too long to be reproduced in the introduction, we want here only to mention
the following main result,  which confirms  a conjecture by
Kharlamov, that more generally for all real K\"ahler surfaces of Kodaira
dimension at most 1  the differentiable type of the pair $(S, \sigma)$ should
determine the deformation type (this is  false  already for complex surfaces
if the  Kodaira dimension  equals 2,  cf. \cite{ca1},\cite{ca2},
\cite{ma}).

\begin{TEO} Fix the topological type of $(S, \sigma)$ corresponding to a real
hyperelliptic surface. Then the moduli space of the real surfaces $(S',
\sigma ')$ with the given topological type is irreducible (and connected).
\end{TEO}

Again, we wish to give the flavour of the argument by outlining it in the much
simpler case of the elliptic curves. Assume for instance that our involution
$\sigma $ acts as follows: $ (x,y) \rightarrow (y,x) $. We look then for  a
translation invariant complex structure $J$ which makes
$\sigma $ antiholomorphic, i.e., we seek for the matrices $J$ with $J^2 = -1$
and with $ J s = - s J$.

The latter condition singles out the matrices
$\left(\begin{array}{cc} a & b \\ -b & -a
\end{array} \right)$

while the first condition is equivalent to requiring that the characteristic
polynomial be equal to $ \lambda^2 + 1$, whence, it is equivalent to the
equation
$b^2 - a^2 =1$.

We get therefore a hyperbola with two branches which are exchanged under the
involution
$ J \rightarrow -J$, but, as we already remarked,
$J$ and $-J$ yield isomorphic real elliptic curves, thus the moduli space
consists of just one branch of the hyperbola.

This example serves also the scope of explaining the statement in the above
theorem that the moduli space is irreducible (and connected): the hyperbola
is an irreducible  algebraic variety, but not an irreducible analytic space,
since it is not connected (in general, moduli spaces of real varieties will
be semianalytic
  spaces or semialgebraic real spaces).

We want to emphasize once more that an interesting question is to determine,
in the realm of the real varieties  whose topological type is determined by
the orbifold  fundamental group, those for which the corresponding moduli
spaces are irreducible (respectively: connected).

It is now however time to recall what the hyperelliptic surfaces are, why they have this name and, last but not least, point out how crucial  is the
role of the hyperelliptic surfaces in the Enriques classification of algebraic
surfaces.

As elliptic curves are exactly the curves such that the homogeneous coordinates of
their points cannnot be uniformized by polynomials, yet they can be
uniformized by entire holomorphic functions on $\Co$, hyperelliptic varieties
of dimension $n$ were generally defined by Humbert and Picard  through the
entirely analogous property that the coordinates of their points, although
  not  uniformizable by rational functions,  could  be uniformized by entire
meromorphic functions on ${\Co} ^n$.

Among these varieties are clearly (nowadays) the Abelian varieties, and the
classification of such varieties in dimension two was achieved by Bagnera and
de Franchis who got the Bordin Prize in 1909 for their important result
(the classification by Enriques and Severi, who got the same prize for it
the year before,  had serious gaps which were corrected only later on).


The missing surfaces, which are now called hyperelliptic,  were described as
quotients of the product of two elliptic curves by the action of a finite
group G. For this reasons, some authors  call these surfaces bielliptic
surfaces (cf. \cite{be}).

The classification is in the end very simple and produces a list of 7 cases
where the Bagnera de Franchis group G and its action is explicitly written
down.

The reason to recall all this is that, as a matter of fact, an important
ingredient in the proof of our theorems is to rerun the arguments of the
proof of the classification theorem, which characterizes the hyperelliptic
surfaces as the algebraic surfaces with nef canonical divisor
$K$, with $K^2= p_g=0$, $q=1$,  and moreover with  Kodaira dimension equal to
$0$. This is done in section $1$.

To keep close in spirit to the beautiful result of Bagnera and de Franchis we
felt compelled to produce  tables exhibiting simple and explicit actions for
the 78 types of the real hyperelliptic surfaces: these are contained in the
last section, and they summarize a lot of information that we could not give
in a more expanded form.

Concerning now the Enriques classification of real algebraic surfaces, it has
been focused up to now mostly on  the classification of the topology of the
real parts, the topological classification of real rational surfaces going
back to Comessatti (\cite{co1}\cite{co2}\cite{co3}), as well as the
classification of real abelian varieties (\cite{co3}, see also \cite{si},
\cite{s-s}).

In the case of real $K3$ - surfaces we have the classification by Nikulin and
Kharlamov (\cite{ni},
\cite{kha}), for the real Enriques surfaces the one by Degtyarev and Kharlamov
(\cite{dekha1}, \cite{dekha3}).

Finally partial results on real ruled and elliptic surfaces have been
obtained by Silhol (\cite{si}) and by Mangolte (\cite{man3}).

The  paper is organized as follows:

in section $1$ we moreover recall the description given by Bagnera and de
Franchis of the hyperelliptic surfaces as quotients of a product of two
elliptic curves $E \times F$ by the product action of a finite group $G$
acting  on $E$ as a group of translations  and on $F$ via
an action whose quotient is ${\proj}^1$.

$G$ is called the Bagnera de Franchis group (or symmetry group) and is a
quotient of the fundamental group of the surfaces $S$.

In section $2$ we observe that the orbifold fundamental group has a finite
quotient  $\hat{G}$ which contains $G$ as a normal subgroup of index
$2$:
$\hat{G}$ is called the extended Bagnera de Franchis group and its structure
will be investigated in detail in section $4$.

In the rest of section $2$ we show that isomorphisms of real hyperelliptic
surfaces lift to isomorphisms of the respective products of elliptic curves,
compatibly with the identifications of the respective extended Bagnera de
Franchis groups.

Section $3$ is devoted first to showing that the  representation of the
orbifold fundamental group as a group of affine transformations of ${\Q}^4$
is uniquely determined, up to isomorphism, by the abstract structure of the
group, what proves  theorem 0.2 (by the way, we show in the course of the
proof a fact hardly mentioned in the literature, namely, that the
differentiable structure of a hyperelliptic surface is determined by the
fundamental group).

Second, after recalling quite briefly the notion of moduli spaces for real
varieties, we show that, once this affine representation is fixed, the moduli
space for the compatible complex structures is irreducible and connected.

Section $4$ recalls some known facts about antiholomorphic maps of elliptic
curves and applies these results to the determination of the possible
extended Bagnera de Franchis group which do in effect occur.

Section $5$ determines the analytical actions of these groups on the two
factors under the condition that the exact sequence

$ 1 \rightarrow G \rightarrow \hat{G} \rightarrow \Z /2 \rightarrow 1 $

splits, while section $6$ deals with the simpler case where there is no
splitting.

Section $7$ gives the recipe to identify the real part $S(\R)$ of our
surfaces as a disjoint union of Klein bottles and tori.

Finally, section $8$ applies the results developed insofar and achieves the
classification of the $78$ components of the moduli space, for which explicit
analytical representations, describing the action of the extended Bagnera de
Franchis group, are given through a series of tables.

\vspace{1cm}

{\bf Acknowledgements.} We would like to thank V. Kharlamov for an
interesting conversation concerning the status of the Enriques classification
of algebraic surfaces, and for communicating his conjecture to us. Finally,
S.T. Yau pointed out Calabi's example (\cite{cal})
  of diffeomorphic complex manifolds such that one is K\"ahler,
while the other is not: with some work, this example shows that,
without the K\"ahler hypothesis, the moduli space of complex manifolds which
are K($\pi , 1$)'s (for a fixed $\pi$) can be disconnected. The present
research was performed in the realm of the Schwerpunkt " Globale Methode in
der Komplexe Geometrie", and partially supported by the italian M.U.R.S.T.
40/100 program "Geometria Algebrica".\\

\section{Basics on hyperelliptic surfaces}

In this section we recall
(\cite{b-df}, \cite{b-df2} see also \cite{be}, Chapitre VI, pp. 91-115,
\cite{bpvdv} pp. 147-149) the definition of hyperelliptic surfaces and their
characterization in the realm of the Enriques classification of complex
algebraic surfaces. We shall also briefly recall the main lines of the proof
of the Bagnera - de Franchis classification theorem, since we shall
repeatedly need modified or sharper versions of the arguments used therein.

\begin{DEF} A complex surface $S$ is said to be hyperelliptic if $S \cong (E
\times F)/G$, where $E$ and $F$ are elliptic curves and $G$ is a finite group
of translations of $E$ with a faithful action on $F$ such that $F/G \cong
{\proj}^1$.
\end{DEF}

$G \subset Aut(F)$, so $G = T \rtimes G'$ (semidirect product), where $T$ is a
group of translations and $G' \subset Aut(F)$ consists of group automorphisms.
Since $F/G \cong {\proj}^1$, then $G' \neq 0$, hence
$G' \cong {\Z}/m$, with $m = 2,3,4,6$, by the following well known result.

\begin{FACT} Let $F$ be an elliptic curve. Every automorphism of $F$ is the
composite of a translation and a group automorphism. The non trivial group
automorphisms are the symmetry $x \mapsto -x$ and also:\\  for the curve $F_i
= {\Co}/({\Z } \oplus {\Z} \cdot i)$, $x \mapsto
\pm ix$.\\  For the curve $F_{\rho} = {\Co}/({\Z} \oplus {\Z} \cdot
\rho)$, where $\rho^3 =1 \neq \rho$, $x \mapsto \pm \rho x$, and $x \mapsto
\pm \rho^2 x$.
\end{FACT}

Since $G$ is abelian, as a group of translations of $E$, the product $T
\rtimes G'$ must be direct. We have the following result:

\begin{TEO}
\label{bdf} (Bagnera - de Franchis) Every hyperelliptic surface is one of the
following, where $E$, $F$ are  elliptic curves and $G$ is a group of
translations of $E$ acting on $F$ as specified:
\begin{enumerate}
\item $(E \times F)/G$, $G = {\Z}/2$ acting on $F$ by symmetry.
\item $(E \times F)/G$, $G = {\Z}/2 \oplus {\Z}/2$ acting on $F$ by $x
\mapsto -x$, $x \mapsto x +\epsilon$, where
$\epsilon$ belongs to the group $F_2$ of points of $F$ of order $2$.
\item $(E \times F_i)/G$, $G = {\Z}/4$ acting on $F_i$ by $x \mapsto ix$.
\item $(E \times F_i)/G$, $G = {\Z}/4 \oplus {\Z}/2$ acting on $F_i$ by
$x \mapsto ix$, $x \mapsto x + (1+i)/2$.
\item $(E \times F_{\rho})/G$, $G = {\Z}/3$ acting on $F_{\rho}$ by $x
\mapsto \rho x$.
\item $(E \times F_{\rho})/G$, $G = {\Z}/3 \oplus {\Z}/3$ acting on
$F_{\rho}$ by $x \mapsto \rho x$, $x \mapsto x + (1- \rho)/3$.
\item $(E \times F_{\rho})/G$, $G = {\Z}/6$ acting on $F_{\rho}$ by
$x \mapsto -\rho x$.
\end{enumerate}
\end{TEO}

Hyperelliptic surfaces are algebraic surfaces with $p_g=0$, $q =1$, $K^2 =0$,
$K$ nef.\\  We have the following basic result of the Enriques classification
of surfaces (Kodaira \cite{koiv} indeed proved that the same result holds more
generally for compact complex surfaces provided one replaces the hypothesis
$q = 1$ by $b_1 = 2$).
\begin{TEO}
\label{class}  The complex surfaces $S$ with $K$ nef, $K^2 =0$, $p_g =0$ and
such that either $S$ is algebraic with $q = 1$, or more generally $b_1 = 2$,
are hyperelliptic surfaces if and only if $kod(S) = 0$ (this is equivalent to
requiring that all the fibers of the Albanese map be smooth of genus 1).
\end{TEO}   {\bf Proof.}  Let $\alpha: S \longrightarrow A$ be the Albanese
map,
$q =1$, so $A$ is a curve of genus 1. Let $\pi: {\Co} \rightarrow A$ be the
universal covering and let us consider the pull - back diagram:

$$\diagram \tilde{S} \dto^{\tilde{\alpha}} \rto      & S \dto^{\alpha}  \\
            {\Co}  \rto^{\pi} & A
\enddiagram$$

All the fibers of $\alpha$ are smooth, therefore also the fibers of
$\tilde{\alpha}$ are smooth. So $\alpha$ and $\tilde{\alpha}$ are ${\cal
C}^{\infty}$ bundles, by Ehresmann's theorem. ${\Co}$ is contractible, so
$\tilde{S}$ is diffeomorphic to the product ${\Co} \times F$. Thus we obtain a
holomorphic map (since we have in fact a locally liftable holomorphic map $f :
{\Co} \rightarrow {\cal H}/PSL(2,{\Z})$, and ${\Co}$ is simply connected)
$$f: {\Co} \rightarrow {\cal H} = \{\tau \in {\Co} \ | \ Im(\tau) >0\},$$
$$ t \mapsto \tau,$$ where  $$\tilde{\alpha}^{-1}(t) \cong {\Co}/{\Z}
\oplus \tau {\Z}.$$   By Liouville's theorem we see that $f$ is constant,
therefore  $$\tilde{S}
\cong {\Co} \times F.$$  If $A = {\Co}/\Lambda$, then $S \cong ({\Co}
\times F)/\Lambda$, where
$\Lambda$ acts on ${\Co}$ by translations and on $F$ by a map $\mu: \Lambda
\rightarrow Aut(F)$. \\   Let $ \Gamma := Aut(F)/Aut^0(F)$, where
$Aut^0(F)$ are the automorphisms which are homotopic to the identity, and let
$\nu:
\Lambda \rightarrow \Gamma$ be the induced map. We set $\Lambda' := ker \nu$,
$\Gamma' := \nu(\Lambda)$. Then we get the following pull-back diagram:

$$\diagram  S':= \tilde{S}/\Lambda' \dto^{\alpha'} \rto^{\psi}    & S =
S'/\Gamma' \dto^{\alpha}  \\
             A' := {\Co}/\Lambda'  \rto & A = A'/\Gamma'
\enddiagram$$

$\Lambda'$ acts by translations on $F$, so it acts as the identity on
$H^0(\Omega^1_F)$. Then there exist $\eta, \eta' \in H^0(\Omega^1_{S'})$ such
that $\eta \wedge \eta' \not \equiv 0$, thus $K_{S'} \equiv 0$ and $q(S')
=2$, so
$S'$ is a complex torus. The map $\psi: S' \rightarrow S$ is an unramified
covering of degree $m$, where $m \in \{2,3,4,6\}$. Then $mK_S =
\psi_{*} \psi^{*} K_S = \psi_{*} K_{S'} = 0$. In particular $12K_S \equiv 0$.
  The map $\alpha' : S' \rightarrow A'$ is a fibre bundle on an elliptic curve
and
$S'$ is a complex torus of dimension 2. Since $b_1 = 2$, $S$ is algebraic
hence
$S'$ is algebraic too and, by Poincar\'e's reducibility theorem there exists a
finite unramified covering $A'' \rightarrow A'$ yielding a product structure
on the pull back $S''$ of $S'$.

$$\diagram  S'':= A'' \times_{A'} S' \cong  A'' \times F \dto^{\alpha''}
\rto    & S' \dto^{\alpha'}  \\
             A'' =: {\Co}/\Lambda''  \rto & A' = {\Co}/\Lambda'
                                       \enddiagram$$

Then we have found that $S \cong (A'' \times F)/G$, where
$G=\Lambda/\Lambda''$. \\ Moreover, by choosing the covering $A''
\rightarrow A'$ minimal with the above property, one sees that $E:=A''$, $F$,
$G$ are as in the list by Bagnera de Franchis.
\hfill Q.E.D. \\

\section{Real conjugations on hyperelliptic surfaces}

Let us now suppose that $S$ is a real hyperelliptic surface, i.e. there is an
antiholomorphic involution $\sigma: S \rightarrow S$, and we consider the
isomorphism class of the pair $(S, \sigma)$.\\ Since, by definition of the
Albanese map $\alpha$, fixed a point $x_0 \in S$,
$\alpha(x) = \int_{x_0}^{x}$, we obtain
$$ \alpha(\sigma(x)) = \int_{x_0} ^{\sigma(x)} = \int_{x_0}^{\sigma(x_0)} +
\int_{\sigma(x_0)}^{\sigma(x)},$$ if we define
$$\bar{\sigma}(\gamma) := \int_{x_0}^{\sigma(x_0)} + \sigma_* (\gamma),$$ we
get an induced antiholomorphic map on the Albanese variety $\bar{\sigma} : A
\rightarrow A$ with the property that the following diagram commutes

$$\diagram  S \dto^{\alpha} \rto^{\sigma}    & S \dto^{\alpha}  \\
             A = {\Co}/\Lambda  \rto^{\bar{\sigma}} & A = {\Co}/\Lambda
                                       \enddiagram$$

A direct calculation, or the remark that $\alpha(S)$ generates $A$ and
$\bar{\sigma}^2$ is the identity on $\alpha(S)$, assures that $\bar{\sigma}$
is an antiholomorphic involution on $A$. \\ Notice that, is $S({\R}) \neq
\emptyset$, we may choose a point $x_0$ with
$\sigma(x_0) = x_0$, and then $\bar{\sigma}$ will be a group homomorphism.
\\ We want to prove that $\sigma$ lifts to a map $\tilde{\sigma}: A'' \times F
\rightarrow A'' \times F$, where $A''$ and $F$ are as in the proof of the
previous theorem. \\ Observe that since we have the pull-back diagram

$$\diagram  S'' \cong  A'' \times F \dto^{\alpha''} \rto^{\phi}    & S
\dto^{\alpha}  \\
             A'' = {\Co}/\Lambda''  \rto & A = {\Co}/\Lambda
  \enddiagram$$

it suffices to prove that the involution $\bar{\sigma}$ on $A$ lifts to
$A''$.\\ In fact then $\sigma$ lifts to $S''$ as a fibre product and so we
have an induced action on  $A'' \times F$, preserving $\alpha''$.\\  We need
the following
\begin{LEM}  Let $\pi: Y \longrightarrow X$ be a connected covering space, and
let $g$ be a homeomorphism of $X$. Choose $x_0 \in X$, $y_0 \in Y$ with
$\pi(y_0) = x_0$, and let $z_0 = g(x_0)$, $w_0 \in Y$ with $\pi(w_0) = z_0$.
\\  Then there exists a lift $\tilde{g}$ of $g$ with $\tilde{g}(z_0) = w_0$ if
and only if, given a path $\tilde{\delta}$ from $y_0$ to $w_0$ and setting
$\delta = \pi \circ \tilde{\delta}$, and considering the isomorphism $\Delta:
\pi_1(X, z_0) \stackrel{\cong} \rightarrow \pi_1(X,x_0)$ such that
$\Delta(\gamma) = \delta \gamma \delta^{-1}$, and similarly
$\tilde{\Delta}$, we have
$$\Delta g_* (H) = H, \ where \ H = H_{y_0} = \pi_*(\pi_1(Y, y_0)).$$
\end{LEM}  {\bf Proof.}  Consider the diagram of pointed spaces

$$\diagram  (Y,y_0) \dto^{\pi}     & (Y,w_0) \dto^{\pi}  \\
             (X,x_0)  \rto^{g} & (X,z_0)
                                       \enddiagram$$ Then $\tilde{g}$ exists
and it is unique if and only if $g_*(H_{y_0}) = H_{w_0}$; applying $\Delta$,
iff
$$\Delta g_* (H_{y_0}) = \pi_*(\tilde{\Delta} \pi_1(Y,w_0)).$$  But
$\tilde{\Delta} \pi_1(Y,w_0) = \pi_1(Y,y_0)$, whence $\tilde{g}$ exists and
it is unique iff $\Delta g_*(H) = H$.
\hfill Q.E.D. \\

\begin{COR} Under the above notation there exists a lift $\tilde{g}$ of $g$ if
and only if
$\Delta g_* (H)$  is a conjugate of $H$.
\end{COR} {\bf Proof.} There exists $\tilde{g}$ if and only if there exists
$w_0 ' \in \pi^{-1}(z_0)$ such that $\Delta' g_{*}(H) = H$. This is
equivalent to say that $\Delta g_{*} (H)$ is a conjugate of $H$ (we conjugate
by $\delta {\delta'}^{-1}$, where $\delta' = \pi \circ \tilde{\delta}'$ and
$\tilde{\delta}' $ is a path from $y_0$ to $w'_0$).
\hfill Q.E.D. \\

\begin{COR}
\label{fixlift} There exists a lift $\tilde{g}$ with a fixed point if and
only if there exists $x_0 \in Fix(g)$ and a conjugate subgroup $H'$ of $H$
such that
$g_*(H') =H'$. \\   Furthermore then, if $g$ has order $n$ then also
$\tilde{g}$ has order $n$ (since $\tilde{g}^n$ is a lift of the identity and
has a fixed point).
\end{COR}   Now we go back to our situation. \\ Since the fundamental groups
of
$S'$, resp. $S''$ give rise to subgroups of
$\pi_1(S)$ which are the pull backs of $\pi_1(A')$, resp. $\pi_1(A'')$ under
the bundle homotopy exact sequence
$$1 \rightarrow \pi_1(F) \rightarrow \pi_1(S) \rightarrow \pi_1(A) \rightarrow
1,$$ we obtain that they are preserved under $\sigma_*$ if and only if the
corresponding subgroups of $\pi_1(A)$ are preserved under
$\bar{\sigma}_*$.\\ In the latter case all the fundamental groups are abelian
and we need to prove that
$\bar{\sigma}$ lifts to $A''$, or equivalently that $\Lambda'' $ is
$\bar{\sigma}$ - invariant. \\ We observe that $\Lambda'$ is invariant by
$\bar{\sigma}$, since $\Lambda'$ is the kernel of the topological monodromy
$\nu: \Lambda \rightarrow \Gamma = Aut(F)/Aut^0(F)$, which is induced by
$\sigma_*$, hence it is $\bar{\sigma} $ - equivariant. \\ We observe that
while
$\Lambda'$ is canonically defined as the kernel of the topological monodromy
of
$\alpha$, it is a priori not clear that $\Lambda''$ be canonically defined (we
shall indeed prove later that $\Lambda''$ is the centre of the fundamental
group of $S$).\\ However, using the list of Bagnera - de Franchis, case by
case, we see that
$\Lambda''$ is a characteristic subgroup of $\Lambda$ and therefore that it is
$\bar{\sigma}$ - invariant.  \\

In the cases 1, 3, 5, 7 of the list of Bagnera - de Franchis (\ref{bdf})
there is nothing to prove, since we have $\Lambda'' = \Lambda'$. \\ In case 2
we find
$\Lambda'' = 2 \Lambda$. In fact  we have $G = {\Z}/2 \times {\Z}/2 $, acting
on
$A'' \times F$ as:
$$(x_1,x_2) \mapsto (x_1 +\eta, -x_2),$$
$$(x_1,x_2) \mapsto (x_1 + \eta', x_2 + \epsilon).$$

Analogously in case 6, $\Lambda'' = 3 \Lambda$, in  case 4, $\Lambda'' =
2\Lambda
\cap \Lambda'$. Therefore we always find that $\sigma(\Lambda'') =
\Lambda''$ and thus we get a lift
$$\diagram  A'' \times F  \dto
\rto^{\tilde{\sigma}}   & A'' \times F \dto  \\             S  \rto^{\sigma} &
S                                                 \enddiagram$$

\begin{DEF}
\label{hatg} The extended symmetry group $\hat{G}$ is the group generated by
$G$ and
$\tilde{\sigma}$.
\end{DEF}
$\hat{G}$ is the group of homeomorphisms of $S''$ which lift the group $\{1,
\sigma\}$: hence we have the following extension that will be studied in the
next section
$$(*) \ \ 0 \rightarrow G \stackrel{i} \rightarrow \hat{G} \stackrel{\pi}
\rightarrow {\Z}/2 \rightarrow 1$$ With the same arguments as above we obtain
the following

\begin{TEO}
\label{ISO1}  Let $(S, \sigma)$, $(\hat{S}, \hat{\sigma})$ be isomorphic real
hyperelliptic surfaces (i.e. there exists $\psi: S  \stackrel{\cong}
\rightarrow \hat{S}$ such that $\psi^{-1} \hat{\sigma} \psi = \sigma$). Then
the respective extended symmetry groups $\hat{G}$ are the same for $S$ and
$\hat{S}$. Moreover let  $S = (E \times F)/G$, $\hat{S} = (\hat{E} \times
\hat{F})/G$ be two Bagnera - De Franchis realizations. Then there exists an
isomorphism
$\Psi: E \times F \rightarrow \hat{E} \times \hat{F}$, of product type (i.e.
$\Psi = \Psi_1 \times \Psi_2$) commuting with the action of $\hat{G}$, and
inducing the given isomorphism $\psi: S \stackrel{\cong} \rightarrow \hat{S}$.
\end{TEO} {\bf Proof.}
$\psi$ induces an isomorphism $\psi_*$ of the Albanese varieties, which is
compatible with the antiholomorphic involutions $\bar{\sigma}$ and
$\hat{\bar{\sigma}}$ and such that we have a commutative diagram

$$\diagram   S \dto^{\alpha} \rto^{\psi} & \hat{S}  \dto^{\hat{\alpha}}  \\ A
\rto^{\psi_*} & \hat{A}
\enddiagram$$

Whence, by taking the coverings associated to the subgroups and points
corresponding under $\psi$ and $\psi_*$, we obtain isomorphisms
$\tilde{\psi}$, $\psi_1$ and a commutative diagram

$$\diagram   S'' \dto^{\alpha''} \rto^{\tilde{\psi}} & \hat{S}''
\dto^{\hat{\alpha}''}  \\ A'' \rto^{\psi_1} & \hat{A}''
\enddiagram$$

We observe that if $(E \times F)/G$ is a Bagnera - de Franchis realization of
$S$, then $E \cong A''$ and there is an isomorphism of $(S'' \rightarrow A'')$
with $(E \times F \rightarrow E)$, commuting with the action of $G$. \\ We
obtain therefore a commutative diagram

$$\diagram   E \times F \dto \rto^{\tilde{\psi}} & \hat{E} \times \hat{F}
\dto  \\  E \rto^{\psi_1} & \hat{E}
\enddiagram$$ where, moreover, $\tilde{\psi}$ and $\psi_1$ are real. \\ Since
$\tilde{\psi}$ preserves the fibres of the two projections, we have
$$\tilde{\psi} (e,f) = (\psi_1(e), \psi_2(e,f)).$$ Let us fix an origin $0
\in E$, then, since $\psi_2$ is an affine map, we can write
$$\psi_2(e,f) = \psi_2(0,f) + r(e),$$ where $r$ is a holomorphic homomorphism
$r: E \rightarrow Pic^0(\hat{F})$ ($Pic^0(\hat{F})$ is the group of
translations of $\hat{F}$). \\ By $G$ - equivariance, for all $g \in G
\subset Pic^0(E)$, we have
\begin{equation}
\label{pro}
\psi_2(e + g, g(f)) = g(\psi_2(0,f) + r(e)).
\end{equation}   But the left hand side of \ref{pro} equals
$$\psi_2(0, g(f)) + r(e + g) = \psi_2(0, g(f)) + r(e) + r(g),$$ thus if we let
$g_*$ be the linear part of $g: \hat{F} \rightarrow \hat{F}$, and we look at
the linear part of \ref{pro} with respect to $e$, we obtain
$$g_*(r(e)) \equiv r(e).$$ Since there is a $g \in G$ such that $g_* \neq
Identity_{\hat{F}}$, we infer that $r(e)$ is constant, or equivalently that
$\psi_2(e,f) = \psi_2(f)$.
\hfill Q.E.D. \\

\begin{PROP}
\label{separa} Let $(S = (E \times F)/G, \sigma)$ be a real hyperelliptic
surface, and let
$\tilde{\sigma} : E \times F \rightarrow E \times F$ be a lift of $\sigma$.
Then the antiholomorphic map $\tilde{\sigma}$ is of product type.
\end{PROP}  {\bf Proof}.   As in the proof of the previous theorem, we have
$$\tilde{\sigma}(e,f)  = (\sigma_1(e), \sigma_2(e,f))$$ since
$\tilde{\sigma}$ preserves the fibration onto $E$. Then, after choosing an
origin
$0 \in E$, we have
$$\sigma_2(e,f) = \sigma_2(0,f) + r(e),$$  where $r : E \rightarrow Pic^0(F)$
is an antiholomorphic homomorphism. \\ We know that $\tilde{\sigma}$
normalizes the group $G$, in particular  if we take some element  $g \in G$,
such that
$g_*
\neq Id_{F}$, $\tilde{\sigma}$ has a matrix
$\left( \begin{array}{cc} a & 0 \\ b & c
\end{array} \right)
$ while $g_* =
\left( \begin{array}{cc} 1 & 0 \\ 0 & \xi
\end{array} \right)
$ with $\xi \neq 1$. We must have

$$
\left( \begin{array}{cc} a & 0 \\ b & c
\end{array} \right)
\left( \begin{array}{cc} 1 & 0 \\ 0 & \bar{\xi}
\end{array} \right)
\left( \begin{array}{cc} a & 0 \\ b & c
\end{array} \right)^{-1} =
\left( \begin{array}{cc} 1 & 0 \\ 0 & \theta
\end{array} \right)
$$ whence $b = 0$.
\hfill Q.E.D. \\\\

Now we would like to see when the lift $\tilde{\sigma}$ is an involution. Let
us denote by $\phi: A'' \times F \rightarrow S$ the map in the diagram.
\begin{REM}
\label{liftnonempty} Let $\tilde{\sigma}$ be a lift of $\sigma$, we have two
different cases:
\begin{enumerate}
\item $\exists z \in S$ such that $\sigma(z) =z$. Then, since the covering is
Galois, for all $z' \in \phi^{-1}(z)$, there exists a lift $\tilde{\sigma}$ of
$\sigma$ such that $\tilde{\sigma}(z')= z'$. But then $\tilde{\sigma}^2$ lifts
the identity and has a fixed point, therefore it is the identity. \\
\item $Fix(S) = \emptyset$.\\  Let $z,w \in S$, $\forall z' \in \phi^{-1}(z)$,
$\forall w' \in \phi^{-1}(w)$, $\exists ! $ lift $\tilde{\sigma}$, with
$\tilde{\sigma}(z') = w'$.  \\
\end{enumerate}
\end{REM}

\section{Orbifold fundamental groups and the topological type of a real
hyperelliptic surface}

Let $(X, \sigma)$ be a smooth real variety (i.e. a smooth complex variety
together with an antiholomorphic involution $\sigma$) of complex dimension
$n$. Then we have a double covering $\pi : X \rightarrow Y = X/<\sigma>$
ramified on $X' = X({\R}) = Fix(\sigma)$. $X'$ is a real submanifold of
codimension $n$, hence the map
$$ \pi_1(X - X') \rightarrow \pi_1(X)$$ is surjective if $n \geq 2$ and it is
an isomorphism for $n \geq 3$. The singularities of $Y$ are contained in $Y'
= \pi (X')$ and there we have a local model ${\R}^n \times ({\R}^n/(-1))$.
Therefore
$Y$ is smooth for $n = 2$ and topologically singular for $n \geq 3$.  The
local punctured fundamental group $\pi_1(Y - Y')_{loc}$ is isomorphic to
${\Z}$ for $n = 2$, while it is isomorphic to ${\Z}/2$ for $n \geq 3$. This
means that the kernel of the surjection
$$\pi_1(Y - Y') \rightarrow \pi_1(Y)$$ is normally generated by loops
$\gamma$ around the components of $Y'$. If
$n \geq 3$, then we have automatically $\gamma^2 = 1$.

\begin{DEF} Assume $dim_{{\Co}} X = n \geq 2$. Then the orbifold fundamental
group of
$Y$ (or of $(X, \sigma)$) is defined to be $\pi_1(Y - Y')/<\gamma^2>$.
\end{DEF}

\begin{REM}
\begin{enumerate}
\item If $n \geq 3$, the orbifold fundamental group of $Y$ coincides, as
observed, with the fundamental group of $Y-Y'$.
\item If $X' = \emptyset$, $\pi_1 ^{orb}(Y) = \pi_1(Y)$ and we have an exact
sequence

\begin{equation}
\label{orb} 1 \rightarrow \pi_1(X) \rightarrow \pi_1^{orb}(Y) \rightarrow
{\Z}/2
\rightarrow 1
\end{equation}

\item For $n \geq 2$ we have the exact sequence \ref{orb}, since we have

$$\diagram   1 \rto & \pi_1(X - X')  \dto \rto   & \pi_1(Y - Y')  \dto \rto &
{\Z}/2 \dto  \rto & 1  \\ 1 \rto & \pi_1(X) \rto & \pi_1^{orb}(Y) \rto &
{\Z}/2
\rto & 1
\enddiagram$$

and the kernel of the map $\pi_1(X - X') \rightarrow \pi_1(X)$ is normally
generated by the above $\gamma^2$'s.
\end{enumerate}
\end{REM}

Let $\tilde{X}$ be the universal covering of $X$, so that $X =
\tilde{X}/\pi_1(X)$. The exact sequence \ref{orb} defines a group which is the
group of liftings of $\sigma$ to $\tilde{X}$, so that $Y =
\tilde{X}/\pi_1^{orb}(Y)$. \\  If $X' = X({\R}) \neq \emptyset$, let us
choose a base point $x_0 \in X'$. Then $\sigma$ acts on $\pi_1(X, x_0)$ by an
automorphism of order 2, which is obtained by conjugation, since we have the
following

\begin{LEM} If $X' \neq \emptyset$, then the exact sequence \ref{orb} always
splits.
\end{LEM} {\bf Proof.} The proof follows immediately from corollary
\ref{fixlift}.
\hfill Q.E.D. \\

\begin{REM} The topological invariants of $(X, \sigma)$ are: the topological
invariants of
$X({\R})$, the topological invariants of $Y$ and $\pi_1^{orb}(Y)$.
\end{REM}

We have the following

\begin{TEO}
\label{orbifold}
  Let $(S, \sigma)$ be a real hyperelliptic surface. Then the topological
  and also the differentiable type
of the pair $(S, \sigma)$ is completely determined by the orbifold fundamental
group exact sequence.
\end{TEO}

{\bf Proof.}  We want first of all to show how $\Pi := \pi_1(S)$
determines the topological (actually differentiable) type of the hyperelliptic
surface $S$. \\  Consider the exact homotopy sequence associated with the
covering $\psi : S'
\rightarrow S$ described in the proof of theorem \ref{class}.

\begin{equation}
\label{pig'} 1 \rightarrow \Omega'  = \pi_1(S') \cong {\Z}^4 \rightarrow
\Pi \rightarrow G' \rightarrow 1
\end{equation} where $G' = G/_{translations}$. The exact sequence
\begin{equation}
\label{albanese} 1 \rightarrow \Gamma = [\Pi, \Pi] \rightarrow \Pi
\rightarrow \Lambda =
\Pi/[\Pi,\Pi] \rightarrow 1
\end{equation} is given by the induced map of fundamental groups associated to
the Albanese map and $\Gamma \cong {\Z}^2$ is the fundamental group of the
fibre of
$\alpha$.\\
$\Omega'$ is then the kernel of the action of $\Pi$ on $\Gamma$ by conjugation
($\Gamma \subset \Omega'$ and $\Omega'/\Gamma = \Lambda'$ in our previous
notation).\\ Now $\Omega'$ is a representation of $G'$ and $\Gamma$ is a
subrepresentation. An easy calculation shows that $(\Omega')^{G'} =
\Lambda ''$, which yields a direct sum $\Omega = \Lambda'' \oplus \Gamma$,
such that
$\Pi/ \Omega = G$. Notice that since $\Omega'$ is abelian and by \ref{pig'},
we know that $(\Omega')^{G'}$ is the centre of $\Pi$, so $\Lambda''$ is the
centre of $\Pi$. \\   Since $\Omega \subset \Omega'$, the universal cover
$\tilde{X} \cong {\R}^4$ of $S$ is homeomorphic to $\Omega' \otimes {\R}$ on
which $\Omega'$ acts freely by translations. We have $S \cong (\Omega'
\otimes {\R})/ \Pi$, thus it suffices to show that the exact sequence
\ref{pig'} determines the action of $\Pi$ on the universal covering $\Omega'
\otimes {\R}$. \\   The action is given by a group of affine transformations
of
${\R}^4$ and since $\Omega'$ is the sugroup of translations, the action of
$G'$ on the torus $\Omega' \otimes {\R}/\Omega'$ has a linear part which is
determined by the conjugation action of $G'$ on $\Omega'$. \\  Since $G'$ is
cyclic, we may find a lift $\lambda (g)$ to $\Pi$ of its generator $g$. Now,
for all $g'
\in G'$ we have a lift $\lambda (g')$ in $\Pi$ and if $m = |G'|$, we have
$\lambda (g)^m
\in \Omega'$. Every element $\gamma
\in \Pi$ can be uniquely written as $\omega' \lambda (g')$, with $\omega' \in
\Omega'$, $g' \in G'$. \\ Since $\lambda (g)^m \in \Omega'$ and it is
invariant by conjugation by
$\lambda(g)$, we have $\lambda (g)^m \in (\Omega')^{G'} = \Lambda''$,
therefore we let $g$ act on $(\Lambda'' \otimes {\R}) \oplus (\Gamma
\otimes {\R})$ by $$(e,f) \mapsto (e + \frac{1}{m} \lambda (g)^m, g(f))$$
where
$g$ acts on $\Gamma$ by conjugation. \\  It is immediate that the above
action is precisely the one yielding $S$ as $(E
\times F )/G = (\Lambda'' \otimes {\R} \oplus \Gamma \otimes {\R})/\Pi$ and
that the way we described it is completely dictated by $\Pi$ as an abstract
group. \\ Thus we have proven that $\Pi$ determines the topological type of
$S$. \\ Let us now consider the ramified covering $E \times F
\rightarrow E \times F/\hat{G}$, where
$\hat{G}$ is the extended symmetry group defined in \ref{hatg}. Then we have
the following exact sequences

$$
\begin{array}{cccccccccc} & 1 & \rightarrow &   \Omega   &\rightarrow&
\Pi      & \rightarrow & G &\rightarrow&1\\ &   &             & \Arrowvert
&           &     \cap     &             & \cap &&\\ & 1 & \rightarrow &
\Omega   &\rightarrow&   \hat{\Pi}  & \rightarrow &\hat{G}&\rightarrow&1
\end{array}
$$

where $\hat{\Pi}$ is the orbifold fundamental group,

$$
\begin{array}{cccccccccc} & 1 & \rightarrow &   \Omega'   &\rightarrow&
\Pi      & \rightarrow & G' &\rightarrow&1\\ &   &             & \Arrowvert
&           &     \cap     &             & \cap &&\\ & 1 & \rightarrow &
\Omega'   &\rightarrow&   \hat{\Pi}  & \rightarrow &\hat{G}'&\rightarrow&1
\end{array}
$$

We want to describe how $\hat{\Pi}$ acts on $\Omega' \otimes {\R}$. In order
to understand the action of $\hat{\Pi}$ it suffices to describe the action of
a suitable element $\sigma$ in $\hat{\Pi} - \Pi$. $\sigma$ acts on
$\Omega' \otimes {\R}$ by an affine transformation: $v \mapsto Av + b$, of
which we know the linear part $A$, that is determined by the conjugation
action of
$\sigma$ on $\Omega'$. \\ On the other hand, we have $\sigma^2 \in
\Pi$, therefore we know the affine map $\sigma^2$. Since we have
$$\sigma^2(v) = A^2v + Ab + b =: A^2 v + b',$$ we are able to determine $b$
uniquely from $A$ and $b'$, in the case in which
$(A + I)$ is invertible. \\ We also know that $\sigma$ is of product type, and
that $A^2 = A_1^2 \oplus A_2^2$ has the property that $A_1^2 = I$. Thus
$W = {\R}^2_1 = \Lambda'' \otimes {\R}$ splits as a direct sum of eigenspaces
$W^+ \oplus W^-$ and we can recover the translation vector
$b_1^+$ by what we have remarked above. Whereas on $W^- = Im(A_1 - I)$ we can
change coordinates by a translation in such a way that the action of
$\sigma$ is linear on $W^-$.\\ Therefore we can restrict our attention to the
second component $A_2$. If
$(-1)$ is not an eigenvalue of $A_2$, we are done. Since $\sigma$ is
antiholomorphic, $(-1)$ is an eigenvalue of $A_2$ if and only if $A_2^2 = I$.
Therefore it remains to treat the case in which $A_2^2 =I$ and $(A_2 g)^2 =
I$,
$\forall g \in G$. \\ In this case, the second component of $\sigma$ is given
by
$\sigma_2(x) = A_2 x + c$. Let us consider another  element in $\hat{\Pi} -
\Pi$, $s = g \circ \sigma$, where $g \in G$ is not a translation.  We know the
action of $\sigma \circ s$ and of $s \circ \sigma$, since they are in $\Pi$.

$$\sigma_2 \circ s_2(x) = A_2(gA_2 x + c') + c,$$
$$s_2 \circ \sigma_2 (x) = gA_2(A_2 x + c) + c' = gx + gA_2 c + c'.$$ Looking
at the translation terms, we get knowledge of
$gA_2 c + c'$ and $A_2 c' + c$.\\ We argue by looking at the rank of the
matrix
$
\left( \begin{array}{cc} gA_2 & I \\ I & A_2
\end{array} \right)
$ which is applied to the vector
$
\left( \begin{array}{c} c \\ c'
\end{array} \right).
$

It has the same rank as

$$rank
\left( \begin{array}{cc} g & I \\ A_2 & A_2
\end{array} \right)
  = rank
\left( \begin{array}{cc} g-I & I \\ 0 & A_2
\end{array} \right)
  = 2 + rk(g -I).$$ Then we are done, since $(g - I)$ is invertible, because
$g$ is holomorphic and it is not a translation, thus it does not have $1$ as
eigenvalue.   So the proof of the theorem is concluded, since we have shown
that the exact sequence of the orbifold fundamental group completely
determines the action of $\hat{\Pi}$ on $\Omega \otimes {\R}$.
\hfill Q.E.D. \\

\begin{COR}
\label{toptype}  The topological type of a real hyperelliptic surface $(S,
\sigma)$ determines:
\begin{enumerate}
\item the Bagnera - de Franchis group $G$,
\item the extension $1 \rightarrow G \rightarrow \hat{G} \rightarrow {\Z}/2
\rightarrow 1$
\item the topological type of the action of ($\hat{G} \supset G$) on $E$ and
$F$, or equivalently the affine equivalence class of the representation of
$\hat{\Pi}/\Gamma
\supset \Pi/\Gamma$ on $\Lambda'' \otimes {\Q}$ and of $\hat{\Pi}/\Lambda''
\supset \Pi/\Lambda''$ on $\Gamma \otimes {\Q}$. In particular it determines
the topological type of the real elliptic curve $(E/G,
\hat{G}/G)$.
\end{enumerate}
\end{COR}
\begin{REM}
\label{repre} Notice that the representation $\rho$ of $\hat{\Pi}$ in the
group of the affine tranformations of $\Omega \otimes {\Q}$ takes values in
$A(2, {\Z}, {\Q}) \times A(2, {\Z}, {\Q})$, where
$A(2, {\Z}, {\Q}) = \{v \mapsto Bv + \beta, \ | B \in GL(2, {\Z}), \
\beta \in {\Q}^2\}$.
\end{REM}

We want now to discuss in some greater generality
than needed for our present purposes the notion of moduli space of real
varieties. Since we will only assume that $X$ is a complex manifold,
we have to adopt the point of view of Kodaira-Spencer-Kuranishi
(cf. \cite{ca3}, \cite{k-m}, \cite{ko}).

For a general complex manifold $X$ we have the Kuranishi family of
deformations of $X$, $ \phi : \mathcal X \rightarrow \mathcal B  $. Here,
the base $\mathcal B  $ of the Kuranishi family is a complex  analytic subset
of the vector space $ H^1 (X, \Theta_X ) $ corresponding to the complex
structures $J$ which satisfy Kuranishi's integrability equation.

If now $(X, \sigma)$ is real, we want to see when a neighbouring complex
structure $J$ is such that the differentiable map $\sigma$  remains
antiholomorphic (whence we get a deformation $X_t$ of the complex manifold
$X$ such that the new pair $( X_t , \sigma)$ is still real).

The corresponding equation writes down simply as

$\sigma_* J = - J \sigma_* $, or equivalently,

$ - \sigma_* J \sigma_*  =  J  $.

We see immediately that this condition means that $J$ lies in the fixed locus
of the  involution induced by $\sigma$.

Interpreting now $J$ as a harmonic representative of a Dolbeault cohomology
class $\theta \in  H^1 (X, \Theta_X ) $, we are going to show more precisely
that $\sigma$ induces a complex antilinear  involution $\sigma^*$ on the
vector space
$ H^1 (X, \Theta_X ) $, such that  the real part ${\cal B }( \R )$ consists of
the deformations of $X$ for which $\sigma$ remains real.

Recall that the complex structure $J_0$ of $X$ induces a splitting of
the complexified real tangent bundle of $X$, $T_X \otimes {\Co} = T^{0,1}
\oplus T^{1,0} $ where $T^{0,1} = \bar T^{1,0} $.

The form $\theta = \Sigma_{i \bar j} \theta_{i \bar j} ( \partial / \partial
z_i \otimes d \bar y_j )$ can be interpreted as yielding a linear map
$ \theta' : T^{0,1} \rightarrow T^{1,0} $ such that the new subbundle
$T^{1,0}_{\theta} = \{ (u,v) \in T^{0,1}
\oplus T^{1,0} \ | \ v = \theta' (u) \} $.

$\theta'$ is completely determined by the matrix $ \phi = \theta_{i \bar j}$,
and then the  subspace $T^{0,1}_{\theta} = \{ (u,v) \in T^{0,1}
\oplus T^{1,0} \ | \ u = \bar \theta' (v) \} $, where $ \bar \theta' $ is
determined by the matrix $ \bar \phi = \bar \theta_{\bar i  j}$ in the
chosen coordinate basis.

Saying that $\sigma$ is antiholomorphic amounts to saying that  its
differential induces complex linear isomorphisms $\sigma_* :  T^{0,1}
\rightarrow T^{1,0} $ and  $\sigma_* :  T^{1,0}
\rightarrow T^{0,1} $ which are inverses and conjugate to each other.

Whence, if $A$ is the matrix of $\sigma_* :  T^{0,1}
\rightarrow T^{1,0} $, $A^{-1} = \bar A $ is the matrix of
$\sigma_* :  T^{1,0} \rightarrow T^{0,1} $.

We want now to write down the condition that $\sigma$ be antiholomorphic
for the new complex structure induced by the form $\theta$.

Again, this means that $\sigma_* :  T^{1,0}_{\theta}
\rightarrow    T^{0,1}_{\theta} $, i.e., that for each vector $u$, the image
$\sigma_* (u , \theta' (u)) = ( \bar A \phi (u) , A (u))$
satisfies the equation of $T^{0,1}_{\theta} $, namely, we have

$ \bar A \phi (u) = \bar \phi A (u)$.

Again we can equivalently  define

$ \sigma^* (\phi) : =  A \bar \phi A$

write the above (since it must hold for each vector $u$) as

$ \phi = \sigma^* (\phi) $.

It is now obvious that $ \sigma^*$ is complex antilinear, and it is an
involution since

$ {\sigma^*}^2 (\phi) = A \bar A  \phi \bar A A = \phi$.

By what we have shown above, it follows right away that ${\cal B }( \R )$,
being the intersection of a complex analytic space with the fixed part of a
complex antilinear involution, is a real analytic space.

In general, by a result of Wavrik (cf. \cite{wav}),
if the function $ h : {\cal B } \rightarrow  \Z $ given
by

$ h(t) = h^0 (X_t, \Theta_{X_t} ) $

is constant on the germ ${\cal B }$, then the quotient

  ${\cal B }/ Aut (X)$ is a local moduli space.

  The group $Aut(X)$ does not however act on the set ${\cal B }( \R )$,
  therefore we take the smaller subgroup

  $ Aut^{\sigma}(X)$ : = $ \{  \phi \in Aut(X) | \phi^* (\sigma) = \sigma \}$,
  where $\phi^* (\sigma) := \phi \sigma \phi^{-1}  $.

  We finally can define the real local moduli space as

  \begin{DEF}
  Under the assumption that the function $ h^0 (X_t, \Theta_{X_t} ) $
  is constant on the base of the Kuramishi family, the real local moduli
  space of $(X,\sigma)$is defined as the quotient ${\cal B }( \R ) /
Aut^{\sigma}(X) $.

  \end{DEF}

\begin{REM} Therefore, the real local moduli space  is just a real semi
algebraic space, and it maps neither  surjectively  nor one to one to the
real part of the complex local moduli space.

\end{REM}

\begin{REM} Once we have  local moduli spaces for the varieties or manifolds
under consideration, the standard procedure is to consider the global
moduli space as the set of isomorphism classes of such varieties and to use
the local moduli spaces as giving local charts. For instance, if we have local
moduli spaces for a certain class of complex manifolds, these charts provided
by the local moduli spaces yield for our global moduli space the structure of
a complex analytic space, possibly non Hausdorff.

Likewise, for real algebraic varieties such that Wavrik's condition
holds, the global moduli space is a semianalytic space.

\end{REM}
Now, one can give different definitions for polarized algebraic varieties, and
again one has to distinguish between the real part of the quotient of a
Hilbert scheme, and the quotient of its real part. However, the approach via
the Kuranishi family is particularly suitable for the case of hyperelliptic
surfaces, as we are going now to see.

\begin{TEO} The moduli space of real hyperelliptic surfaces of a given
topological type is irreducible.
\end{TEO}

{\bf Proof.} Let us fix an orbifold fundamental group
$\hat{\Pi}$, then theorem \ref{orbifold} and corollary \ref{toptype} tell us
that $\hat{\Pi}$ determine the representation
$\rho: \hat{\Pi} \rightarrow {A(2, \Z, \Q)} ^2$. Such a representation
induces a representation of $\hat{G}'$ into $GL(2, \Z)^2$:
$$\hat{G}' \stackrel{(\lambda_1, \lambda_2)} \longrightarrow GL(2, \Z)^2,$$
where $\hat{G}'$ is the group fitting into the exact sequence
$$ (\#) \ 1 \rightarrow G' \rightarrow \hat{G}' \rightarrow \Z/2
\rightarrow 1$$

We observe that giving a hyperelliptic surface with the given topological type
amounts to giving two complex structures $J_1$, $J_2$ on $\Z^2 \otimes \Q$
such that $\lambda_i(G')$ consists of $\Co$ - linear maps, whereas
$\lambda_i(\hat{G}' - G')$ consists of $\Co$ - antilinear maps ($i =1,2$).

Therefore we have to solve the following problem:

given a representation $\lambda: \hat{G}' \rightarrow GL(2, \Z)$, find all
the complex structures $J$ such that for a generator $g'$ of $G'$  we have

$\lambda(g') J = J \lambda(g') $,

and for a $\sigma \not \in G'$ we have

$\lambda(\sigma) J = - J \lambda(\sigma).$

We have of course to keep in mind that, if $M$ is a fixed differentiable
manifold with an involution $\sigma :  M \rightarrow  M$, and we look for the
complex structures $J$ which make $\sigma$ an antiholomorphic involution, for
each solution $J$ we shall also find the solution $-J$.

But, as we already remarked in the introduction, $\sigma$ provides an
isomorphism between the complex manifolds $ (M, J)$ and $ (M,-J)$, and
clearly it conjugates the involution $\sigma$ to itself.

Therefore, if we shall see that the parameter space for our complex
structures $J$ will consist of exactly two irreducible components, exchanged
by the involution $ J \to -J$, it will follow that the moduli space is
irreducible.

We shall prove in the next section that the extension $(\#)$ always splits
(see corollary \ref{gcyc}), therefore we may assume that $\sigma^2 = 1$.  In
particular we can find a basis in such a way that either
$\zeta:= \lambda(\sigma) =
\left( \begin{array}{cc} 1 & 0\\ 0 & -1
\end{array} \right)
$
  or $\zeta =
\left( \begin{array}{cc} 0 & 1\\ 1 & 0
\end{array} \right)
$.

If we now set $B := \lambda(g')$, since $G' \cong \Z/q$ for $q = 2, 3, 4, 6$,
we have either $\zeta B \zeta = B$, or $\zeta B \zeta = B^{-1}$. We have two
different cases:

Case I: $\zeta:= \lambda(\sigma) =
\left( \begin{array}{cc} 1 & 0\\ 0 & -1
\end{array} \right)
$.

Then, if we set $B = \lambda(g') =
\left( \begin{array}{cc} b_{11} & b_{12}\\ b_{21} & b_{22}
\end{array} \right)
$, we have again two cases:
\begin{enumerate}
\item $\zeta B \zeta = B$,
\item $\zeta B \zeta = B^{-1}$.
\end{enumerate} In case 1 we find
$
\left( \begin{array}{cc} b_{11} & b_{12}\\ b_{21} & b_{22}
\end{array} \right) = \left( \begin{array}{cc} b_{11} & -b_{12}\\ -b_{21} &
b_{22}
\end{array} \right)$, whence, since $det(B) = 1$, we have $B = \pm Id$. But we
shall see in \ref{action} that this case occurs if and only if $q = 2$ and $B
= -Id$.

In case 2 we find $
\left( \begin{array}{cc} b_{11} & -b_{12}\\ -b_{21} & b_{22}
\end{array} \right) = \left( \begin{array}{cc} b_{22} & -b_{12}\\ -b_{21} &
b_{11}
\end{array} \right)$, whence $b_{11} = b_{22} =: b$, $b^2 - b_{12}b_{21} = 1$,
the characteristic polynomial is $\lambda^2 - 2b \lambda +1$. But since
either $B = -Id$, or $B$ has no real eigenvalues, $b^2 -1 <0$, thus $b = 0$,
whence $B =
\pm
\left( \begin{array}{cc} 0 & -1\\ 1 & 0
\end{array} \right)
$ and we may assume $B =
\left( \begin{array}{cc} 0 & -1\\ 1 & 0
\end{array} \right)
$ by changing generator of $G'$.

Case II: $\zeta = \left( \begin{array}{cc} 0 & 1\\ 1 & 0
\end{array} \right)
$.  We have again the two cases $1$ and $2$.

Case II 1: We find $
\left( \begin{array}{cc} b_{11} & b_{12}\\ b_{21} & b_{22}
\end{array} \right)  = \left( \begin{array}{cc} b_{22} & b_{21}\\ b_{12} &
b_{11}
\end{array} \right) $, thus $B =\left( \begin{array}{cc} b & c\\ c & b
\end{array} \right)$ and since $det(B) = b^2 - c^2 =  1$, we have $b+c =
\pm 1$,
$b -c = \pm 1$, whence $c = 0$ and $B = -Id$.

Case II 2: We have $\left( \begin{array}{cc} b_{22} & b_{21}\\ b_{12} & b_{11}
\end{array} \right) = \left( \begin{array}{cc} b_{22} &- b_{12}\\ -b_{21} &
b_{11}
\end{array} \right)$ which yields $B = \left( \begin{array}{cc} b_{11} & c\\
-c & b_{22}
\end{array} \right)$. But $b_{11}b_{22} + c^2 = 1$, thus either $c =0$ and $B
= -Id$, or $b_{11}b_{22} \leq 0$.

If $c = \pm 1$, $b_{11} b_{22} = 0$. By exchanging the two basis vectors $e_1$
and $e_2$ we may assume $c = -1$ and by replacing $B$ with $B^{-1}$ we can
assume
$b_{22} = 0$.

Thus $B = \left( \begin{array}{cc} b & -1\\ 1 & 0
\end{array} \right)$ and either $b = \pm 1$ if $G' = \Z/6$ or $G' = \Z/3$, or
$b = 0$, $B = \left( \begin{array}{cc} 0 & -1\\ 1 & 0
\end{array} \right)$ if $G' = \Z/4$.

Assume then $|c| \geq 2$. Since $|b_{11} + b_{22}| <2$, we have two possible
cases:
$b_{11} + b_{22} =0$, $b_{11} + b_{22} =\pm 1$.

In the first case we find $B =\left( \begin{array}{cc} b & c\\ -c & -b
\end{array} \right)$ and $-b^2 + c^2 = 1$, thus $b = 0$, $c = \pm 1$, a
contradiction.

In the second case the matrix $B$ has the form
$\left( \begin{array}{cc} b & c\\ -c & \pm1 -b
\end{array} \right)$, we may assume $b >0$, and we have $\pm b - b^2 + c^2 =
1$.

The equation $\pm b - b^2 + c^2 = 1$ is equivalent to $c^2 = 1 + b^2 \mp b =
(b
\mp 1)^2 \pm b$, whence $\pm b = (c -b \pm 1)(c + b \mp 1)$. Taking absolute
values in the last equation, we see that we cannot have $|c|\geq 2$.

We have therefore obtained the following possibilities for $\zeta$ and $B$:

\begin{itemize}
\item Case I 1:
$\zeta = \left( \begin{array}{cc} 1 & 0\\ 0 & -1
\end{array} \right)$,
$B = -Id$.
\item  Case I 2:
$\zeta = \left( \begin{array}{cc} 1 & 0\\ 0 & -1
\end{array} \right)$, $B = -Id$, or $B = \left( \begin{array}{cc} 0 & -1\\ 1
& 0
\end{array} \right)$.
\item Case II 1:
$\zeta = \left( \begin{array}{cc} 0 & 1\\ 1 & 0
\end{array} \right)$,
$B = -Id$.

\item Case II 2:
$\zeta = \left( \begin{array}{cc} 0 & 1\\ 1 & 0
\end{array} \right)$,
$B = -Id$, or $B = \left( \begin{array}{cc} 1 & -1\\ 1 & 0
\end{array} \right)$, or $B = \left( \begin{array}{cc} -1 & -1\\ 1 & 0
\end{array} \right)$, or $B = \left( \begin{array}{cc} 0 & -1\\ 1 & 0
\end{array} \right)$.
\end{itemize}

We now look for a complex structure
$J = \left( \begin{array}{cc} d_{11} & d_{12}\\ d_{21} & d_{22}
\end{array} \right)$ such that

$$\zeta J \zeta = -J,$$
$$BJ = JB.$$

Case I: the first equation reads $ \left( \begin{array}{cc} d_{11} & -d_{12}\\
-d_{21} & d_{22}
\end{array} \right) = -J$, thus $d_{11} = d_{22}=0$, and since the
characteristic polynomial for $J$ is $\lambda^2 +1$, we have $det(J) = 1$,
equivalently
$d_{12} d_{21} = -1$, whence $J$ has the form
$ \left( \begin{array}{cc} 0 & -d\\ 1/d & 0
\end{array} \right)$.

Whence, as promised, the parameter space consists of the two branches of an
hyperbola, which are exchanged upon multiplication by $-1$.

Case II: the equation $\zeta J \zeta = -J$ reads
$ \left( \begin{array}{cc} d_{22} & d_{21}\\ d_{12} & d_{11}
\end{array} \right) = -J$, thus we obtain $J = \left( \begin{array}{cc} a &
b\\ -b & -a
\end{array} \right)$, with $a^2 - b^2 = -1$, whence we can write $a = (c -
1/c)/2$, $b = (c + 1/c)/2$ and the conclusion is exactly as before.

Let us now consider the commutation relation $BJ = JB$.

If $B = -Id$, it is obviously verified.

If $B = \left( \begin{array}{cc} 0 & -1\\ 1 & 0
\end{array} \right)$, we have cases I and II:

In case I the commutation relation yields
$ \left( \begin{array}{cc} -d & 0\\ 0 & -1/d
\end{array} \right) = \left( \begin{array}{cc} -1/d & 0\\ 0 & -d
\end{array} \right)$, whence $d^2 = 1$, i.e. $d = \pm 1$,
$J = \pm \left( \begin{array}{cc} 0 & 1\\ -1 & 0
\end{array} \right)$.

In case II the commutation relation yields
$ \left( \begin{array}{cc} b & -a\\ -a & b
\end{array} \right) = \left( \begin{array}{cc} b & a\\ a & b
\end{array} \right)$, whence $a =0$ and $b^2 = 1$, i.e. $b = \pm 1$, $J =\pm
\left( \begin{array}{cc} 0 & 1\\ -1 & 0
\end{array} \right)$.

Assume now $B = \left( \begin{array}{cc} 1 & -1\\ 1 & 0
\end{array} \right)$, then we are in case II and the commutation relation
reads
$\left( \begin{array}{cc} a+b & a+b\\ a & b
\end{array} \right) = \left( \begin{array}{cc} a+b & -a\\ -a-b & b
\end{array} \right)$, whence $b = -2a$ and since $b^2 - a^2 = 1$, we have $a =
\pm 1/\sqrt{3}$.

We observe that if $J$ commutes with $B$, it commutes also with $B^2$, thus,
up to a base change, we have solved also the case $G' = \Z/3$.

\hfill Q.E.D. \\

\section{The extended symmetry group}

In this section we want to describe the extended
  Bagnera - de Franchis group, fitting into the exact sequence:

$$(*) \ \ 0 \rightarrow G \stackrel{i} \rightarrow \hat{G}
\stackrel{\pi} \rightarrow {\Z}/2  \rightarrow 1$$.

Since $G$ is a normal Abelian subgroup of $\hat{G}$, conjugation induces an
action of $\hat{G}/G = {\Z}/2$ on $G$.

Recall that, once such action is specified, the equivalence classes of such
extensions are in bijective correspondence with the elements of
$$H^2({\Z}/2, G) = \frac{ker(\sigma - 1)}{(1 + \sigma)G} =
\frac{G^{{\Z}/2}}{(1 +
\sigma)G},$$ where ${\Z}/2 = <\sigma>$. \\

Let us consider the antiholomorphic elements of the extended Bagnera- de
Franchis group, which will be denoted by $\tilde{\sigma}$ (they are the lifts
of
$\sigma$ in $\hat{G}$). We showed in proposition 2.6 that their action is
also of product type, whence we can restrict our preliminary investigation to
the question: which extended Bagnera- de Franchis groups act as group of
diholomorphic automorphisms of  an elliptic curve?

\begin{DEF} A diholomorphic action of an extended Bagnera- de Franchis group
$\hat{G}$ on an elliptic curve is an action such that
$G$ is precisely  the subgroup of holomorphic automorphisms in $\hat{G}$,
while the elements in $\hat{G} - G$ act as antiholomorphic automorphisms.
\end{DEF}

\begin{LEM}
\label{quatra} Consider a diholomorphic action of an extended Bagnera- de
Franchis group $\hat{G}$ on an elliptic curve.  Then the square of an
antiholomorphic map $\tilde{\sigma}$ in
$\hat{G}$ is a translation.
\end{LEM}  {\bf Proof.} Let $\tilde{\sigma}^2 = g \in G$. Passing to the
universal cover we can write
$\tilde{\sigma}(z) = a \bar{z} + b$. We have then
$ \tilde{\sigma}^2(z) = a \bar{a} z + a \bar{b} + b = g(z)$. We know that
$g(z) = \xi z + c$ with  $\xi^n = 1$, $n = |G|$. \\  Since $\xi = |a|^2 \in
{\R}_{>0}$, we obtain $\xi =1$.

\hfill Q.E.D. \\

\begin{COR}
\label{gcyc} If the group $G$ is cyclic, then the extension $(*)$ always
splits.
\end{COR} {\bf Proof.} Since $G = G'$, for any lift $\tilde{\sigma}$ of
$\sigma$, we must have
$\tilde{\sigma}^2 = Id \in G$, so $(*)$ splits.
\hfill Q.E.D. \\

  We want now to recall the description of the antiholomorphic
  maps ${\sigma}$ acting
  on an elliptic curve $C$, whose square ${\sigma}^2$ is a translation.

  Write $C = {\Co}/\Gamma$,
  $\sigma(z) = a \bar{z} + b$, $\sigma^2(z) = a
\bar{a} z + a \bar{b} + b$.

Whence
\begin{equation}
\label{A} |a|^2 =1,
\end{equation}

The following condition must be verified:
$$a \bar{\Gamma} = \Gamma$$ which is clearly equivalent to the existence of
integers $m,n,m',n'$ such that:

\begin{equation}
\label{B} a = m + n \tau \in \Gamma,
\end{equation}
\begin{equation}
\label{C} a \bar{\tau} =  m' + n' \tau, \ m',n' \in {\Z}, \ with \  m n' - nm'
= - 1.
   \end{equation}

We may rewrite \ref{A} as

\begin{equation}
\label{A'} ( m + n \tau )( m + n \bar{\tau} ) = m^2 + n^2 |\tau|^2 + mn (\tau
+\bar{\tau}) = 1
\end{equation}

We may assume that $\tau$ lies in the modular triangle, i.e., that

$ | \tau| \geq 1 , | Re \tau| \leq 1/2$.

Then $ Im(\tau) \geq  \sqrt{3}/ 2$, and since $a = m + n \tau$,
$ Im(a) = n \ Im(\tau)  \leq 1$, and we conclude that $|n| \leq 1$.

We have then the following cases:

\begin{itemize}
\item
$n=0$: then, since $|a|=1$, $m =a = - n'=  \pm 1$.

Moreover, \ref{C} tells us that $ 2(Re) \ \tau = - n' m'$, whence we either
have
$ Re(\tau) = 0$ or  $ Re(\tau) = - 1/2$.
\item observing that if $|n|=1$, by \ref{A'} we infer

$ m^2 + 1 \leq m^2 + |\tau|^2= 1 \pm (\tau +\bar{\tau}) \leq 1 + |m|$,

thus $|m| \leq 1$, $|\tau| = 1$, giving rise to the following two cases:
\item
$|n|=1, m = 0 , |\tau| = 1$ (thus $ a = \pm \tau$)

\item
$|n|=1, |m|=1$, $|\tau| = 1$ and since $ 1 + |\tau|^2= 1 - mn (\tau
+\bar{\tau})$,
  we may also assume $ Re (\tau) = - 1/2$, thus $m=n$, $ a= m (1+\tau)$.
\end{itemize}

\begin{LEM}
\label{tra} Assume now that $\sigma$ is an antiholomorphism of an elliptic
curve whose square is a translation of finite order $d$.
  Then, we may choose the origin in the universal cover $\Co$ in such a way
that
$\sigma(z) = a \bar{z} + b$, with
$a$ as above,
$ b \in  (1/d) \Z$ for $ a \neq \pm 1$,
$b \in (1/2d) \Z$ if $a =1$,
$ b \in  (1/2d) \Z (\tau )$ for $a = -1$, $ Re \ \tau =0$,
$ b \in  (1/2d) \Z ( 2 \tau + 1)$ for $a=-1,  Re \ \tau =-1/2$.

Assume further that $d = 1$, i.e. that $\sigma$ is an involution. Then,
obviously, we may get $b = 0$ if and only if $Fix(\sigma) \neq
\emptyset$.

  $Fix(\sigma) = \emptyset$ if and only if  $Re(\tau) = 0$ and we may choose
the
  origin in such a way that $ b = 1/2$ for $a=1$, or $ b =\tau/2$
  for $a=-1$.

Moreover $\sigma$ normalizes a finite group of translations $T$, if and only
if, identifying
$T$ with a subgroup of $C$, $ a \bar{T} = T $.

\end{LEM}

{\bf Proof.} We look first for a vector $w$ such that $ a
\bar{w} - w + b : =
\beta$ be either  a real vector, or an imaginary vector.

We look therefore at the image of the linear map $ w \rightarrow  a
\bar{w} - w $. Its complexification $ \Co \otimes \Co \rightarrow
\Co \otimes \Co$ has a matrix

$
\left( \begin{array}{cc} -1 & a\\
\bar{a} & -1
\end{array} \right)
$   whose determinant is zero, whence, the image of the above linear map
equals $
\R (a-1)$ for $a \neq 1$, and
  $ \R  i \ (Im \ \tau)$ for $ a=1$.

We can therefore achieve that $\beta$ be a real vector, unless $a$ is real,
$ a \neq 1$, i.e., $ a=-1$, in which case we can achieve that
$\beta$ be an imaginary vector.

Since $\sigma^2(z) = z + a \bar{\beta} + \beta$, if $\beta$ is real we get $2
\beta \in (1/d) \Z$, if $a = 1$.

If $a \neq \pm 1$ and $\beta$ is real, then we have $(1 \pm \tau)
\beta \in (1/d) \Gamma$, so $\beta \in (1/d) \Z$. If $a = -1$ and $\beta$ is
imaginary, then $2 \beta \in (1/d) \Gamma$ and the first assertion follows.

If $d = 1$, we  observe that if $Re(\tau) = -1/2$, the involutions $z \mapsto
\bar{z} + 1/2$, $z \mapsto - \bar{z} + 1/2 +
\tau =  - \bar{z} + i Im(\tau)$ have respectively $i/2 Im(\tau)$,
$1/4$ as fixed points, thus in both cases we can assume $\beta = 0$.

It is immediate to verify that if $Re(\tau) = 0$ and
  $a=1,  b = 1/2$, or $a=-1,  b =\tau/2$
  there are no fixed points.

The third assertion follows from the fact that $ \sigma^{-1} (z) = a
\bar{z} -a \bar{b}$, whence, for a translation
$ z \rightarrow z + c$, conjugation by $\sigma$ yields
$ z \rightarrow z + a \bar{c}$.
\hfill Q.E.D.\\

\begin{REM} Assume now that $\sigma$ is an antiholomorphic involution of an
elliptic curve $C$ as in the previous lemma. Then there are only three
possible topological types for the action of $\sigma$ on $C$;

\begin{itemize}
\item $ Fix(\sigma) = \emptyset$
\item $Fix (\sigma)$ is homeomorphic to ${\bf S}^1$:

this occurs if $|\tau| = 1$, $a = \pm \tau$, or if $Re(\tau) = -1/2$, $a = \pm
1$, $b = 0$.
\item $Fix (\sigma)$ has two connected components homeomorphic to ${\bf S}^1$:
this occurs if $Re(\tau) = 0$, $a = \pm 1$, $b = 0$.
\end{itemize}
\end{REM} {\bf Proof.}
$Fix (\sigma) = \emptyset$: there are two cases which are obviously
topological equivalent (exchange the two basis vectors of $\Gamma$, $1$, and
$\tau$).

If $Fix (\sigma) \neq \emptyset$, then we may assume $b = 0$, and the
topological type is completely determined by the integral conjugacy class of
the matrix
$A =
\left( \begin{array}{cc} m & m'\\ n & n'
\end{array} \right)
$, whose square is the identity and which has $1$ and $-1$ as eigenvalues.

If $A$ is diagonalizable then $Fix (\sigma)$ has two connected components
homeomorphic to ${\bf S}^1$, otherwise $A$ is conjugated to the matrix
$
\left( \begin{array}{cc} 0 & 1\\ 1 & 0
\end{array} \right)
$ and $Fix (\sigma)$ is homeomorphic to ${\bf S}^1$.
\hfill Q.E.D.\\

\begin{REM} Let $G$ be a Bagnera- de Franchis group: then $G$ has a first
incarnation as a group of translations of $E$, and a second one as a direct
product
$ G = G' \times T$, where $T$ is a group of translation and $G'$ is cyclic
and a subgroup of the multiplicative group. If $\hat{G}$ is an extended
Bagnera- de Franchis group, we let
$\sigma$ be an element in $\hat{G} - G$: it conjugates $G$, sending
$T$ to itself. In the first incarnation, group conjugation is given, as we
saw in the previous lemma, by complex conjugation followed by multiplication
by
$a$ (if $\sigma$ acts by
$ \sigma (z) = a \bar{z} + b$). 

It follows that the extension $ \ 0
\rightarrow G \rightarrow \hat{G}
\rightarrow {\Z}/2 \rightarrow 0 $ splits if and only if, $c'$ being the translation vector
of
$\sigma^2$, then $c'$ lies in the image of the endomorphism of $G$ given by
$s + Id$, $s$ being the action of $\sigma$. In the second incarnation, let $ z
\rightarrow \xi z$ be a generator
$g'$ of $G'$, and let $ \sigma (z) = a \bar{z} + b$: then $\sigma$ conjugates
$g'$ to the transformation $ z \rightarrow \bar{\xi} z + a (\bar{\xi}- 1)
\bar{b}$.

Whence, if the action is trivial then $ G' \cong \Z / 2 \Z$ and $ 2 b \in
\Gamma$.
\end{REM} {\bf Proof.} We need only to remark that  $ a \bar{c} \in \Gamma$
holds if and only if $c \in \Gamma$.

\hfill Q.E.D.\\

We can now give the list of all the possible groups $\hat{G}$.
\begin{LEM}
\label{action}

Let us consider the extension
$$(*) \ 0 \rightarrow G \rightarrow \hat{G} \rightarrow {\Z}/2 = <\sigma>
\rightarrow 0.$$

We have the following possibilities for the action of $\sigma$ on $G$:
\begin{enumerate}
\item If $G = {\Z}/2$, then  $\hat{G} = {\Z}/2 \times {\Z}/2$.
\item If $G = {\Z}/2 \times {\Z}/2$, then either $\sigma$ acts  as the
identity on $G$ and $\hat{G} = {\Z}/2 \times {\Z}/2 \times {\Z}/2$, if $(*)$
splits,
$\hat{G} = {\Z}/4 \times {\Z}/2$ if $(*)$ does not split, and in this latter
case the square of $\sigma$ is the generator of $T$.

Or $\sigma$ acts as $ \left( \begin{array}{cc} 1 & 0\\ 1 & 1
\end{array} \right)
$,  the sequence splits, $\hat{G} = D_4$, the dihedral group, and again the
  square of the generator of $\Z/4 \Z$ is the generator of $T$.
\item If $G = {\Z}/4$, then $\sigma$ acts as  $-Id$ on $G$ and $\hat{G} =
D_4$.
\item If $G = {\Z}/4 \times {\Z}/2$, then either
$\hat{G} = T \times D_4 \cong {\Z}/2 \times D_4$, or
$\hat{G} $ is isomorphic to the group $G_1:= <\sigma, g, t, \ | \sigma^2 =1,
\ g^4 = 1, \ t^2 = 1, \ t \sigma = \sigma t, \ tg = gt, \ \sigma g = g^{-1} t
\sigma >$, and its action on the second elliptic curve $F$ is generated by
the following  transformations:
$ \sigma (z) =  \bar{z} + 1/2$, $ g(z) = iz$, $t(z) = z + 1/2 (1+i)$. The
group $G_1$ is classically denoted by $c_1$ (cf. \cite{hs} p. 39).

In particular, in both cases  $(*)$ splits.

\item If $G = {\Z}/3$, then $\sigma$ acts as $-Id$ on $G$ and
$\hat{G} = {\cal S}_3$.
\item If $G = {\Z}/3 \times {\Z}/3$, then we may choose $G'$ so that $\sigma$
acts as
$-Id \times Id$ on $G = G' \times T$ and $\hat{G} = {\cal S}_3 \times {\Z}/3$.
\item If $G = {\Z}/6$, then $\sigma$ acts as $-Id$ on $G$ and $\hat{G} = D_6$.
\end{enumerate}
\end{LEM}  {\bf Proof.}  Thanks to \ref{gcyc} we know that in the cases
1,3,5,7 the sequence $(*)$ splits.\\  Case 1: if $G = {\Z}/2$, then $\sigma$
acts as the identity on $G$ whence $\hat{G} = {\Z}/2 \times {\Z}/2$. \\

Cases 3,5,7: from the previous remark we can immediately determine $\hat{G}$
in these cases in which $G$ is respectively equal to $  \Z/4, \Z/3, \Z/6$. In
fact
$|\xi| = 1$, whence $\bar{\xi} = \xi^{-1}$. Thus $\sigma$ acts as
$-Id$ on $G = G'$, the extension splits and $\hat{G}$ is a dihedral group:
respectively $D_4$, $D_3 = {\cal S}_3$, $D_6$.\\

Case 6: $G = \Z/3 \times \Z/3$. Any element $\sigma$ of order $2$ makes
$(*)$ split.\\  Since $T = \Z/3$ is invariant and $\sigma$ acts as $-1$ on $
G/T$, we can take eigenspaces $T$ and $G'$ such that
$\sigma$ acts on $G'$ as $-Id$ and as $+ Id$ or $ - Id$ on $T$. The second
case is not possible, since, looking at the first incarnation, we obtain that
the two eigenvalues of the action of $\sigma$ on the lattice
$\Gamma$ of $E$ (either $=1$ or $-1$) reduce to $-1$ modulo $3$, thus they
equal $-1$ and $\sigma$ acts holomorphically, a contradiction.

Case 2: $G = {\Z}/2 \times {\Z}/2$.

$\sigma$ acts on $G$ and trivially on $T= {\Z}/2$, whence either the action is
the identity, or is given by
$\left( \begin{array}{c c}  1 & 0\\ 1 & 1
\end{array} \right)
$. In the latter case, since as we noticed
  $H^2({\Z}/2, G) = \frac{ker(\sigma - 1)}{(1 + \sigma)G} $,  then
$H^2(\Z/2, G) = 0$, so the extension splits and
$\hat{G} = D_4$ (and an element of order $4$ is given by $\sigma e_1$ whose
square is indeed $e_2$, the generator of $T$).
\\

If $\sigma$ acts as the identity on $G$,
$\hat{G}$ is abelian, and the exact sequence splits if and only if there is no
element of order $4$. Hence the only possibilities for $\hat{G}$ are the
following:
\begin{itemize}
\item $\hat{G} = {\Z}/2 \times {\Z}/2  \times {\Z}/2$, and the sequence
splits,
\item $\hat{G} = {\Z}/4 \times {\Z}/2$.
\end{itemize}

In the latter case, we recall that by lemma 4.2 the square of
$\sigma$ is the generator of $T$.

Case 4: $G$ is, non canonically, isomorphic to ${\Z}/4 \times T$, with
$T \cong {\Z}/2$.

$\sigma^2$ belongs to $T$, moreover $\sigma$ acts as $Id$ on $T$ and as $-
Id$ on
$ G/ T$.

Case 4.I: there is an element $g$ such that $\sigma g \sigma^{-1} = g^{-1}$.
In this case, for each $ g \in G$ we get $(\sigma g)^2 = \sigma g \sigma^{-1}
g (\sigma)^2 = (\sigma)^2$ because $ Id + s = 0$.

Whence, we distinguish two cases:

Case 4.I.a): $ (\sigma)^2 = 1$. Here, it follows rightaway that

$\hat{G} \cong T \oplus D_4$.

Case 4.I.b): $ (\sigma)^2 = t$, $t$ being the generator of $T$.

This case can be excluded as follows: first of all, notice that for the
elliptic curve $F$ we have $\tau = i$, and $ t$ equal to the translation by
the semiperiod
$ 1/2(1+i)$.

We have $\sigma g \sigma^{-1} = g^{-1}$ for each element of order
$4$, in particular for the multiplication by $i$. Furthermore,  up to
multiplying
$\sigma$ by a power of the above element, we  may assume (cf. the discussion
preceding lemma \ref{tra}) that $a=1$, i.e., that $ \sigma (z) = \bar{z} + b$.

Then, our conjugation relation reads out as follows:

$ a (-1-i) \bar{b} \equiv 0 \ (mod \ \Gamma)$, or, more simply, as

$ (1+i) \bar{b} \equiv 0 \ (mod \ \Gamma)$.

Whence, $b$ is a semiperiod $ b = 1/2 ( x + iy)$ ($ x,y \in \Z$) and the above
condition means that $ x+y \equiv 0 \ ( mod \ 2)$.

But a previous calculation shows that the square of $\sigma$ is  a
translation by
$ a \bar{b} + b = x  \ \equiv 0 \ (mod  \ \Gamma)$, contradicting $
(\sigma)^2 = t$.

There remains the following

Case 4.II: $\sigma g \sigma^{-1} = g^{-1} t$.

In this case we can always choose $ \sigma$ such that
$ (\sigma)^2 = 1$. In fact, if $ (\sigma)^2 = t$, then
  $(\sigma g)^2 = \sigma g \sigma^{-1} g (\sigma)^2 = g^{-1} t g (\sigma)^2 =
g^{-1} t g t = 1$.

Again, here $\sigma g \sigma^{-1} = g^{-1} t$ for each element of order $4$ in
$G$. As in case 4.I.b, we choose $g$ as the element given (on $F$) by
multiplication by $i$, and we have that  $ t$ is the translation by the
semiperiod $ (1+i)/2$.

In this case, however, we can only multiply $\sigma$ by the square of
$g$, whence we may only assume $ a= 1$ or $ a=i$.

Since we are assuming $ (\sigma)^2 = 1$, we get $a \bar{b} + b
\equiv 0 \ ( mod \ \Gamma)$, furthermore $\sigma g \sigma^{-1} = g^{-1} t$
reads out, if $ \sigma (z) = a \bar{z} + b$, as follows:

(**) $ - a (1+i) \bar{b} \equiv 1/2(1+i) \ (mod \ \Gamma)$, i.e.,

$ (1+i) [  a \bar{b} + 1/2 ] \in \Gamma$ whence, as before,
$b$ is a semiperiod $ b = 1/2 ( x + iy)$ ($ x,y \in \Z$),

and (**) simply means $ x + y \equiv 1 \ ( mod \ 2)$.

Up to exchanging $ \sigma$  with $ \sigma t$, we may assume

$ x=1, y=0$, and either

\begin{itemize}
\item
$ \sigma (z) =  \bar{z} + 1/2$
\item
$ \sigma (z) = i \bar{z} + 1/2$
\end{itemize}

However, the condition $a \bar{b} + b
\equiv 0 \ (mod \ \Gamma)$ is only verified in the former case.

Our group $\hat{G}$ is generated by elements $ t,\sigma$ of order
$2$, $g$ of order $4$, $t, g^2$ are in the centre, moreover the relation
$\sigma g \sigma^{-1} = g^{-1} t$ holds: $\hat{G}$  is thus isomorphic to the group
$G_1$.

\hfill Q.E.D. \\

\section{The split case}

In this section we want to treat the case in which the extension

$(*) \ \ 0 \rightarrow G \stackrel{i} \rightarrow \hat{G}
\stackrel{\pi} \rightarrow {\Z}/2  \rightarrow 1,$

that we studied in the last section, splits.

Let $\sigma$ be an antiholomorphic involution on an elliptic curve
$C$.  We recall the analysis of the antiholomorphic maps
$\sigma$ with $\sigma^2$ equal to a translation given in the last section, it
can be summarized as follows
$$
\begin{tabular}{|c|c|}
\hline Case & $a$    \\
\hline
$Re(\tau) =0, \ |\tau|> 1 $ & $\pm 1$    \\
$\tau =i$ & $1 \equiv -1$  \\
$\tau =i$ & $i \equiv -i$  \\
$|\tau| =1, \ -1/2 < Re(\tau) <0$ & $\pm \tau$ \\
$\tau =\rho$ & $1 \equiv \rho \equiv \rho^2$ \\
$\tau =\rho$ & $-1 \equiv -\rho \equiv - \rho^2$ \\
$Re(\tau) = -1/2, \ |\tau| >1$ & $\pm 1$ \\
\hline

\end{tabular}
$$ In the above table we have used the fact that if $\tau = i$, by conjugating
with the automorphism of $C$ given by multiplication by
$i$ we can reduce the case $a = -1$ to the case  $a = 1$, the case $a = i$ to
the case $a = -i$.

Analogously, if $\tau = \rho$, by conjugating with the automorphisms of $C$
given by multiplication by $\rho$, or by $\rho^2$, we can reduce the cases $a
= \rho,
\rho^2$ to the case  $a = 1$, the cases $a = -\rho, - \rho^2$ to the case $a =
-1$.

We consider first of all the action of $\hat{G}$ on the curve $E$. Let $\sigma
(z) = a \bar{z} + b$ be any antiholomorphic involution of $\hat{G}$ acting on
$E$. By lemma \ref{tra} we know that we have the following possibilities for
$b$: if $Re(\tau) = 0$, $a = 1$, we have $b = 0, 1/2$; if $Re(\tau) = 0$, $a =
-1$, we have $b = 0,
\tau/2$; if $|\tau| = 1$, or $Re(\tau) = -1/2$, we can assume $b = 0$.

\begin{DEF} We call cases $1,3, 5, 7$ of \ref{action} the simple dihedral
cases: they are characterized by the property that
  $G = \Z/q$, $\hat{G} = D_q$ (where  we have respectively $q = 2,4,3,6$).
\end{DEF}

\begin{REM} In a dihedral case with $G = \Z/q$, $\hat{G} = D_q$,
$\hat{G}$ is generated by elements $g$ and $\sigma$, where
$g$ acts on $E$ by $g(z) = z + c$,
  $c$ being a torsion element in $Pic^0(E)$ of order precisely $q$
  and such that  $a \bar{c} \equiv -c \ (mod \ \Gamma)$.
\end{REM}

{\bf Proof.} We choose generators  $\sigma$ and $g$ such that $\sigma \circ g
\equiv g^{-1}
\circ \sigma$, equivalently $a \bar{c} \equiv -c \ (mod \ \Gamma)$.
\hfill Q.E.D. \\

\begin{LEM}
\label{dihedral} The subgroup $H_q$ of $Pic^0(E)_q$ on which
$\sigma$ acts as $-Id$ is isomorphic to $\Z/q$, if  and only if either
$|\tau| = 1$, $a = \pm \tau$, $ q =2,4,3,6$ or $Re(\tau) = -1/2$, $a = \pm
1$, $q = 2,4,3,6$, or $Re(\tau) = 0$, $a = \pm 1$ and $q = 3$.

If $Re(\tau) = 0$, $a = \pm 1$, $q = 2,4,6$, then $H_q$ is isomorphic to $\Z/q
\times \Z/2$.
\end{LEM}

{\bf Proof.} We know that the map $z \mapsto a \bar{z}$ is represented by the
matrix $
\left( \begin{array}{cc}  m & m' \\  n & n'
\end{array} \right)
$  on $\Z \oplus \Z \tau$, hence the action on $Pic^0(E)$ with basis
$1/q$, $\tau/q$ is given by the reduction modulo $q$ of the above integral
matrix.

We have to solve the equation $a \bar{c} + c \equiv 0 \ (mod \
\Gamma)$, with $c \in Pic^0(E)_q$, so we consider the kernel of the reduction
of the matrix
$ M = \left( \begin{array}{cc}  m+ 1 & m' \\  n & n'+1
\end{array} \right)
$ modulo $q$.

For $|\tau| = 1$, $a = \pm \tau$, we have: $m = n' = 0$, $n = m' =
\pm 1$, while for $Re(\tau) = -1/2$, $a = \pm 1$, we have: $n = 0$,
$m = \pm 1$, $n' = m' = -m$, therefore the kernel of the linear map
$a \bar{z} + z$ on $Pic^0(E)_q$ is isomorphic to $\Z/q$.

It remains to consider the case $Re(\tau) = 0$, $a = \pm 1$. We have: $n = m'
= 0$, $m = \pm 1$, $n' = -m$.  So the matrix
$ M = \left( \begin{array}{cc}  2 & 0 \\  0 & 0
\end{array} \right)
$  if $a=1$, $ M = \left( \begin{array}{cc}  0 & 0 \\  0 & 2
\end{array} \right)
$ if $a = -1$. This tells us that the kernel of $M$ is isomorphic to
$\Z/3$ if $q = 3$, while it is isomorphic to $\Z/2 \times \Z/q$, if
$q = 2, 4, 6$.
\hfill Q.E.D. \\

\begin{LEM}
\label{z2z2} If $G = \Z/2 \times \Z/2$ and $\hat{G} = \Z/2 \times
\Z/2 \times \Z/2$, we must have $Re(\tau) = 0$, $a = \pm 1$, and $b = 0,1/2$
in case $a =1$, $b = 0, \tau/2$, in case $a = -1$.

If $G = \Z/2 \times \Z/2$ and $\hat{G} = D_4$, we have either
$|\tau| = 1$, $a = \pm \tau$, or $Re(\tau) = -1/2$, $a = \pm 1$.

\end{LEM}

{\bf Proof.} If $G = \Z/2 \times \Z/2$ and $\hat{G} = \Z/2
\times \Z/2 \times \Z/2$, since $\hat{G}$ is abelian, the map $d
\mapsto a \bar{d}$ is the identity on the points of two torsion and one can
easily see (cf. the proof of the previous lemma) that this can happen if and
only if
$Re(\tau) = 0$ and $a = \pm 1$.

If $G = \Z/2 \times \Z/2$, $\hat{G} = D_4$, we know from
\ref{action} that $\sigma$ acts as $
\left( \begin{array}{cc}  1 & 0\\  1 & 1
\end{array} \right)
$ on the points of two torsion of $E$, then we have either $|\tau| =1$, $a =
\pm
\tau$, or $Re(\tau) = -1/2$, $a = \pm 1$.

\hfill Q.E.D. \\

\begin{LEM} If $G = \Z/4 \times \Z/2$, then again we must have
$Re(\tau) = 0$, $a = \pm 1$, and $b = 0,1/2$ in case $a =1$, $b = 0,
\tau/2$, in case $a = -1$.

In the case $\hat{G} \cong T  \times D_4 \cong \Z/2 \times D_4$, the subgroup
$G
\subset Pic^0(E)_4$ consists of the points
$c$ of four torsion satisfyng the equation $a \bar{c} \equiv -c$.

In the case $\hat{G}= G_1$, we can choose generators of $\hat{G}$,
$\sigma(z) = a \bar{z}+ b$,
$g(z) = z + \eta$, $t(z) = z + \epsilon$, with $\eta$ of order
$4$, $\epsilon$ of order $2$ such that:

\begin{itemize}
\item if $a = 1$, $\epsilon = 1/2$ and  $\eta = 1/4 + \tau/4$, $b = 0, 1/2$;
\item if $a = -1$, $\epsilon = \tau/2$, $\eta =  1/4 +  \tau/4$, $b = 0,
\tau/2$.
\end{itemize}

\end{LEM}

{\bf Proof.} In both cases the action of $\sigma$ on $Pic^0(E)_2 \subset G$ is
trivial, whence we get the same conditions upon $a,b,\tau$ as in the previous
lemma.

Let  $g$ be any element in $G$ of order $4$: then the following relations hold

  $\sigma \circ g = g^{-1} \circ \sigma$, in the case $\hat{G} = \Z/2 \times
D_4$,

$\sigma \circ g = g^{-1} t\circ \sigma$, in the case $\hat{G} = G_1$ (where
$t$ is the generator of $T$).

If $\hat{G} \cong \Z/2 \times D_4$, then the above equation shows that $G$ is
contained in the subgroup of $Pic^0(E)_4$ where $\sigma$ acts as $-Id$. But we
have seen in the proof of \ref{dihedral} that this subgroup is isomorphic to
$\Z/2 \times \Z/4$ if $Re(\tau) = 0$, $a = \pm 1$, whence it equals $G$.

If $\hat{G} = G_1$, then the relation
$a \bar{\eta} + \eta \equiv \epsilon \ (mod \ \Gamma)$, gives
$2Re(\eta) \equiv \epsilon  \ (mod \ \Gamma)$ if $a = 1$, $2i Im(\eta) \equiv
\epsilon  \ (mod \ \Gamma)$ if $a = -1$, therefore
  $\epsilon \equiv 1/2$  if $a=1$,
$\epsilon \equiv \tau/2$ if $a = -1$.

If $a=1$, it follows that $Re(\eta) = \pm 1/4$, and since $2 \eta$ is of two
torsion but distinct from $\epsilon$, we get that
$\eta = \pm 1/4 + \pm \tau/4$. Up to replacing $g$ with its inverse, or up to
composing with $t$, we achieve $\eta = 1/4 + \tau/4$, and an entirely similar
argument works in the case $a=-1$.

\hfill Q.E.D. \\

\begin{LEM} If $G = \Z/3 \times \Z/3$, $\hat{G} = {\cal S}_3 \times
\Z/3$, then $G$ is isomorphic to the group $Pic^0(E)_3$ and the action of any
antiholomorphic involution $\sigma$ on $Pic^0(E)_3$ has $1$ and $-1$ as
eigenvalues.

Therefore the datum of a $\hat{G}$- action is equivalent to the datum of an
isomorphism class of an antiholomorphic involution.

We have thus the usual following possibilities for
$\sigma$:
\begin{itemize}
\item $Re(\tau) = 0$, $a = 1$, $b = 0, 1/2$,
\item $Re(\tau) = 0$, $a = -1$, $b = 0, \tau/2$,
\item $|\tau| = 1$, $a = \pm \tau$,
\item $Re(\tau) = -1/2$, $a = \pm 1$.
\end{itemize}
\end{LEM}

{\bf Proof.} Since $\sigma$ is an antiholomorphic involution, the eigenvalues
of the action of $\sigma$ on the lattice
$\Gamma$ of $E$ are $1$ and $-1$, thus the same holds for their reduction
modulo
$3$, therefore we have all the possible values of
$a$ and $b$ (see lemma \ref{tra}).
\hfill Q.E.D. \\

\bigskip

Let us now consider the action of $\hat{G}$ on $F$.

\begin{LEM}
\label{fixg} Consider a cyclic group $G'$ of automorphisms on an elliptic
curve
$F$, generated by a transformation $g$ having a fixed point $0$.

Then the order  $d$ of $g$ equals $2, 4 , 3, 6$ and $Fix(g)$ is as follows:
\begin{itemize}
\item If $d = 2$, $Fix(g)$ is the subgroup $F_2$ of the 2-torsion points of
$F$.
\item If $d =4$, $Fix(g)$ is the subgroup of $F_2$ isomorphic to
$\Z/2$ generated by $(1+i)/2$.
\item If $d =3$, $Fix(g)$ is the subgroup of $F_3$ isomorphic to
$\Z/3$ generated by $(1-\rho)/3$.
\item If $d =6$, $Fix(g)$ is only the origin.
\end{itemize}
\end{LEM} {\bf Proof.} The proof is a simple computation, we notice that if $d
=3$, or $d = 6$, $\Gamma \cong \Z[\rho]/(\rho^2 + \rho + 1)$.
\hfill Q.E.D. \\

\begin{REM}
\label{fixggg}

Let $\sigma$ be an antiholomorphic involution in
$\hat{G}$ acting on $F$.

If $\sigma \circ g = g^{-1} \circ \sigma$, then $\sigma$ acts on the sets
$Fix(g)$, $Fix(g^2) - Fix(g)$, $Fix(g^3) - Fix(g)$.

Moreover, $Fix(g)$ is a subgroup of $Pic^0(F)$ and  $\sigma$ acts on
$Fix(g)$ by an affine transformation.

Notice that $g$ acts trivially on $ Fix(g)$, whence the action on $ Fix(g) =
Fix (g^{-1})$ is independent on the choice of
$\sigma$ in the simple dihedral case, when we take for $g$ a generator of the
subgroup $G$: we thus get in this case a topological invariant of the real
hyperelliptic surface.
\end{REM}
\begin{COR} In the cases $5,7$ of \ref{action} we may pick as the origin in
$F$ a common fixed point of $g$ and $\sigma$, therefore if
$G = \Z/3$ we may assume $\sigma(z) = \pm \bar{z}$, if $G = \Z/6$ we may
assume
$\sigma(z) = \bar{z}$.
\end{COR} {\bf Proof.} We have already observed that we may choose
$a = \pm 1$. If $G = \Z/6$, by exchanging $\sigma$ with $g^n \circ
\sigma$ we see that we can assume $a = 1$.

If $G = \Z/3$, $g(z) = \rho z$, $Fix(g^2) - Fix(g)$ is only the origin, and
$\sigma$ must act on $Fix(g^2) - Fix(g)$ by the previous remark.

If $G = \Z/6$, $g(z) = - \rho z$, $Fix(g)$ is only the origin and
$\sigma$ must act on $Fix(g)$ by the previous remark.
\hfill Q.E.D. \\
\begin{LEM} If $G = \Z/4$, $\hat{G} = D_4$ we may assume $a = 1$, and either
$\sigma(z) = \bar{z}$, or $\sigma(z) = \bar{z} + (1+i)/2$.
\end{LEM}

{\bf Proof.} We know that $a = \pm 1, \pm i$, but by exchanging $\sigma$ with
$g^n \circ \sigma$ we can assume $a = 1$.

By \ref{fixg} and \ref{fixggg} we know that $\sigma$ acts on the subgroup of
$F_2$ generated by $(1+i)/2$, that concludes the proof.
\hfill Q.E.D. \\
\begin{LEM} If $G = \Z/3 \times \Z/3$, $\hat{G} = D_3 \times \Z/3$, we may
choose generators of $\hat{G}$ acting on $F$ given by
$\sigma(z) = - \bar{z}$, $g(z) = \rho z$, $t(z) = z + (1 - \rho)/3$.
\end{LEM} {\bf Proof.} By \ref{action} we know that we may assume that
$\hat{G}$ is generated by $\sigma(z) = a \bar{z} + b$, $g(z) =
\rho z$, $t(z) = z + (1- \rho)/3$ such that $\sigma \circ g = g^{-1}
\circ \sigma$, $\sigma \circ t = t \circ \sigma$.  We have already remarked
that we may choose $a = \pm 1$, and one can easily see that the relation
$\sigma
\circ t = t \circ \sigma$ holds only for $a = -1$.

Since $\sigma \circ g = g^{-1} \circ \sigma$, we can assume $b =0$, by
\ref{fixggg} as in the case $G = \Z/3$.

\hfill Q.E.D. \\

\begin{LEM} If $G = \Z/4 \times \Z/2$ we have the following cases:
\begin{itemize}
\item If $\hat{G} = \Z/2 \times D_4$, we can choose the following set of
generators of $\hat{G}$ acting on $F$: $g(z) = i z$,
$t(z) = z + (1+i)/2$, $\sigma(z) = \bar{z}$.
\item If $\hat{G} =G_1$, we can choose the following set of generators of
$\hat{G}$ acting on $F$: $g(z) = i z$, $t(z) = z + (1+i)/2$, $\sigma(z) =
\bar{z} + 1/2$.

\end{itemize}
\end{LEM} {\bf Proof.} For the case $\hat{G} = G_1$ see lemma
\ref{action}.

Assume $\hat{G} = \Z/2 \times D_4$, then we have already observed in
\ref{action} that we may assume that $\hat{G}$ is generated by the
transformations $g(z) = i z$, $t(z) = z + (1+i)/2$, $\sigma(z) = a
\bar{z} + b$, such that $\sigma \circ g = g^{-1} \circ \sigma$,
$\sigma \circ t = t \circ \sigma$.

Therefore by composing $\sigma$ with a power of $g$ we may assume $a = 1$, and
the relation $\sigma \circ g = g^{-1} \circ \sigma$ gives by \ref{fixggg}
either
$b =0$, or $b = (1+i)/2$. The second case can be excluded by exchanging
$\sigma$ with $t \circ \sigma$.
\hfill Q.E.D. \\

\begin{LEM}  If $G = \Z/2$, $\hat{G} = \Z/2 \times \Z/2$ generated by
$\sigma(z) = a \bar{z} + b$, $g(z) = -z$, we have the following possibilities
for
$F$, $a$, $b$:
\begin{itemize}
\item $Re(\tau) = 0$, $a = 1$, $b$ an element in $F_2$.
\item $|\tau| = 1$, $a = \tau$, $b = 0$.
\item $Re(\tau) = -1/2$, $a = 1$, $b = 0$.
\end{itemize} Moreover there are three different topological types of the
action of the affine transformation $\sigma$ on $Fix (g)$, where $g$ is a
generator of
$G$ (see \ref{fixggg}), namely:
\begin{itemize}
\item If $Re(\tau) = 0$, $a = 1$, $b =0$, then $\sigma$ acts as the identity
on
$Fix(g)$.
\item If $Re(\tau) = 0$, $a = 1$, $b \in Pic^0(E)_2$ of order 2, then $\sigma$
acts as a translation on $Fix(g)$.
\item If $|\tau| = 1$, $a = \tau$, or $Re(\tau) = -1/2$, $a = 1$, then
$\sigma$ acts as a linear map with matrix
$
\left( \begin{array}{cc}  1 & 1\\  0 & 1
\end{array} \right)
$ on $Fix(g)$.
\end{itemize}

\end{LEM}       {\bf Proof.}  Since $\hat{G}$ is abelian, we have the
condition
$2 b \equiv 0 \ (mod \ \Gamma)$, while the condition
$\sigma^2 = 1$ reads $ a \bar{b} + b \equiv 0 \ (mod \ \Gamma)$.

By exchanging $\sigma$ with $- \sigma$, $a$ gets multiplied by $-1$, so we
obtain the statement on the possible values of $a$.

Assume now that $Re(\tau) = 0$, $a = 1$.

The conditions $2 b \equiv 0 \ (mod \ \Gamma)$, $a \bar{b} + b =
\bar{b} + b = 2 Re(b) \equiv 0 \ (mod \ \Gamma)$ give us the following
possibilities for $b$: $b = 0, \ 1/2, \ \tau/2 \ (1 +
\tau)/2$.

Assume now that $|\tau| = 1$, $a = \tau$.  The condition $2 b \equiv 0 \ (mod
\
\Gamma)$, allows us to write $b = x/2 + \tau y/2$, with
$x,y \in \{0,1\}$.

The condition $a \bar{b} + b = \tau  \bar{b} + b \equiv 0 \ (mod \
\Gamma)$ implies $x \equiv y  \ (mod \ 2)$, so either $ b = 0$, or
$b = (1 + \tau)/2$, but by conjugation  with the translation
$\phi(z) = z + 1/2$, we can assume $b = 0$, $\sigma(z) = \tau
\bar{z} $.

Assume now that $Re(\tau)  = -1/2$.

We may assume $a = 1$, and the condition $2 b \equiv 0 \ (mod \
\Gamma)$ allows us to write $b = x/2 + \tau y/2$, with $x, y \in
\{0,1\}$, while the condition $a \bar{b} + b \equiv 0 \ (mod \
\Gamma)$ reads $2 Re(b) = x - y/2 = n \in \Z$. Thus $y \equiv 0 \ (mod \ 2)$
and either $b = 0$, or $b = 1/2$, but by translating by
$\phi(z) = z + \tau/2$, we can assume $b = 0$ and $\sigma(z) =
\bar{z}$.

Finally we know from \ref{fixggg} that $\sigma$ acts as an affine
transformation on $Fix(g) = Pic^0(E)_2$, that this action is a topological
invariant, and we observe that there are the following 4 different types of
affine transformations acting on $Pic^0(E)_2$:

the identity, a translation, the map $A =
\left( \begin{array}{cc}  1 & 1\\  0 & 1
\end{array} \right)
$, the map $A$ composed with a translation.

Then the statement follows by an easy computation.
\hfill Q.E.D. \\

There remains to treat the case $G = \Z/2 \times
\Z/2$.

\begin{LEM} Let $G = \Z/2 \times \Z/2$, then either $\hat{G} = \Z/2
\times \Z/2 \times \Z/2$, or $\hat{G} = D_4$ and we can choose generators of
$\hat{G}$ as follows:
$g(z) = -z$, $t(z) = z + \epsilon$, $\epsilon$ of order $2$,
$\sigma(z) = a \bar{z} + b$, where
\begin{itemize}
\item if $\hat{G}  = \Z/2 \times \Z/2 \times \Z/2$ we have the following
cases:

$Re(\tau) = 0$, $a = 1$: if $\epsilon = 1/2$, $b =0, \tau/2$; if
$\epsilon = \tau/2$, or $\epsilon = 1/2 + \tau/2$, we can choose $b = 0, 1/2$.

$|\tau | = 1$, $a = \tau$, $\epsilon = (1+ \tau)/2$, $b = 0$.

$Re(\tau) = -1/2$, $a = 1$, $\epsilon = 1/2$, $b = 0$.

\item if $\hat{G}  = D_4$ we have the following cases:

$Re(\tau) = 0$, if $a = 1$, $ b = \tau/4, \tau/4 + 1/2$, $\epsilon =
\tau/2$; if $a =- 1$, $ b = -1/4, -1/4 + \tau/2$, $\epsilon = 1/2$;

$|\tau| =1$, if $a = \tau$, $ b= 1/4 - \tau/4$, $\epsilon = 1/2 +
\tau/2$; if $a = -\tau$, $ b= 1/4 + \tau/4$, $\epsilon = 1/2 +
\tau/2$;

$Re(\tau) = -1/2$, if $a =  1$, $b = 1/4 + \tau/2$, $\epsilon = 1/2$; if $a =
-1$, $b = -1/4$, $\epsilon = 1/2$;
\end{itemize}
\end{LEM} {\bf Proof.}  Assume first of all $G = \Z/2 \times \Z/2$,
$\hat{G} = \Z/2 \times \Z/2 \times \Z/2$.

By exchanging $\sigma$ with $- \sigma$, $a$ gets multiplied by $-1$, thus we
have the statement for $a$.

Since $\hat{G}$ is abelian, we have the condition $2 b \equiv 0 \ (mod \
\Gamma)$, while the condition $\sigma^2 = 1$ reads $ a
\bar{b} + b \equiv 0 \ (mod \ \Gamma)$.

Assume first of all $Re(\tau) = 0$, $a =1$, then as in the case $G =
\Z/2$, the conditions $2 b \equiv 0 \ (mod \ \Gamma)$, and $ a
\bar{b} + b \equiv 0 \ (mod \ \Gamma)$, give us  the following possibilities
for
$b$: $b = 0, \ 1/2, \ \tau/2, \ (1+\tau)/2$.

Now, by composing $\sigma$ with $t$, we obtain the statement.

If $|\tau| =1$, $a = \tau$, as in the case $G = \Z/2$, we find
$\sigma(z) = \tau \bar{z}$.  Since $\sigma$ must commute with $t$, we must
have $ a \bar{\epsilon} = \tau \bar{\epsilon} = \epsilon$, which yields
$\epsilon = (1 +
\tau)/2$.

If $Re(\tau) = -1/2$, we can assume $a = 1$, and as in the case $G =
\Z/2$, we find $\sigma(z) = \bar{z}$.  The condition $ a
\bar{\epsilon} = \bar{\epsilon} = \epsilon$ yields $\epsilon = 1/2$.

If $G = \Z/2 \times \Z/2$, $\hat{G} = D_4$, we know by \ref{action} that we
may choose generators $g$, $t$, $\sigma$ such that
$\sigma(-z) = - \sigma(z) + \epsilon$ or equivalently $2 b \equiv
\epsilon \in 1/2 \Gamma - \Gamma$, $\sigma(z + \epsilon) = \sigma(z) +
\epsilon$, or equivalently $a \bar{\epsilon} \equiv \epsilon \ (mod \
\Gamma)$, while
$\sigma^2 = 1$ yields $a \bar{b} + b \in
\Gamma$. An element of order $4$ in $\hat{G}$ is $\sigma \circ g$, whose
square is $t$.

Assume $Re(\tau) = 0$, $a = 1$, then the conditions $2b \in 1/2
\Gamma - \Gamma$ and $2 Re(b) \equiv 0 \ (mod \ \Gamma)$ imply $2b
\equiv \tau/2 \   (mod \ \Gamma)$, therefore we can choose $b = \pm
\tau/4$, or $b = \pm \tau/4 + 1/2$, $\epsilon = \tau/2$. But by composing
$\sigma$ with $t$, we can assume $b = \tau/4$, or $b =
\tau/4 + 1/2$.

A similar computation gives the result for $a = -1$.

If $|\tau| = 1$, $a = \tau$, then the conditions $2b \in 1/2 \Gamma -
\Gamma$ and
$\tau \bar{b} + b  \equiv 0 \ (mod \ \Gamma)$ tell us that we can choose $b =
1/4 - \tau/4$, $\epsilon = (1 + \tau)/2$.

In fact by the first condition we can write $2b = n/2 + m/2 \tau$, with $n, m
\in
\{0,1\}$.
$2(a \bar{b} + b) = 2(\tau \bar{b} + b) = (n + m)/2 + \tau (n + m)/2
\in 2 \Gamma$, thus $n + m \equiv 0 \ (mod \ 4)$, and we obtain $b
\equiv 1/4 - \tau/4 \ (mod \ \Gamma)$.  A similar computation gives the result
for $a = -\tau$.

If $Re(\tau) = -1/2$, $a = 1$, the conditions $2b \in 1/2 \Gamma -
\Gamma$ and $\bar{b} + b  = 2 Re(b) \equiv 0 \ (mod \ \Gamma)$ tell us that we
can choose $b \equiv  1/4 + \tau/2$, $\epsilon = 1/2$.

In fact we can write $2b = n/2 + m/2 \tau$, with $n, m \in \{0,1\}$, then $2
Re(b) =  n/2 - m/4 \equiv 0  \ (mod \ \Gamma)$ and we have
$m \equiv 2n \ (mod \ 4)$, whence we may take $b  \equiv  1/4 +
\tau/2$, and therefore $\epsilon = 1/2$.

A similar computation gives the result for $a = -1$.

\hfill Q.E.D. \\

\section{The non split case}

In this section we want to treat the  case where the exact sequence $(*)$ does
not split. We have

$G = T  \times {\Z}/2 \cong {\Z}/2 \times {\Z}/2$,
$$0 \rightarrow  T \times {\Z}/2 \stackrel{j} \rightarrow {\Z}/4 \times {\Z}/2
\rightarrow {\Z}/2 \rightarrow 1$$  and
${\sigma}^2 = t$,  where $t$ is the generator of $T$, a translation.

We get as generators
${\sigma}$,
  and $g$, where $g$ acts on $F$ by multiplication by $-1$.

  We consider first the action of $\hat{G}$ on $F$.

  Let $ \sigma(z) = a \bar{z} + b$, then $b$ is a half-period
  for the following reasons:  $ \sigma, g$ commute, a condition which is
equivalent
  to $ 2 b \equiv 0 \ (mod \ \Gamma)$, and moreover $b \neq 0$, else
  $\sigma^2$ is the identity.

  If $b$ is a half-period, then necessarily $\sigma^2$ has order at most
  two, and it has order exactly two iff $(a \bar{b} + b) \not \in \Gamma$.

\begin{LEM} The case  $ Re (\tau) =0 , |\tau| > 1$ is impossible for the curve
$F$. In the case $ \tau = i$ for the curve $F$ we can assume $ a =i$.
\end{LEM}  {\bf Proof.}  In fact, assume $a= \pm 1$: then, since $(a \bar{b}
+ b)
\in  1/2 \Gamma -\Gamma$, either $ 2 Re (b) \in  1/2 \Gamma -\Gamma$ or
$ 2i  Im (b) \in  1/2 \Gamma -\Gamma$, contradicting $ 2 b \in   \Gamma$
which in this case is equivalent to $ 2 Re (b) \in \Gamma$ and
$ 2i  Im (b) \in   \Gamma$.

Finally, if $ \tau = i$ the case $ a= -i$ can be reduced to the case
$ a=i$ modulo changing coordinates via the automorphism of $F$ given by
multiplication by $i$.

\hfill Q.E.D. \\

\begin{PROP} In the non split case, let $|\tau| =1$ for $F$. Then we can
choose an element
$\sigma$ of $\hat{G}$ of order 4 such that $\sigma(z) = \tau \bar{z} + 1/2$.
In this case $\sigma^2$ is the translation by the vector $t = (1 + \tau)/2$.

If $Re(\tau) = - 1/2$, we can choose an element $\sigma$ of $\hat{G}$ of
order 4 such that $\sigma(z) = \bar{z} + \tau/2$. In this case $\sigma^2$ is
the translation by the vector $t = 1/2$.

If $\tau = \rho $, the case $a = \rho$ can be reduced to the case $a =1$, by
conjugating with the automorphism of $F$ given by multiplication by $\rho$.
Thus we can always assume $\sigma(z) = \bar{z} + \tau/2$, $\sigma^2(z)  =  z
+ 1/2$.

\end{PROP}  {\bf Proof.}  We observe first of all that by exchanging $\sigma$
with
$\sigma g$, $a$ gets multiplied by $-1$, while  conjugating $\sigma$ with an
automorphism of $F$ given by multiplication by $\lambda$
$a$ gets multiplied  by $\lambda^2$.

Therefore, we may assume $a = \tau $ if
$|\tau| = 1$, and
$a = 1$ if
$Re(\tau) = -1/2$.

Furthermore, by exchanging $\sigma$ with
$\sigma^{-1}$, we can substitute $b$ with $-a \bar{b}$.

{\bf Case $|\tau| =1$}.

Since $2b \in \Gamma$, we may write $b = x/2 + y \tau/2$, with $x,y \in
\{0,1\}$.

The condition $ a \bar{b} + b = \tau  \bar{b} + b \in 1/2 \Gamma - \Gamma$
yields
$x + y \equiv 1 \ (mod \ 2)$, so either $b = 1/2$, or $b = \tau/2$. But by
exchanging $b$ with $-\tau \bar{b}$, we can assume $b = 1/2$.

Whence we have: $\sigma(z) = \tau \bar{z} + 1/2$, and therefore
$\sigma^2(z) = z+ 1/2 + \tau/2$.

{\bf Case $Re(\tau) = -1/2$}.

The condition $a \bar{b} + b = \bar{b} + b \in 1/2 \Gamma - \Gamma$ is
equivalent to the condition $2 Re(b) \in 1/2 \Gamma - \Gamma$, that implies
$Re(b) \equiv \pm 1/4 \ (mod \Gamma)$.

Since $2 b \in  \Gamma$, we easily see that $b$ is either congruent to
$\tau/2$, or to $1/2 + \tau/2$. But by exchanging $b$ with $-a \bar{b} = -
\bar{b}$, we can assume $b = \tau/2$.

Thus we may assume $\sigma(z) = \bar{z} + \tau/2$, and $\sigma^2(z) = z -1/2
\equiv z + 1/2  \ (mod \ \Gamma)$.

In the case $\tau = \rho$, we observe that the two maps $z \mapsto \tau
\bar{z} + 1/2$ and $z \mapsto \bar{z} + \tau/2$ are conjugated by the
automorphism given by multiplication by $\rho$.
\hfill Q.E.D. \\

\begin{PROP}   In the non split case for the curve $E$, only the case $ Re
(\tau) =0 ,  |\tau| \geq  1$, $a = \pm 1$ occurs.

If $a =1$, $\sigma(z) = \bar{z} + b$, we  have
$\sigma^2(z) = z + 1/2$, $b \equiv 1/4 $.

If $a = -1$, we have
$\sigma(z) = - \bar{z} +b$, $\sigma^2(z) = \tau/2$, $b \equiv 1/4 \tau $.

\end{PROP}

{\bf Proof.} The fact that $\hat{G}$ is abelian implies that for every
$d \in 1/2 \Gamma - \Gamma$, we must have
$a \bar{d} \equiv  d \ (mod \ \Gamma)$, whence the map $d \mapsto a \bar{d}$
is the identity on  the points of two torsion of $E$.

If $|\tau| =1$, we know that $a = \pm \tau$, but then for $d = 1/2$ we have
$a \bar{d} = \pm \tau/2 \not \equiv 1/2 \ (mod \ \Gamma)$, a contradiction.

If $Re(\tau) = - 1/2$, then $a = \pm 1$ and if we take $d = \tau/2$,
$a \bar{d}  = \pm \bar{\tau}/2 = \mp (\tau + 1)/2 \not \equiv \tau/2 \ (mod \
\Gamma)$, again a contradiction.

Therefore we can assume $Re(\tau) = 0$, $a = \pm 1$.  In this case the map $d
\mapsto a \bar{d}$ is the identity on the points  of two torsion of $E$.

If $a =1$,  the condition
$\sigma^2 (z) = a \bar{b} + b = \bar{b} + b = 2 Re(b) \in 1/2 \Gamma - \Gamma$
implies $Re(b) \equiv \pm 1/4 $,  and
$\sigma^2(z) = z + 1/2$. By \ref{tra}, we can assume $b \equiv \pm 1/4$ and by
exchanging $\sigma$ with $\sigma^{-1}$ we can assume $b = 1/4$.

If $a =-1$,  the condition $\sigma^2 (z) = a \bar{b} + b =  - \bar{b} + b = 2i
Im(b) \in 1/2 \Gamma - \Gamma$ implies  $Im(b) \equiv
\pm 1/4 Im(\tau) $ and $\sigma^2(z) = z + \tau/2$, again we conclude by
\ref{tra} and by substituting $\sigma$ with $\sigma^{-1}$.

\hfill Q.E.D. \\

\begin{REM} Observe that, by remark \ref{liftnonempty}, 1, in the non split
case we must have $Fix(\sigma) = S(\R) = \emptyset$.
\end{REM}

\section{Topology of the real part of $S$}

In this section we want to  describe the topology of the fixed point locus of
the involution $\sigma$ acting on a real hyperelliptic surface $S$.\\

\begin{REM}  The fixed point locus of $\sigma$ can only be a disjoint union of
tori and Klein bottles.\\  In fact we have the Albanese map

$$\diagram  S \dto^{\alpha} \rto^{\sigma}    & S \dto^{\alpha}  \\
             A = {\Co}/\Lambda  \rto^{\bar{\sigma}} & A = {\Co}/\Lambda
                                       \enddiagram$$ and each component of the
real locus is a $S^1$ bundle on $S^1$.
\end{REM}

We have the following

\begin{LEM} Let $\sigma$ be as above with $Fix(\sigma) \neq \emptyset$ and
assume that the group $G$ is of odd order, i.e., either ${\Z}/3$ or ${\Z}/3
\times {\Z}/3$. Then the connected components of $Fix(\sigma)$ are
homeomorphic to
$S^1 \times S^1$.
\end{LEM} {\bf Proof.}

Let $\pi: E \times F \rightarrow (E \times F)/G$ be the projection.  We notice
(see \ref{liftnonempty}) that if $z \in Fix(\sigma)$, then $\forall
\tilde{z} \in
\pi^{-1}(z)$, $\exists
\ \tilde{\sigma}: E \times F \rightarrow E \times F$ such that
$\tilde{\sigma}(\tilde{z}) = \tilde{z}$.  We have the following commutative
diagram,
$$\diagram  E \times F  \dto^{\pi} \rto^{\tilde{\sigma}}   & E \times F
\dto^{\pi}  \\
             S  \rto^{\sigma} & S
                                       \enddiagram$$ where $\tilde{\sigma}$ is
an antiholomorphic involution on $E \times F$.

Let $C$ be the connected component of $Fix(\sigma)$ that contains $z$ and let
$\tilde{C}$ be the connected component of $\pi^{-1}(C)$ that contains
$\tilde{z}$.

Clearly $\tilde{\sigma} (\tilde{C}) = (\tilde{C})$, moreover
$\tilde{\sigma}$ is a lift of the identity of $C$ and with $\tilde{\sigma}
(\tilde{z}) = (\tilde{z})$, whence $\tilde{\sigma}$ is the identity on
$\tilde{C}$.

  We know that the components of
$Fix(\sigma)$ can be either tori or Klein bottles.

Since $\tilde{C} \cong S^1
\times S^1$ in the latter case $\tilde{C}$ would be an oriented covering of
odd degree (equal to $|G|$) of a Klein bottle. But this is impossible since
every oriented covering of a Klein bottle factors through the orientation
covering which has degree 2.
\hfill Q.E.D. \\

Using the Harnack - Thom - Krasnov inequality we are able to bound the number
of connected components of the real part of $S$ in all the cases given by the
list of Bagnera - de Franchis. We will show later that for all $G$ we can
find real hyperelliptic surfaces $S$ that reach these bounds (i.e., we can
find M - surfaces).
\begin{REM}

\begin{enumerate}

\item If $G = {\Z}/2$, then $h^0(S({\R}), {\Z}/2) \leq 4$;
\item If $G = {\Z}/2 \times {\Z}/2 $, then $h^0(S({\R}), {\Z}/2)
\leq 3$;
\item If $G = {\Z}/4$, then $h^0(S({\R}), {\Z}/2) \leq 3$;
\item If $G = {\Z}/4 \times {\Z}/2$, then $h^0(S({\R}), {\Z}/2) \leq 2$;
\item If $G = {\Z}/3$, then $h^0(S({\R}), {\Z}/2) \leq 2$;
\item If $G = {\Z}/3 \times {\Z}/3$, then $h^0(S({\R}), {\Z}/2) \leq 2$;
\item If $G = {\Z}/6$, then $h^0(S({\R}), {\Z}/2) \leq 2$.
\end{enumerate}
\end{REM}  {\bf Proof.}  We recall the Harnack-Thom-Krasnov  inequality:
$$\sum_m (dim_{{\Z}/2} H^m(S({\R}), {\Z}/2) ) \leq \sum_m (dim_{{\Z}/2} H^m(S,
{\Z}/2) -2 \lambda_m),$$  where $$\lambda_m := dim_{{\Z}/2} (1 + \sigma)
H^m(S, {\Z}/2).$$  Now we have
$$\sum_m (dim_{{\Z}/2} H^m(S({\R}), {\Z}/2)) = 2 h^0(S({\R}),{\Z}/2) +
h^1(S({\R}), {\Z}/2),$$  and since the components of
$S({\R})$ are either Tori or Klein bottles, we have $h^1(S({\R}), {\Z}/2) = 2
h^0(S({\R}), {\Z}/2)$. Thus we find
\begin{equation}
\label{htk} 4 h^0(S({\R}), {\Z}/2) \leq \sum_m (dim_{{\Z}/2} H^m(S,{\Z}/2)).
\end{equation}  We notice that $b_2(S) =2$, and we have the following list for
$H_1(S, {\Z})$ (see e.g. \cite{su}):
\begin{itemize}
\item $G = {\Z}/2,  \ H_1(S, {\Z}) = {\Z}^2 \oplus ({\Z}/2)^2$;
\item $G = {\Z}/2 \oplus {\Z}/2,  \ H_1(S, {\Z}) = {\Z}^2 \oplus {\Z}/2$;
\item $G = {\Z}/4,  \ H_1(S, {\Z}) = {\Z}^2 \oplus {\Z}/2$;
\item $G = {\Z}/4 \oplus {\Z}/2,  \ H_1(S, {\Z}) = {\Z}^2 $;
\item $G = {\Z}/3,  \ H_1(S, {\Z}) = {\Z}^2 \oplus {\Z}/3$;
\item $G = {\Z}/3 \oplus {\Z}/3,  \ H_1(S, {\Z}) = {\Z}^2$;
\item $G = {\Z}/6,  \ H_1(S, {\Z}) = {\Z}^2$.
\end{itemize}  The Universal Coefficients Theorem and Poincar\'e duality
allow us to compute all the Betti numbers of $S$ with coefficients in
${\Z}/2$. Then by
\ref{htk} we obtain the given bounds on the number of connected components of
$S({\R})$.
\hfill Q.E.D. \\

We want now to show how one can compute the number $k$ of Klein bottles,
respectively the number $t$ of 2-dimensional tori in the real part
$S({\R})$ of a real hyperelliptic surface.

For every connected component $V$ of $Fix(\sigma) = S({\R})$ the inverse image
$\pi^{-1}(V)$ splits as the $G$-orbit of any of its connected components, let
$W$ be one such.

We have already observed that there is a lift $\tilde{\sigma}$ of $\sigma$
such that
$W$ is in the fixed locus of $\tilde{\sigma}$.

If $\tilde{\sigma_1}$, $\tilde{\sigma_2}$ are two distinct antiholomorphic
involutions, then
$\pi^{-1}(V) \cap Fix(\tilde{\sigma_1}) \cap Fix(\tilde{\sigma_2}) =
\emptyset$, in fact otherwise there would exist a component $W_h$ such that
$\tilde{\sigma_1}_{|W_h} = \tilde{\sigma_2}_{|W_h} = Id$, thus
$\tilde{\sigma_1} \tilde{\sigma_2} = Id$, a contradiction.

We conclude from the above argument that the connected components of
$Fix(\sigma)$ correspond bijectively to the set ${\cal C}$ obtained as
follows: consider all the lifts $\tilde{\sigma}$ of $\sigma$ which are
involutions and pick one representative $\tilde{\sigma_i}$ for each conjugacy
class.

Then we let ${\cal C}$ be the set of equivalence classes of connected
components of
$\cup Fix(\tilde{\sigma_i})$, where two components $A$, $A'$ of
$Fix(\tilde{\sigma_i})$ are equivalent if and only if  there exists an
element $g
\in G$, such that $g(A) = A'$.

Let $\tilde{\sigma}$ be an antiholomorphic involution which is a lift of
$\sigma$ and such that $Fix(\tilde{\sigma}) \neq \emptyset$.

Since $\tilde{\sigma}$ is of product type, $Fix(\tilde{\sigma})$ is a disjoint
union of $2^{a_1 + a_2}$ copies of ${\bf S}^1 \times {\bf S}^1$, where $a_i
\in
\{0,1\}$.

In fact if an antiholomorphic involution $\hat{\sigma}$ on an elliptic curve
$C$ has fixed points, then $Fix(\hat{\sigma})$ is a disjoint union of $2^a$
copies of ${\bf S}^1$, where $a =1$ if the matrix of the action of
$\hat{\sigma}$ on $H_1(C, \Z)$ is diagonalizable, else $a =0$.

What said insofar describes the number of such components $V$:

in order to determine their nature, observe that $V = W/H$, where $W$ is as
before and $H \subset G$ is the subgroup such that $HW = W$.

Since the action of $G$ on the first curve $E$ is by translations, the action
on the first ${\bf S}^1$ is always orientation preserving, whence
$V$ is a Klein bottle if and only if $H$ acts on the second ${\bf S}^1$ by
some orientation reversing map, or equivalently $H$ has some fixed point on
the second
${\bf S}^1$.

Let $h \in H$ be a transformation having a fixed point on the second ${\bf
S}^1$: since the direction of this ${\bf S}^1$ is an eigenvector for the
tangent action of $h$, it follows that the tangent action is given by
multiplication by
$-1$.

The existence of such an element $h$ for the second ${\bf S}^1$ is obvious if
$a_1 = a_2 =0$ and $\tilde{\sigma}$ commutes with $h$ (this always occurs
except in the last case of $2$ of \ref{action}) and in this case we get
$1$ Klein bottle.

In the other cases one needs a more delicate analysis whose details we omit
here.

In the next section we shall give a precise description of all the isomorphism
classes of  real hyperelliptic surfaces and we shall also determine the real
part
$S({\R})$ for each isomorphism class.

\section{Moduli space of real hyperelliptic surfaces}

In this section we describe all the isomorphism classes of of real
hyperelliptic surfaces. Each slot in the following tables corresponds to a
fixed topological type.

In the column of the values of $\sigma_1$, we always have two possible
values, that correspond to two antiholomorphic maps which are topologically
equivalent, but not analytically equivalent. Nevertheless they are in the
same connected component of the moduli space of real elliptic curves, since
the curve with $\tau_1 =i$ has an isomorphism given by the multiplication by
$i$ that conjugates the two different antiholomorphic maps.

Assume first of all that we are in the split case for the extension $(*)$.

We set
$E = {\Co}/{\Z}+ \tau_1 {\Z}$, $F = {\Co}/{\Z}+ \tau_2 {\Z}$,
$\tau_j = x_j + i y_j$,
$\tilde{\sigma} = (\sigma_1, \sigma_2): E \times F \rightarrow E \times F$.
The action of $G$ is given as in the Bagnera - de Franchis list.

The first case that we will consider is the one with
$G = \Z/2$, $\hat{G} = \Z/2 \times \Z/2$. Let
$g$ be a generator of
$G$: its action  on $E$ is a translation by an element $\eta \in Pic^0(E)_2$,
thus
$g(z_1,z_2) = (z_1 + \eta, -z_2)$, $\forall (z_1,z_2) \in E \times F$.

Let us now give a list of the topological invariants that distinguish the
different cases presented in the forthcoming table.

Recall that for an antiholomorphic involution $\sigma$ on an elliptic curve
$C$, the number $\nu(\sigma)$ of connected components of $Fix(\sigma)$, which
is either
$0$, or $1$ or $2$, is a topological invariant of the involution.

1) In our case, given  an antiholomorphic involution on a hyperelliptic
surface
$S$, since we have exactly two lifts of $\sigma$ on $E \times F$,
$\tilde{\sigma} =  (\sigma_1,
\sigma_2)$, $\tilde{\sigma} \circ g$ the sets

1)(a): $\{\nu(\sigma_1), \nu(\sigma_1 \circ g)\}$,

1)(b): $\{\nu(\sigma_2), \nu(\sigma_2 \circ g)\}$,

are topological invariants of $\sigma$.

Since $G$ acts on $E$ by translation and on $F$ by multiplication by
$-1$, also the parity of $\nu(\sigma_i)$ (determined by the alternative: the
action on the first homology group is diagonalizable or not) is a topological
invariant.

Other topological invariants are:

2) the homology of $S(\R)$;

3) the action of the involution $\sigma_2$ on the fixed point locus of the
action of
$g$ on $F$ (see \ref{fixggg});

4) the action of $\hat{G'}$ on $\Lambda'' \otimes \R$, where $\Lambda'' =
\pi_1(E)
\cong {\Z}^2$, given by the orbifold fundamental group exact sequence:

$1 \rightarrow \Lambda'' \oplus \Gamma \rightarrow \hat{\Pi}  \rightarrow
\hat{G}
\rightarrow 1$,

where $\Gamma = \pi_1(F)$ (see section  2).

In particular, given a transformation in $G$, we have the condition whether
its first action is representable by a translation in $\Lambda''
\otimes \R$ which is an eigenvector for $\sigma_1$.

$$\begin{tabular}{|c|c|c|c|c|c|}
\hline
$G$, $\hat{G}$ & $\tau_1$ & $\tau_2$ & $\sigma_1(z_1), \ \eta$ &
$\sigma_2(z_2)$ & $S({\R})$  \\
\hline
\hline
${\Z}/2$ &  $x_1 =0$ & $x_2 =0$ & $\bar{z_1}, \eta = \frac{1}{2}$ &
$\bar{z_2}$ & $4K$ \\
$({\Z}/2)^2$ & & &  or $-\bar{z_1},
\eta = \frac{\tau_1}{2}$ &  & \\
\hline
${\Z}/2$ & $|\tau_1| = 1,$  & $x_2 =0$  & $\pm \tau_1 \bar{z_1}$,
$\eta = \frac{1 + \tau_1}{2}$,  & $\bar{z_2}$ & $4K$\\
$({\Z}/2)^2$ & &  & if
$|\tau_1| =1$, &  & \\

& $\cup \ x_1 = -\frac{1}{2}$, & & $\pm \bar{z_1}, \ \eta = \frac{1}{2}$ & &
\\

& & &  if $x_1 = -\frac{1}{2}$ & & \\
\hline

${\Z}/2$ & $x_1 =0$, & $x_2 = 0$, & $\bar{z_1}$, $\eta = \frac{1}{2}$ &
$\bar{z_2}  + \frac{\tau_2}{2}$ & $2T$\\

$({\Z}/2)^2$ &  & & or $-
\bar{z_1}$, $\eta =
\frac{\tau_1}{2}$ & & \\

\hline

${\Z}/2$ & $x_1 =0$, & $x_2 = 0$, & $\bar{z_1}$, $\eta =
\frac{\tau_1}{2}$ &
$\bar{z_2} + \frac{\tau_2}{2}$ & $2T$\\

  $({\Z}/2)^2$ &  & & or $-
\bar{z_1}$, $\eta =
\frac{1}{2}$ &  & \\
\hline

${\Z}/2$ & $x_1 =0$, & $| \tau_2| = 1$, & $\bar{z_1}$, $\eta =
\frac{\tau_1}{2}$ & $\tau_2 \bar{z_2}$ if & $2T$\\

  $({\Z}/2)^2$ & &or & or $- \bar{z_1}$,
$\eta = \frac{1}{2}$ & $|\tau_2| = 1$& \\   & & $\ x_2 = - \frac{1}{2}$  &
&$\bar{z_2}$ if & \\    & & & & $x_2 = - \frac{1}{2}$ & \\  & & & & & \\
\hline

${\Z}/2$ & $x_1 = 0$, & $x_2 =0$, & $\bar{z_1}$, $\eta = \frac{1+
\tau_1}{2}$ & $\bar{z_2}$ & $2T$\\
  $({\Z}/2)^2$ & & & or $- \bar{z_1}$, $\eta = \frac{1 + \tau_1}{2}$ & & \\
\hline

${\Z}/2$ & $x_1 = 0$, & $x_2 =0$, & $\bar{z_1}$, $\eta = \frac{1+
\tau_1}{2}$ & $\bar{z_2} + \frac{\tau_2}{2}$ & $2T$\\
$({\Z}/2)^2$  & & & or $- \bar{z_1}$,
$\eta = \frac{1 + \tau_1}{2}$ & & \\
\hline

${\Z}/2$ & $x_1 = 0$, & $x_2 =0$, & $\bar{z_1}$, $\eta = \frac{1}{2}$ &
$\bar{z_2} +  \frac{1}{2}$ & $\emptyset$\\
  $({\Z}/2)^2$ & & &or $- \bar{z_1}$, $\eta = \frac{\tau_1}{2}$  & & \\
\hline

${\Z}/2$ & $x_1 = 0$, & $x_2 =0$, & $\bar{z_1}$, $\eta = \frac{1}{2}$ &
$\bar{z_2} +  \frac{1+ \tau_2}{2}$ & $\emptyset$\\
$({\Z}/2)^2$ & & & or $- \bar{z_1}$, $\eta = \frac{\tau_1}{2}$   & & \\
\hline

${\Z}/2$ & $x_1 = 0$, & $x_2 =0$, & $\bar{z_1}$, $\eta =
\frac{\tau_1}{2}$ &
$\bar{z_2} + \frac{1+ \tau_2}{2}$ & $\emptyset$\\
$({\Z}/2)^2$ & & & or $- \bar{z_1} $, $ \eta = \frac{1}{2}$ & & \\
\hline

${\Z}/2$ & $x_1 = 0$, & $x_2 =0$, & $\bar{z_1}$, $\eta = \frac{1 +\tau_1}{2}$
&
$\bar{z_2} + \frac{1}{2}$ & $\emptyset$\\
  $({\Z}/2)^2$ & & & or $- \bar{z_1}$,
$\eta = \frac{1 + \tau_1}{2}$ & & \\
\hline

${\Z}/2$ & $x_1 = 0$, & $x_2 =0$, & $\bar{z_1}$, $\eta = \frac{1 +\tau_1}{2}$
&
$\bar{z_2} + \frac{1 + \tau_2}{2}$ & $\emptyset$\\
$({\Z}/2)^2$  & & &or $- \bar{z_1}$,
$\eta = \frac{1 + \tau_1}{2}$ & & \\
\hline

${\Z}/2$ & $|\tau_1| = 1$, & $x_2 =0$, & $\pm \tau_1 \bar{z_1}$, $\eta =
\frac{1 +\tau_1}{2}$ & $\bar{z_2} + \frac{1 + \tau_2}{2}$ & $\emptyset$\\
$({\Z}/2)^2$ & $\cup \ x_1 = -
\frac{1}{2}$ & & if
$|\tau_1| = 1$ & & \\  &$|\tau_1| \geq 1$ & & $\pm \bar{z_1},  \ \eta =
\frac{1}{2}$ & & \\   & & & if $x_1 = - \frac{1}{2}$ & & \\
\hline
\end{tabular}
$$

$$\begin{tabular}{|c|c|c|c|c|c|}
\hline

${\Z}/2$ & $x_1 = 0$, & $x_2 =0$, & $\bar{z_1} + \frac{1}{2}$, $\eta =
\frac{\tau_1}{2}$ & $\bar{z_2}$ & $\emptyset$\\
$({\Z}/2)^2$  & & &or $- \bar{z_1} +
\frac{\tau_1}{2}$,
$\eta = \frac{1}{2}$ & & \\
\hline

${\Z}/2$ & $x_1 = 0$, & $x_2 =0$, & $\bar{z_1} + \frac{1}{2}$, $\eta =
\frac{\tau_1}{2}$ & $\bar{z_2} + \frac{\tau_2}{2}$ & $\emptyset$\\
$({\Z}/2)^2$  & & &or $- \bar{z_1} +
\frac{\tau_1}{2}$,
$\eta = \frac{1}{2}$ & & \\
\hline

${\Z}/2$ & $x_1 = 0$, & $x_2 =0$, & $\bar{z_1} + \frac{1}{2}$, $\eta =
\frac{\tau_1}{2}$ & $\bar{z_2} + \frac{1 + \tau_2}{2}$ & $\emptyset$\\
$({\Z}/2)^2$  & & &or $- \bar{z_1} +
\frac{\tau_1}{2}$,
$\eta = \frac{1}{2}$ & & \\
\hline

${\Z}/2$ & $x_1 =0$, & $|\tau_2| =1$, & $\bar{z_1} + \frac{1}{2}$, $\eta =
\frac{\tau_1}{2}$ & $\tau_2 \bar{z_2}$   & $\emptyset$ \\
$({\Z}/2)^2$ &  &$ \cup \ x_2 = - \frac{1}{2}$,      & or $- \bar{z_1} +
\frac{\tau_1}{2}$, $ \eta =
\frac{1}{2}$ & if
$|\tau_2| =1$, & \\

& & & & $\bar{z_2}$ & \\  & & & & if $x_2 = - \frac{1}{2}$ & \\
\hline

${\Z}/2$ & $x_1 =0$, & $|\tau_2| =1$ & $\bar{z_1}$, $\eta = \frac{1}{2}$ &
$\tau_2 \bar{z_2}$, & $2K$ \\
$({\Z}/2)^2$ &  &$ \cup \ x_2 = - \frac{1}{2}$    & or $- \bar{z_1}$, $
\eta =
\frac{\tau_1}{2}$ & if
$|\tau_2| =1$, & \\

& & & & $ \bar{z_2}$, & \\  & & & & if $x_2 = - \frac{1}{2}$ & \\
\hline

${\Z}/2$ & $|\tau_1| =1$, & $|\tau_2| =1$, & $\pm \tau_1 \bar{z_1}$,
$\eta = \frac{1 + \tau_1}{2}$ &  $\tau_2 \bar{z_2}$, & $2K$ \\
$({\Z}/2)^2$  &$\cup \ x_1 = - \frac{1}{2}$  &$ \cup \ x_2 = -
\frac{1}{2}$        & if $|\tau_1| =1$ & if $|\tau_2| =1$, & \\

& &  & $\pm
\bar{z_1}$, $  \eta = \frac{1}{2}$ & $\bar{z_2}$, & \\  &  & & if $x_1 = -
\frac{1}{2}$ & if $x_2 = -
\frac{1}{2}$ & \\
\hline

${\Z}/2$ & $x_1 =0$, & $x_2  =0$, & $\bar{z_1}$, $\eta =
\frac{\tau_1}{2}$ &
  $\bar{z_2}$, & $4T$ \\
$({\Z}/2)^2$  & & &or $-
\bar{z_1}$, $\eta = \frac{1}{2}$ & & \\

\hline

${\Z}/2$ & $x_1 =0$, & $|\tau_2| =1$, & $\bar{z_1}$, $\eta = \frac{1 +
\tau_1}{2}$ &  $\tau_2 \bar{z_2}$, & $T$ \\
  $({\Z}/2)^2$ &  &$ \cup \ x_2 = - \frac{1}{2}$,       & or $- \bar{z_1}
$,
$\eta =
\frac{1 +
\tau_1}{2}$ & if $|\tau_2| =1$, & \\

& & & & $\bar{z_2}$, & \\  &  & & & if $x_2 = - \frac{1}{2}$ & \\
\hline

${\Z}/2$ & $|\tau_1| =1$, & $x_2 =0$ & $\pm \tau_1 \bar{z_1}$, $\eta =
\frac{1 + \tau_1}{2}$ &  $\bar{z_2} + \frac{\tau_2}{2}$, & $T$ \\
$({\Z}/2)^2$  & $\cup \ x_1 = - \frac{1}{2}$ & & if $|\tau_1| =1$ & & \\

& & & $\pm \bar{z_1}$, $ \eta = \frac{1}{2}$ & &\\  & & & if $x_1 = -
\frac{1}{2}$ &  & \\
\hline

\end{tabular}
$$

Let us now explain how to distinguish the different cases listed in the
  above table:

\bigskip

$S(\R) = 4K$: the first two cases of the list are distinguished by the parity
of
$\nu(\sigma_1)$ ($ \nu = 2$, respectively $ \nu= 1$).

\bigskip

$S(\R) = 2T$: the third case listed with $S(\R) = 2T$ is distinguished from
all the others by the parity of $\nu(\sigma_2) = 1$.

The fourth case is distinguished from all the others by the invariant 3), or
by 1)(b) (the set of values of $\nu(\sigma_2)$ equals $\{ 2,2 \}$).

The first is distinguished from the second by the invariant 1)(a) (the
respective sets are $\{ 2,0 \}$, $\{ 2,2 \}$)) and from the fifth by the
invariant 4): in fact the action of $g$ on
$\Lambda''\otimes \R$ is given by $z_1
\mapsto z_1 + 1/2$ in the first case, by $z_1 \mapsto z_1 + (1+ \tau_1)/2$ in
the fifth case, $\sigma_1(z_1) = \bar{z_1}$ in both cases, but $1/2$ is an
eigenvector of the action of $\sigma_1$ on $\Lambda''$, while $(1 +
\tau_1)/2$ is not.

The second is distinguished from the fifth by the invariant 1)(a) (the
respective sets are $\{ 2,2 \}$, $\{ 2,0 \}$)).

\bigskip

$S(\R) = \emptyset$: the sixth is distinguished from all  others by the
parity of
$\nu(\sigma_1) = 1$.

The last one is distinguished from all  others by the parity of
$\nu(\sigma_2) = 1$.

The seventh is distinguished from all  others by the invariant 3), or by
1)(b) (the set of values of $\nu(\sigma_2)$ equals $\{2,2\}$, while in all
the other cases it equals either $\{0,2\}$, or $\{0,0\}$, or
$\{1,1\}$ in the last case).

The first case is distinguished from the second, the third, the fifth, the
ninth case by the invariant 1)(b): in the first case the set of values of
$\nu(\sigma_2)$ equals $\{0,2\}$, in the others $\{0,0\}$.

The first case is distinguished from the the eighth case by the invariant
1)(a):  in the first case the set of values of $\nu(\sigma_1)$ equals
$\{2,0\}$, in the eight case it equals $\{0,0\}$.

The first case is distinguished from the fourth by the invariant 4):
  in fact the action of $g$ on
$\Lambda''\otimes \R$ is given by $z_1
\mapsto z_1 + 1/2$ in the first case, by $z_1 \mapsto z_1 + (1+ \tau_1)/2$ in
the fourth case, $\sigma_1(z_1) = \bar{z_1}$ in both cases, but $1/2$ is an
eigenvector of the action of $\sigma_1$ on $\Lambda''$, while $(1 +
\tau_1)/2$ is not.

The second case is distinguished from the third by the invariant 4):
  in fact the action of $g$ on
$\Lambda''\otimes \R$ is given by $z_1
\mapsto z_1 + 1/2$ in the second case, by $z_1 \mapsto z_1 + \tau_1/2$ in the
third case, $\sigma_1(z_1) = \bar{z_1}$ in both cases, but $1/2$ is the $+1$
eigenvector of the action of $\sigma_1$ on $\Lambda''$, while
$\tau_1/2$ is the $-1$ eigenvector.

The second case is distinguished from the fourth and the eighth case by the
invariant 1)(b): in the second case the set of values of
$\nu(\sigma_2)$ equals $\{0,0\}$, in the other cases $\{0,2\}$.

The second case is distinguished from the fifth by the invariant 4):
  in fact the action of $g$ on
$\Lambda''\otimes \R$ is given by $z_1
\mapsto z_1 + 1/2$ in the second case, by $z_1 \mapsto z_1 + (1+
\tau_1)/2$ in the fifth case, $\sigma_1(z_1) = \bar{z_1}$ in both cases, but
$1/2$ is an eigenvector of the action of $\sigma_1$ on $\Lambda''$, while
$(1+ \tau_1)/2$ is not.

The second case is distinguished from the the ninth case by the invariant
1)(a):  in the second case the set of values of $\nu(\sigma_1)$ equals
$\{2,0\}$, in the ninth case it equals $\{0,0\}$.

The third case is distinguished from the fourth by the invariant 1)(b): in
the third case the set of values of $\nu(\sigma_2)$ equals $\{0,0\}$, in the
fourth case $\{0,2\}$.

The third case is distinguished from the fifth by the invariant 4):
  in fact the action of $g$ on
$\Lambda''\otimes \R$ is given by $z_1
\mapsto z_1 + \tau_1/2$ in the third case, by $z_1 \mapsto z_1 + (1+
\tau_1)/2$ in the fifth case, $\sigma_1(z_1) = \bar{z_1}$ in both cases, but
$\tau_1/2$ is an eigenvector of the action of $\sigma_1$ on
$\Lambda''$, while $(1+ \tau_1)/2$ is not.

The third case is distinguished from the eighth and the ninth case by the
invariant 1)(a):  in the third case the set of values of $\nu(\sigma_1)$
equals $\{2,2\}$, in the eighth and the ninth case it equals $\{0,0\}$.

The fourth case is distinguished from the fifth and the ninth case by the
invariant 1)(b): in the fourth case the set of values of $\nu(\sigma_2)$
equals $\{0,2\}$, in the other cases $\{0,0\}$.

The fourth case is distinguished from the eighth case by the invariant
1)(a):  in the fourth case the set of values of $\nu(\sigma_1)$ equals
$\{0,2\}$, in the eighth case it equals $\{0,0\}$.

The fifth case is distinguished from the eighth case by the invariant 1)(b):
in the fifth case the set of values of $\nu(\sigma_2)$ equals
$\{0,0\}$, in the eighth case $\{2,0\}$.

The fifth case is distinguished from the ninth case by the invariant 1)(a):
in the fifth case the set of values of $\nu(\sigma_1)$ equals
$\{2,0\}$, in the ninth case it equals $\{0,0\}$.

The eighth case is distinguished from the ninth case by the invariant 1)(b):
in the eighth case the set of values of $\nu(\sigma_2)$ equals
$\{2,0\}$, in the ninth case $\{0,0\}$.

\bigskip

$S(\R) = 2K$: the parity of $\nu(\sigma_1)$ distinguishes the two cases.

\bigskip

$S(\R) = T$: the parity of $\nu(\sigma_2)$ distinguishes the two cases.

\bigskip

Let us now consider the case $G = \Z/4$. If $g$ is a generator of $G$, we
denote by
$g(z_1, z_2) = (z_1 + \eta, iz_2)$ the action of $g$ on $E \times F$.

The topological invariants that distinguish the different cases here are the
above mentioned invariant 2), and the analogous of the invariants 1)(a) and
1)(b) that we had in the case
$G = \Z/2$:

$1'$)(a) the set $\{\nu(\sigma_1 \circ g^n ), \ for \ n = 0,1,2,3 \}$.

$1'$)(b) the set $\{\nu(\sigma_2 \circ g^n ), \ for \ n = 0,1,2,3 \}$.

$$\begin{tabular}{|c|c|c|c|c|c|}
\hline
${\Z}/4$ & $x_1 =0$, & $i$ & $\bar{z_1}$, $\eta = \frac{1}{2} +
\frac{\tau_1}{4}$ &  $\bar{z_2}$, & $2T$ \\
$D_4$ & &  & or $- \bar{z_1}$, $\eta = \frac{1}{4} +
\frac{\tau_1}{2}$ & & \\

\hline

${\Z}/4$ & $x_1 =0$, & $i$ & $\bar{z_1} + \frac{1}{2}$, $\eta =
\frac{1}{2} +
\frac{\tau_1}{4}$ &  $ \bar{z_2}$, & $T$ \\
$D_4$ & &  & or
$- \bar{z_1} + \frac{\tau_1}{2}$, $\eta = \frac{1}{4} + \frac{\tau_1}{2}$ & &
\\
\hline

${\Z}/4$ & $x_1 =0$, & $i$ & $\bar{z_1}$, $\eta =\frac{\tau_1}{4}$ &
$\bar{z_2} + \frac{1 + i}{2}$, & $T$ \\
$D_4$ & &  & or $- \bar{z_1}$, $\eta = \frac{1}{4}$ & & \\
\hline

${\Z}/4$ & $x_1 =0$, & $i$ & $\bar{z_1} + \frac{1}{2}$, $\eta =\frac{1}{2} +
\frac{\tau_1}{4}$ &  $\bar{z_2} + \frac{1+i}{2} $, & $T$ \\
$D_4$ & & & or $- \bar{z_1} + \frac{\tau_1}{2}$, $\eta =
\frac{1}{4} + \frac{\tau_1}{2}$ & & \\
\hline

${\Z}/4$ & $x_1 =0$, & $i$ & $\bar{z_1}$, $\eta =\frac{1}{2} +
\frac{\tau_1}{4}$ &  $\bar{z_2} + \frac{1+i}{2} $, & $\emptyset$ \\
$D_4$ & & & or $- \bar{z_1}$, $\eta = \frac{1}{4} +
\frac{\tau_1}{2}$ & & \\
\hline

${\Z}/4$ & $x_1 =0$, & $i$ & $\bar{z_1} + \frac{1}{2}$, $\eta
=\frac{\tau_1}{4}$ &
$\bar{z_2}$, & $\emptyset$ \\
  $D_4$ & & & or $- \bar{z_1} + \frac{\tau_1}{2}$, $\eta =
\frac{1}{4}$ & & \\
\hline

${\Z}/4$ & $x_1 =0$, & $i$ & $\bar{z_1} + \frac{1}{2}$, $\eta
=\frac{\tau_1}{4}$ &
$\bar{z_2} + \frac{1+i}{2}$, & $\emptyset$ \\
  $D_4$ & & & or $- \bar{z_1} + \frac{\tau_1}{2}$, $\eta =
\frac{1}{4}$  & & \\
\hline

${\Z}/4$ & $x_1 =0$, & $i$ & $\bar{z_1}$, $\eta =\frac{\tau_1}{4}$ &
$\bar{z_2}$, & $3T$ \\
  $D_4$ &  & & or $- \bar{z_1}$, $\eta =
\frac{1}{4}$  & & \\
\hline

${\Z}/4$ & $|\tau_1|=1$ & $i$ & $\tau_1 \bar{z_1}$, $\eta =\frac{1-
\tau_1}{4}$ &
$\bar{z_2}$, & $3K$ \\
$D_4$ &  $ \cup \ x_1 =-\frac{1}{2}$ & &or $-\tau_1 \bar{z_1}$, $\eta =
\frac{1+
\tau_1}{4}$ & & \\   &  & & if $|\tau_1| = 1$, & & \\

& & & $\bar{z_1}$, $\eta = \frac{1}{4} +
\frac{\tau_1}{2}$ & & \\ &  & &  $- \bar{z_1}$, $\eta = \frac{1}{4}$  & &
\\ &  & &  if $ x_1 = - \frac{1}{2}$ & & \\

\hline

${\Z}/4$ & $|\tau_1|=1$ & $i$ & $\tau_1 \bar{z_1}$, $\eta =\frac{1-
\tau_1}{4}$ &  $\bar{z_2} + \frac{1+i}{2}$ & $K$ \\
$D_4$ &  $ \cup \ x_1 =-\frac{1}{2}$& &or $-\tau_1 \bar{z_1}$, $\eta =
\frac{1+ \tau_1}{4}$ & & \\  &  & & if $|\tau_1| = 1$, & & \\

& & & $\bar{z_1}$, $\eta = \frac{1}{4} +
\frac{\tau_1}{2}$ & & \\ &  & &  $- \bar{z_1}$, $\eta = \frac{1}{4}$  & &
\\ &  & &  if $ x_1 = - \frac{1}{2}$ & & \\

\hline
\end{tabular}
$$

In the list, if $S(\R) = T$, the second and the third cases are distinguished
by invariant $1'$)(a):  in the second case the set of values of
$\nu(\sigma_1)$ equals $\{2,2, 2,2\}$, in the third it equals
$\{0,2,0,2\}$.

The first case is distinguished from the others by the invariant
$1')$(b): in the first case $\{\nu(\sigma_2 \circ g^n ), \ for \ n = 0,1,2,3
\} =
\{2,1,2,1\}$, while in the other cases $\{\nu(\sigma_2 \circ g^n ), \ for
\ n = 0,1,2,3 \} = \{0,1,0,1\}$.

If $S(\R) = \emptyset$, the first case is distinguished from the third case
by the invariant $1'$)(a):  in the first case the set of values of
$\nu(\sigma_1)$ equals $\{2,0, 2,0\}$, in the third it equals
$\{0,0,0,0\}$.

The second case is distinguished from the other cases  by the invariant
$1'$)(b):  in the second case the set of values of $\nu(\sigma_2)$ equals
$\{2,1, 2,1\}$, in the other cases it equals $\{0,1,0,1\}$.
\bigskip

If $G = \Z/4 \times \Z/2$, we take generators $g$, $t$, such that on $F$,
$g(z_2) = i z_2$, $t(z_2) = z_2 + (1+i)/2$.

If $\hat{G} \cong G_1$, we see that the topological type of $S(\R)$
distinguishes the cases. The same holds if
$\hat{G} \cong D_4 \times \Z/2$.

$$
\begin{tabular}{|c|c|c|c|c|c|}
\hline

${\Z}/4 \times {\Z}/2$ & $x_1 =0$ & $i$ & $\bar{z_1}$, &
$\bar{z_2} + \frac{1}{2}$, & $\emptyset$ \\
$G_1$ &  & & $G = <\frac{1 + \tau_1}{4}>
\times <\frac{1}{2}>$&$\epsilon = \frac{1+i}{2}$ & \\   & & & or $-
\bar{z_1}$ & & \\  & & & $G = <\frac{-1 + \tau_1}{4}> \times
<\frac{\tau_1}{2}>$ & & \\
\hline

${\Z}/4 \times {\Z}/2$ & $x_1 =0$ & $i$ & $\bar{z_1}+ \frac{1}{2}$, &
$\bar{z_2}+ \frac{1}{2}$, & $T$ \\
  $G_1$ & & & $G = <\frac{1+\tau_1}{4}>
\times <\frac{1}{2}>$&$\epsilon = \frac{1+i}{2}$ & \\  & & & or $-
\bar{z_1} +
\frac{\tau_1}{2}$ & & \\   &  & & $G = <\frac{-1 + \tau_1}{4}> \times
<\frac{\tau_1}{2}>$ & & \\
\hline
\end{tabular}
$$

$$
\begin{tabular}{|c|c|c|c|c|c|}
\hline

${\Z}/4 \times {\Z}/2$ & $x_1 =0$ & $i$ & $\bar{z_1}$, &
$\bar{z_2}$, & $2T$ \\
$D_4 \times {\Z}/2$ &  & & $G = <\frac{\tau_1}{4}>
\times <\frac{1}{2}>$&$\epsilon = \frac{1+i}{2}$ & \\   & & & or $-
\bar{z_1}$ & & \\  & & & $G = <\frac{1}{4}> \times <\frac{\tau_1}{2}>$ & & \\
\hline

${\Z}/4 \times {\Z}/2$ & $x_1 =0$ & $i$ & $\bar{z_1}$, &
$\bar{z_2}$, & $K \sqcup T$ \\
  $D_4 \times {\Z}/2$ & & & $G = <\frac{\tau_1}{4}>
\times <\frac{1+ \tau_1}{2}>$&$\epsilon = \frac{1+i}{2}$ & \\  & & & or $-
\bar{z_1}$ & & \\   &  & & $G = <\frac{1}{4}> \times <\frac{1+
\tau_1}{2}>$ & & \\

\hline

${\Z}/4 \times {\Z}/2$ & $x_1 =0$, & $i$ & $\bar{z_1} + \frac{1}{2}$, &
$\bar{z_2}$, & $T$ \\
  $D_4 \times {\Z}/2$ & & & $G = <\frac{\tau_1}{4}>
\times <\frac{1}{2}>$& $\epsilon = \frac{1+i}{2}$ & \\   &  & & or $-
\bar{z_1} + \frac{\tau_1}{2}$ & & \\ &  & & $G = <\frac{1}{4}> \times
<\frac{\tau_1}{2}>$ & & \\
\hline

${\Z}/4 \times {\Z}/2$ & $x_1 =0$, & $i$ & $\bar{z_1} + \frac{1}{2}$, &  $
\bar{z_2}$, & $K$ \\
$D_4 \times {\Z}/2$ &  & & $G = <\frac{\tau_1}{4}>
\times <\frac{1+ \tau_1}{2}>$& $\epsilon = \frac{1+i}{2}$& \\   & & & or
$-
\bar{z_1} + \frac{\tau_1}{2}$ & & \\   &  & & $G = <\frac{1}{4}> \times
<\frac{1+
\tau_1}{2}>$& & \\
\hline
\end{tabular}
$$

If $G = \Z/3$, we denote by $g(z_1, z_2) =( z_1 + \eta, \rho z_2)$ the action
of a generator $g$ of $G$ on $E \times F$.

We notice that $Fix(g) = \Z (1 - \rho) /3 \Z (1 - \rho)$ and the action of
$\sigma_2(z_2) =
\bar{z_2}$ on $Fix(g)$ is $-Id$, while the action of $\sigma_2(z_2) =
-\bar{z_2}$ on
$Fix(g)$ is $Id$, so these two actions are topologically different.

This invariant and $S(\R)$ are sufficient to distinguish all the cases.

$$
\begin{tabular}{|c|c|c|c|c|c|}
\hline
${\Z}/3$ & $x_1 =0$, & $\rho$ & $\bar{z_1}$, $\eta = \frac{\tau_1}{3}$ &
$\bar{z_2}$, & $2T$ \\
  ${\cal S}_3$ & & & or $- \bar{z_1}$, $\eta = \frac{1}{3}$ & & \\
\hline

${\Z}/3$ & $x_1 =0$, & $\rho$ & $\bar{z_1}$, $\eta = \frac{\tau_1}{3}$ &
$-\bar{z_2}$, & $2T$ \\
  ${\cal S}_3$ & & & or $- \bar{z_1}$, $\eta =
\frac{1}{3}$& & \\

\hline

${\Z}/3$ & $x_1 =0$, & $\rho$ & $\bar{z_1} + \frac{1}{2}$, $\eta =
\frac{\tau_1}{3}$ &  $\bar{z_2}$, & $\emptyset$ \\
${\cal S}_3$ &  & & or $- \bar{z_1} + \frac{\tau_1}{2}$,
$\eta =
\frac{1}{3}$ & & \\
\hline

${\Z}/3$ & $x_1 =0$, & $\rho$ & $\bar{z_1} + \frac{1}{2}$, $\eta =
\frac{\tau_1}{3}$ &  $-\bar{z_2}$, & $\emptyset$ \\
  ${\cal S}_3$ & & &or $- \bar{z_1} + \frac{\tau_1}{2}$,
$\eta = \frac{1}{3}$ & & \\
\hline

${\Z}/3$ & $|\tau_1| =1$, & $\rho$ & $\tau_1 \bar{z_1}$, $\eta = \frac{1 -
\tau_1}{3}$ &  $\bar{z_2}$ & $T$ \\
  ${\cal S}_3$ & $\cup \ x_1 = - \frac{1}{2}$       & & or $- \tau_1
\bar{z_1}$, $
\eta = \frac{1+ \tau_1}{3}$&  & \\

&   & & if $|\tau_1| = 1$ & & \\

&  & &  $\bar{z_1}$, $\eta = \frac{1 - \tau_1}{3}$ & & \\

&   & & or $-\bar{z_1}$, $\eta = \frac{1}{3}$ & & \\ &   & &  if $ x_1 = -
\frac{1}{2}$ & & \\

\hline

${\Z}/3$ & $|\tau_1| =1$, & $\rho$ & $\tau_1 \bar{z_1}$, $\eta = \frac{1 -
\tau_1}{3}$ &  $-\bar{z_2}$, & $T$ \\
  ${\cal S}_3$ &  $\cup \ x_1 = - \frac{1}{2}$        & & or $- \tau_1
\bar{z_1}$,
$
\eta = \frac{1+ \tau_1}{3}$ &  & \\

&  & & if $|\tau_1| = 1$ & & \\

&  & &  $\bar{z_1}$, $\eta = \frac{1 - \tau_1}{3}$ & & \\

&   & &  or $-\bar{z_1}$, $\eta = \frac{1}{3}$ & & \\  &   & &  if $ x_1 = -
\frac{1}{2}$ & & \\

\hline
\end{tabular}
$$

If $G = \Z/3 \times \Z/3$, with generators $g$, $t$ acting on $F$ by
$g(z_2) =
\rho z_2$, $t(z_2) = z_2 + (1 - \rho)/3$, the homology of $S(\R)$ is the only
topological invariant needed in order to distinguish the three cases.

$$
\begin{tabular}{|c|c|c|c|c|c|}
\hline
${\Z}/3 \times {\Z}/3$ & $x_1 =0$, & $\rho$ & $\bar{z_1}$,&  $-\bar{z_2}$, &
$2T$
\\
${\cal S}_3 \times {\Z}/3$ & & & $G = <\frac{\tau_1}{3}>
\times <\frac{1}{3}>$ &$\epsilon = \frac{1-\rho}{3}$ & \\   &  & &or $-
\bar{z_1}$ & &
\\   &  & & $G = <\frac{1}{3}> \times <\frac{\tau_1}{3}>$ & & \\

\hline

${\Z}/3 \times {\Z}/3$ & $x_1 =0$, & $\rho$ & $\bar{z_1}+
\frac{1}{2}$,&  $-\bar{z_2}$, & $\emptyset$ \\
${\cal S}_3 \times {\Z}/3$ & & & $G = <\frac{\tau_1}{3}> \times
<\frac{1}{3}>$ &
$\epsilon = \frac{1-\rho}{3}$    & \\   & & & or $ - \bar{z_1} +
\frac{\tau_1}{2}$ & &
\\    &  & & $G = <\frac{1}{3}> \times <\frac{\tau_1}{3}>$ & & \\
\hline

${\Z}/3 \times {\Z}/3$ & $|\tau_1| =1$, & $\rho$ & $\tau_1 \bar{z_1}$,&
  $-\bar{z_2}$, & $T$ \\
${\cal S}_3 \times {\Z}/3$ &    $\cup \ x_1 = - \frac{1}{2}$         & &
$G = <\frac{1 -\tau_1}{3}> \times <\frac{1 + \tau_1}{3}>$ &$\epsilon =
\frac{1-\rho}{3}$ & \\

& & & or $- \tau_1 \bar{z_1}$ &  & \\

&  & & $G = <\frac{1 +\tau_1}{3}> \times <\frac{1 -
\tau_1}{3}>$ & & \\

&  & & if $|\tau_1| = 1 $ & & \\

&  & &  $\bar{z_1}$ & & \\ &  & &  $G = <\frac{1 -\tau_1}{3}> \times
<\frac{1}{3}>$ & & \\

&  & &  or $-\bar{z_1}$ & & \\

&  & &  $G = <\frac{1}{3}> \times <\frac{1 - \tau_1}{3}>$ & & \\

&  & & if $x_1 = - \frac{1}{2}$ & & \\
\hline

\end{tabular}
$$

If $G = \Z/6$, we take a generator $g$ of $G$ such that $g(z_1, z_2) = (z_1 +
\eta, -\rho z_2)$ and also here the homology of $S(\R)$ is the only
topological invariant needed to distinguish the four cases.
$$
\begin{tabular}{|c|c|c|c|c|c|}
\hline
${\Z}/6$ & $x_1 =0$, & $\rho$ & $\bar{z_1}$, $\eta = \frac{\tau_1}{6}$ &
$\bar{z_2}$, & $2T$ \\
$D_6$ & & & or $- \bar{z_1}$, $\eta = \frac{1}{6}$ & & \\
\hline

${\Z}/6$ & $x_1 =0$, & $\rho$ & $\bar{z_1}$, $\eta = \frac{1}{2} +
\frac{\tau_1}{6}$ &  $ \bar{z_2}$, & $T$ \\
$D_6$ &  & & or $- \bar{z_1}$, $\eta = \frac{\tau_1}{2} +
\frac{1}{6}$  & & \\
\hline

${\Z}/6$ & $x_1 =0$, & $\rho$ & $\bar{z_1} + \frac{1}{2}$, $\eta =
\frac{\tau_1}{6}$ &  $\bar{z_2}$, & $\emptyset$ \\
$D_6$ &  & & or $- \bar{z_1} + \frac{\tau_1}{2}$, $\eta =
\frac{1}{6}$  & & \\
\hline

${\Z}/6$ & $|\tau_1| =1$, & $\rho$ & $\tau_1 \bar{z_1}$, $\eta = \frac{1 -
\tau_1}{6}$ &  $ \bar{z_2}$, & $2K$ \\
$D_6$ & $\cup \ x_1 = - \frac{1}{2}$        & & or $-\tau_1 \bar{z_1}$,
$\eta =
\frac{1+
\tau_1}{6}$ & & \\    &  & & if $|\tau_1| = 1$, & & \\   & & &
$\bar{z_1}$, $
\eta =
\frac{1}{6}  +
\frac{\tau_1}{3}$  & & \\   &   & &  or $- \bar{z_1}$, $ \eta =
\frac{1}{6}$ & &
\\   &  & &  if $x_1 = - \frac{1}{2}$  & & \\
\hline
\end{tabular}
$$

Let us now consider the case where $G = {\Z}/2 \times {\Z}/2$.

Here, we may choose as generators of $G$ the generator
$t$ of $T$ and another element $g$, which is not canonically defined.

  Let $\eta_1, \epsilon_1, \epsilon_2$ be such that $g(z_1, z_2) = (z_1 +
\eta_1, -z_2)$, $t(z_1, z_2) = (z_1 + \epsilon_1, z_2 + \epsilon_2)$.

Here we have the same topological invariants 1), 2), as in the case $G =
\Z/2$.

Furthermore the normal subgroup $T$ of $G$ is of order $2$, let as usual
$t$ be a generator. We may consider all the possible liftings  of $t$ to a
vector $t'$ in the lattice
$\Omega' = \Lambda' \oplus \Gamma$.

The condition whether there exists such a $t'$ whose two components are
eigenvectors for the action of $\tilde{\sigma}$ on
$\Omega'$ is a topological invariant of the real hyperelliptic surface
(notice that the two possible choices for ${\sigma}_2$ differ just up to
multiplication by $-1$).

$$\begin{tabular}{|c|c|c|c|c|c|}
\hline
$G$, $\hat{G}$ & $\tau_1$ & $\tau_2$ & $\sigma_1(z_1), \ \eta_1, \
\epsilon_1$ &
$\sigma_2(z_2), \ \epsilon_2$ & $S({\R})$  \\
\hline
\hline

$({\Z}/2)^2$ &  $x_1 =0$ & $x_2 =0$ & $\bar{z_1}, \eta_1 =
\frac{1}{2}$, $\epsilon_1 = \frac{\tau_1}{2}$ & $\bar{z_2}$, $\epsilon_2 =
\frac{1 + \tau_2}{2}$  & $2K$ \\
$({\Z}/2)^3$  & & &  or $-\bar{z_1},
\eta_1 = \frac{\tau_1}{2}$, $\epsilon_1 = \frac{1}{2}$ & & \\
\hline
$({\Z}/2)^2$ &  $x_1 =0$ & $x_2 =0$ & $\bar{z_1}, \eta_1 =
\frac{1}{2}$, $\epsilon_1= \frac{\tau_1}{2}$ & $\bar{z_2}$, $\epsilon_2 =
\frac{1}{2}$  & $2K$ \\
$({\Z}/2)^3$ &  & &  or $-\bar{z_1},
\eta_1 =
\frac{\tau_1}{2}$, $\epsilon_1 = \frac{1}{2}$ & & \\
\hline

$({\Z}/2)^2$ & $x_1 =0$, & $|\tau_2| =1$ & $\bar{z_1}$, $\eta_1 =
\frac{1}{2}$,
$\epsilon_1 = \frac{\tau_1}{2}$ & $\tau_2 \bar{z_2}$, $\epsilon_2 =
\frac{1 +
\tau_2}{2}$ &
$2K$\\

$({\Z}/2)^3$ &  & or     & or
$- \bar{z_1}$, $\eta_1 = \frac{\tau_1}{2}$, $\epsilon_1 = \frac{1}{2}$ & if &
\\  & &$\ x_2 = - \frac{1}{2}$ & & $|\tau_2| = 1$ & \\

& & & & $ \bar{z_2}$, $\epsilon_2 = \frac{1}{2}$ & \\   & &  & &if &
\\     & & & &
$x_2 = - \frac{1}{2}$, & \\
\hline

$({\Z}/2)^2$ &  $x_1 =0$ & $x_2 =0$ & $\bar{z_1}, \eta_1 =
\frac{1}{2}$, $\epsilon_1 = \frac{1+ \tau_1}{2}$ & $\bar{z_2}$,
$\epsilon_2 =
\frac{\tau_2}{2}$  & $2K$ \\
  $({\Z}/2)^3$ &  & &  or $-\bar{z_1},
\eta_1 = \frac{\tau_1}{2}$, $\epsilon_1 = \frac{1 + \tau_1}{2}$ & & \\
\hline

$({\Z}/2)^2$ &  $x_1 =0$ & $x_2 =0$ & $\bar{z_1}, \eta_1 =
\frac{1}{2}$, $\epsilon_1 = \frac{1+ \tau_1}{2}$ & $\bar{z_2}$,
$\epsilon_2 =
\frac{1 + \tau_2}{2}$  & $2K$ \\
  $({\Z}/2)^3$ & & &  or $-\bar{z_1},
\eta_1 = \frac{\tau_1}{2}$, $\epsilon_1 = \frac{1 + \tau_1}{2}$ & & \\
\hline

$({\Z}/2)^2$ & $x_1 =0$, & $| \tau_2| = 1$, & $\bar{z_1}$, $\eta_1 =
\frac{1}{2}$,
$\epsilon_1 = \frac{1 + \tau_1}{2}$ & $ \tau_2 \bar{z_2}$,
$\epsilon_2 = \frac{1 + \tau_2}{2}$ & $2K$\\

$({\Z}/2)^3$ & &or     & or
$- \bar{z_1}$,
$\eta_1 = \frac{\tau_1}{2}$, $\epsilon_1 = \frac{1 + \tau_1}{2}$ & if &
\\   & &$\ x_2 = -
\frac{1}{2}$ & & $|\tau_2| = 1$ & \\

& & & & $ \bar{z_2}$, $\epsilon_2 = \frac{1}{2}$ & \\    & & & & if &
\\    & & & & $x_2 = - \frac{1}{2}$, & \\
\hline

\end{tabular}
$$

$$\begin{tabular}{|c|c|c|c|c|c|}
\hline

$({\Z}/2)^2$ &  $x_1 =0$ & $x_2 =0$ & $\bar{z_1}, \eta_1 =
\frac{1}{2}$, $\epsilon_1 = \frac{\tau_1}{2}$ & $\bar{z_2} +
\frac{\tau_2}{2}$,
$\epsilon_2 = \frac{1}{2}$  & $T$ \\
  $({\Z}/2)^3$ & & &  or
  & & \\   & & & $-\bar{z_1}, \eta_1 = \frac{\tau_1}{2}$, $\epsilon_1 =
\frac{1}{2}$ & & \\
\hline

$({\Z}/2)^2$ &  $x_1 =0$ & $x_2 =0$ & $\bar{z_1}, \eta_1 = \frac{1}{2}$,
$\epsilon_1 =
\frac{1 + \tau_1}{2}$ & $\bar{z_2} + \frac{\tau_2}{2}$,
$\epsilon_2 = \frac{1}{2}$  & $T$ \\
$({\Z}/2)^3$ &  & &  or & & \\   & & & $-\bar{z_1}, \eta_1 =
\frac{\tau_1}{2}$,
$\epsilon_1 = \frac{1 + \tau_1}{2}$ & & \\
\hline

$({\Z}/2)^2$ &  $x_1 =0$ & $x_2 =0$ & $\bar{z_1}, \eta_1 = \frac{1 +
\tau_1}{2}$,
$\epsilon_1 = \frac{1}{2}$ & $\bar{z_2}$, $\epsilon_2 =
\frac{\tau_2}{2}$  & $T$ \\
$({\Z}/2)^3$  &  & &  or
  & & \\   & & & $-\bar{z_1}, \eta_1 = \frac{1 + \tau_1}{2}$, $\epsilon_1 =
\frac{\tau_1}{2}$ & & \\
\hline

$({\Z}/2)^2$ &  $x_1 =0$ & $x_2 =0$ & $\bar{z_1}, \eta_1 = \frac{1 +
\tau_1}{2}$,
$\epsilon_1 = \frac{1}{2}$ & $\bar{z_2}$, $\epsilon_2 =
\frac{1 + \tau_2}{2}$  & $T$ \\
$({\Z}/2)^3$ & & &  or
  & & \\   & & & $-\bar{z_1}, \eta_1 = \frac{1 + \tau_1}{2}$, $\epsilon_1 =
\frac{\tau_1}{2}$ & & \\
\hline

$({\Z}/2)^2$ &  $x_1 =0$ & $x_2 =0$ & $\bar{z_1}, \eta_1 = \frac{1 +
\tau_1}{2}$,
$\epsilon_1 = \frac{\tau_1}{2}$ & $\bar{z_2}$, $\epsilon_2 =
\frac{1 + \tau_2}{2}$  & $T$ \\
$({\Z}/2)^3$ &  & &  or
  & & \\   & & & $-\bar{z_1}, \eta_1 = \frac{1 + \tau_1}{2}$, $\epsilon_1=
\frac{1}{2}$ & & \\
\hline

\end{tabular}
$$

$$\begin{tabular}{|c|c|c|c|c|c|}
\hline

$({\Z}/2)^2$ &  $x_1 =0$ & $x_2 =0$ & $\bar{z_1}, \eta_1 = \frac{1}{2}$,
$\epsilon_1 =
\frac{\tau_1}{2}$ & $\bar{z_2} + \frac{1}{2}$, $\epsilon_2 =
\frac{\tau_2}{2}$  &
$\emptyset$ \\
$({\Z}/2)^3$ &  & &  or & & \\   & & & $-\bar{z_1}, \eta_1 =
\frac{\tau_1}{2}$,
$\epsilon_1 = \frac{1}{2}$ & & \\
\hline

$({\Z}/2)^2$ &  $x_1 =0$ & $x_2 =0$ & $\bar{z_1}, \eta_1 = \frac{1}{2}$,
$\epsilon_1 =
\frac{1 + \tau_1}{2}$ & $\bar{z_2} +
\frac{1}{2}$, $\epsilon_2 = \frac{\tau_2}{2}$  & $\emptyset$ \\
$({\Z}/2)^3$ &  & &  or
  & & \\   & & &$-\bar{z_1}, \eta_1 = \frac{\tau_1}{2}$, $\epsilon_1=
\frac{1 +
\tau_1}{2}$& & \\
\hline

$({\Z}/2)^2$ &  $x_1 =0$ & $x_2 =0$ & $\bar{z_1}, \eta_1 =
\frac{1 + \tau_1}{2}$, $\epsilon_1 = \frac{1}{2}$ & $\bar{z_2} +
\frac{1}{2}$, $\epsilon_2 = \frac{\tau_2}{2}$  & $\emptyset$ \\
  $({\Z}/2)^3$ &  & &  or
  & & \\   & & & $-\bar{z_1}, \eta_1 = \frac{1 + \tau_1}{2}$, $\epsilon_1 =
\frac{\tau_1}{2}$ & & \\
\hline

$({\Z}/2)^2$ &  $x_1 =0$ & $x_2 =0$ & $\bar{z_1} + \frac{1}{2}, \eta_1 =
\frac{\tau_1}{2}$, $\epsilon_1 = \frac{1 + \tau_1}{2}$ & $\bar{z_2}$,
$\epsilon_2 =
\frac{1 + \tau_2}{2}$  & $\emptyset$ \\
  $({\Z}/2)^3$ & & &  or
  & & \\   & & & $-\bar{z_1} + \frac{\tau_1}{2}, \eta_1 = \frac{1}{2}$,
$\epsilon_1 = \frac{1 + \tau_1}{2}$ & & \\
\hline

$({\Z}/2)^2$ &  $x_1 =0$ & $x_2 =0$ & $\bar{z_1} + \frac{1}{2}, \eta_1 =
\frac{\tau_1}{2}$, $\epsilon_1 = \frac{1}{2}$ & $\bar{z_2}$, $\epsilon_2 =
\frac{1 +
\tau_2}{2}$  & $\emptyset$ \\
  $({\Z}/2)^3$ &  & &  or
  & & \\  & & & $-\bar{z_1} + \frac{\tau_1}{2}, \eta_1 = \frac{1}{2}$,
$\epsilon_1 = \frac{\tau_1}{2}$ & & \\
\hline

$({\Z}/2)^2$ &  $x_1 =0$ & $x_2 =0$ & $\bar{z_1}, \eta_1 =
\frac{\tau_1}{2}$, $\epsilon_1= \frac{1}{2}$ & $\bar{z_2}$,
$\epsilon_2 = \frac{\tau_2}{2}$  & $3T$ \\
$({\Z}/2)^3$ &  & &  or & & \\    & & & $-\bar{z_1}, \eta_1 =
\frac{1}{2}$,
$\epsilon_1 =
\frac{\tau_1}{2}$ & & \\
\hline

\end{tabular}
$$

$$\begin{tabular}{|c|c|c|c|c|c|}
\hline
$({\Z}/2)^2$ &  $x_1 =0$ & $x_2 =0$ & $\bar{z_1}, \eta_1 =
\frac{\tau_1}{2}$, $\epsilon_1 = \frac{1}{2}$ & $\bar{z_2}$,
$\epsilon_2 = \frac{1 + \tau_2}{2}$  & $2T$ \\
  $({\Z}/2)^3$ &  & &  or  & & \\ & & & $-\bar{z_1}, \eta_1 = \frac{1}{2}$,
$\epsilon_1 =
\frac{\tau_1}{2}$ & & \\
\hline

$({\Z}/2)^2$ &  $x_1 =0$ & $x_2 =0$ & $\bar{z_1}, \eta_1 =
\frac{\tau_1}{2}$,
$\epsilon_1 = \frac{1}{2}$ & $\bar{z_2} + \frac{\tau_2}{2}$,
$\epsilon_2 = \frac{1}{2}$  & $2T$ \\
  $({\Z}/2)^3$ &  & &  or  & & \\ & & & $-\bar{z_1}, \eta_1 = \frac{1}{2}$,
$\epsilon_1 =
\frac{\tau_1}{2}$ & & \\
\hline

$({\Z}/2)^2$ & $x_1 =0$, & $| \tau_2| = 1$, & $\bar{z_1}$, $\eta_1 =
\frac{\tau_1}{2}$, $\epsilon_1= \frac{1}{2}$ & $\tau_2
\bar{z_2}$, $\epsilon_2 = \frac{1 + \tau_2}{2}$ & $2T$\\

$({\Z}/2)^3$ & &  or & or & if & \\   & &$x_2 = - \frac{1}{2}$ &$-
\bar{z_1}$,
$\eta_1 = \frac{1}{2}$, $\epsilon_1 = \frac{\tau_1}{2}$ & $|\tau_2| = 1$ & \\

& & & & $ \bar{z_2}$, $\epsilon_2 = \frac{1}{2}$ & \\   & & & & if & \\ & & &
&
$x_2 = -\frac{1}{2}$, & \\
\hline

$({\Z}/2)^2$ &  $x_1 =0$ & $x_2 =0$ & $\bar{z_1}, \eta_1 =
\frac{\tau_1}{2}$, $\epsilon_1 = \frac{1 + \tau_1}{2}$ & $\bar{z_2}$,
$\epsilon_2 =
\frac{1 + \tau_2}{2}$  & $2T$ \\
  $({\Z}/2)^3$ &  & &  or  & & \\ & & &  $-\bar{z_1}, \eta_1 =
\frac{1}{2}$,
$\epsilon_1 =
\frac{1 + \tau_1}{2}$ & & \\
\hline

\end{tabular}
$$
$$\begin{tabular}{|c|c|c|c|c|c|}
\hline
$({\Z}/2)^2$ &  $x_1 =0$ & $x_2 =0$ & $\bar{z_1}, \eta_1 =
\frac{1}{2}$, $\epsilon_1 = \frac{\tau_1}{2}$ & $\bar{z_2}$, $\epsilon_2 =
\frac{\tau_2}{2}$  & $2K \sqcup T$ \\
$({\Z}/2)^3$ & & &  or & & \\    & & & $-\bar{z_1}, \eta_1 =
\frac{\tau_1}{2}$,
$\epsilon_1 =
\frac{1}{2}$ & & \\
\hline

$({\Z}/2)^2$ &  $x_1 =0$ & $x_2 =0$ & $\bar{z_1}, \eta_1 =
\frac{1}{2}$, $\epsilon_1 = \frac{1 + \tau_1}{2}$ & $\bar{z_2}$,
$\epsilon_2 =
\frac{1}{2}$  & $2K \sqcup T$ \\
  $({\Z}/2)^3$ &  & &  or  & &\\   & & & $-\bar{z_1}, \eta_1 =
\frac{\tau_1}{2}$,
$\epsilon_1 = \frac{1 +
\tau_1}{2}$& &
\\
\hline

\end{tabular}
$$

$$\begin{tabular}{|c|c|c|c|c|c|}
\hline
$G$, $\hat{G}$ & $\tau_1$ & $\tau_2$ & $\sigma_1(z_1), \ \eta_1, \
\epsilon_1$ &
$\sigma_2(z_2), \ \epsilon_2$ & $S({\R})$  \\
\hline
\hline

$({\Z}/2)^2 $ & $|\tau_1| =1$, & $x_2 =0$ & $\tau_1 \bar{z_1}$, $\eta_1 =
\frac{1}{2}$, $\epsilon_1 = \frac{1+\tau_1}{2}$ & $ \bar{z_2} +
\frac{\tau_2}{4}$, &
$2T$
\\
$D_4$ & or & & or $-\tau_1 \bar{z_1}$, $\eta_1 =
\frac{\tau_1}{2}$, $\epsilon_1 = \frac{1 + \tau_1}{2}$ & $\epsilon_2 =
\frac{\tau_2}{2}$    or & \\    & $ x_1 = - \frac{1}{2}$ & & if $|\tau_1| =
1$, &$- \bar{z_2} -
\frac{1}{4}$ & \\   & & & $\bar{z_1}$, $ \eta_1 = \frac{\tau_1}{2}$,
$\epsilon_1 =\frac{1}{2}   $ &$\epsilon_2 = \frac{1}{2}$ & \\    &   & & or
$- \bar{z_1}$, $
\eta_1 = \frac{\tau_1}{2}$, $\epsilon_1 = \frac{1}{2}$ & & \\    &  & & if
$x_1 = -
\frac{1}{2}$  & & \\
\hline

$({\Z}/2)^2 $ & $|\tau_1| =1$, & $x_2 =0$ & $\tau_1 \bar{z_1}$, $\eta_1 =
\frac{1}{2}$, $\epsilon_1 = \frac{1+\tau_1}{2}$ & $ \tau_2 \bar{z_2} +
\frac{1-\tau_2}{4}$, & $2T$ \\
$D_4$ & or & & or $-\tau_1 \bar{z_1}$, $\eta_1 =
\frac{\tau_1}{2}$, $\epsilon_1 = \frac{1 + \tau_1}{2}$ & $\epsilon_2 =
\frac{1+\tau_2}{2}$    or & \\     & $x_1 = - \frac{1}{2}$ & & if
$|\tau_1| = 1$, &$-
\tau_2 \bar{z_2} +\frac{1+ \tau_2}{4}$ & \\   & & & $\bar{z_1}$, $ \eta_1 =
\frac{\tau_1}{2}$, $\epsilon_1 =\frac{1}{2}   $ &$\epsilon_2 = \frac{1+
\tau_2}{2}$ &
\\   &   & &  or $- \bar{z_1}$, $ \eta_1 = \frac{\tau_1}{2}$, $\epsilon_1 =
\frac{1}{2}$ &or
$\bar{z_2} + \frac{1}{4} + \frac{\tau_2}{2}$ & \\   &  & &  if $x_1 = -
\frac{1}{2}$  &$\epsilon_2 = \frac{1}{2}$ & \\ &  & &  & or $- \bar{z_2} -
\frac{1}{4}$& \\ &  & &  & or $\epsilon_2 = \frac{1}{2}$& \\ &  & &  & if
$x_2 = -
\frac{1}{2}$& \\
\hline

$({\Z}/2)^2 $ & $|\tau_1| =1$, & $x_2 =0$ & $\tau_1 \bar{z_1}$, $\eta_1 =
\frac{1}{2}$, $\epsilon_1 = \frac{1+\tau_1}{2}$ & $ \bar{z_2} +
\frac{1}{2} +
\frac{\tau_2}{4}$, & $\emptyset$ \\
$D_4$ & or & & or $-\tau_1 \bar{z_1}$, $\eta_1 =
\frac{\tau_1}{2}$, $\epsilon_1 = \frac{1 + \tau_1}{2}$ & $\epsilon_2 =
\frac{\tau_2}{2}$    or & \\    & $x_1 = - \frac{1}{2}$  & & if $|\tau_1| =
1$, &$- \bar{z_2} -
\frac{1}{4}+ \frac{\tau_2}{2}$ & \\   & & & $\bar{z_1}$, $ \eta_1 =
\frac{\tau_1}{2}$,
$\epsilon_1 =\frac{1}{2}   $ & $\epsilon_2 = \frac{1}{2}$ & \\   &   & & or $-
\bar{z_1}$,
$ \eta_1 = \frac{\tau_1}{2}$, $\epsilon_1 = \frac{1}{2}$ & & \\   &  & & if
$x_1 = -
\frac{1}{2}$  & & \\
\hline
\end{tabular}
$$

Let us now explain how to distinguish the different cases listed in the
tables above.
\bigskip  Assume first of all $\hat{G} = \Z/2 \times \Z/2 \times \Z/2$.

$S(\R) = 2K$: the first case is distinguished from the second and the fourth
since in the first case $\epsilon_2 = \frac{1+ \tau_2}{2}$, therefore any
lifting of $t$ to a vector $t' \in \Omega'$ cannot be an eigenvector for the
action of $\sigma_2$, while in the second and the fourth case $\epsilon_2 =
\frac{1}{2}$, respectively $\frac{\tau_2}{2}$, which are eigenvectors for the
action of $\sigma_2$.

The first case is distinguished from the third and the sixth by the parity of
$\nu(\sigma_2)$. In fact $\nu(\sigma_2) = 2$ in the first case, while in the
third and the sixth case we have $\nu(\sigma_2) = 1$.

The first case is distinguished from the fifth since in the first case
$\epsilon_1$ is an eigenvector for the action of $\sigma_1$ on
$\Lambda'$, while in the fifth case $\epsilon_1 = \frac{1 + \tau_1}{2}$,
which is not an eigenvector for the action of $\sigma_1$.

The second case is distinguished from the third and the sixth by the parity
of $\nu(\sigma_2)$. In fact $\nu(\sigma_2) = 2$ in the second case, while in
the third and the sixth case we have $\nu(\sigma_2) = 1$.

The second case is distinguished from the fourth since in the second case
$\epsilon_1$ is an eigenvector for the action of $\sigma_1$ on
$\Lambda'$, while in the fourth case $\epsilon_1 = \frac{1 + \tau_1}{2}$,
which is not an eigenvector for the action of $\sigma_1$.

The second case is distinguished from the fifth since in the second case
$\epsilon_2$ is an eigenvector for the action of $\sigma_2$ on $\Gamma$,
while in the fifth case $\epsilon_2 = \frac{1 + \tau_2}{2}$, which is not an
eigenvector for the action of $\sigma_2$.

The third case is distinguished from the fourth and the fifth by the parity
of $\nu(\sigma_2)$. In fact $\nu(\sigma_2) = 1$ in the third case, while in
the other two cases we have $\nu(\sigma_2) = 2$.

The third case is distinguished from the sixth since in the third case
$\epsilon_1$ is an eigenvector for the action of $\sigma_1$ on
$\Lambda'$, while in the sixth case $\epsilon_1 = \frac{1 + \tau_1}{2}$,
which is not an eigenvector for the action of $\sigma_1$

The fourth and the fifth cases are distinguished from the sixth by the parity
of $\nu(\sigma_2)$. In fact $\nu(\sigma_2) = 1$ in the sixth case, while in
the other two cases we have $\nu(\sigma_2) = 2$.

The fourth case is distinguished from the fifth since in the fourth case
$\epsilon_2$ is an eigenvector for the action of $\sigma_2$ on $\Gamma$,
while in the fifth case $\epsilon_2 = \frac{1 + \tau_2}{2}$, which is not an
eigenvector for the action of $\sigma_2$.

\bigskip

$S(\R) = T$: the first two cases are distinguished from the the last three
cases by the invariant 1)(b):  in the first two cases the set of values of
$\nu(\sigma_2)$ equals $\{2,0\}$, in the last three cases it equals $\{2,2\}$.

The first case is distinguished from the second since in the first case
$\epsilon_1$ is an eigenvector for the action of $\sigma_1$ on
$\Lambda'$, while in the second case $\epsilon_1 = \frac{1 + \tau_1}{2}$,
which is not an eigenvector for the action of $\sigma_1$.

The third case is distinguished from the fourth and the fifth since in the
third case $\epsilon_2$ is an eigenvector for the action of
$\sigma_2$ on $\Gamma$, while in the other two cases $\epsilon_2 =
\frac{1 + \tau_2}{2}$, which is not an eigenvector for the action of
$\sigma_2$.

The fourth case is distinguished from the fifth since in the fourth case
$\epsilon_1$ is the $+1$ eigenvector for the action of $\sigma_1$ on
$\Lambda'$, while in the fifth case $\epsilon_1$ is the $-1$ eigenvector for
the action of $\sigma_1$.

\bigskip

$S(\R) = \emptyset$: the first three cases are distinguished from the the
last two cases by the invariant 1)(b):  in the first three cases the set of
values of $\nu(\sigma_2)$ equals $\{0,2\}$, in the last two cases it equals
$\{2,2\}$.

The second case is distinguished from the first and the third since in the
first and the third case $\epsilon_1$ is an eigenvector for the action of
$\sigma_1$ on $\Lambda'$, while in the second case $\epsilon_1 = \frac{1 +
\tau_1}{2}$, which is not an eigenvector for the action of
$\sigma_1$.

The first case is distinguished from the third since in the first case
$\epsilon_1$ is the $-1$ eigenvector for the action of $\sigma_1$ on
$\Lambda'$, while in the third case $\epsilon_1$ is the $+1$ eigenvector for
the action of $\sigma_1$.

The fourth case is distinguished from the fifth since in the fifth case
$\epsilon_1$ is an eigenvector for the action of $\sigma_1$ on
$\Lambda'$, while in the fourth case $\epsilon_1 = \frac{1 + \tau_1}{2}$,
which is not an eigenvector for the action of $\sigma_1$.

\bigskip

$S(\R) = 2T$: the third case is distinguished from all the others by the
parity of $\nu(\sigma_2) = 1$.

The second case is distinguished from all the others by the invariant 1)(b).
In fact the set of values of $\nu(\sigma_2)$ is $\{2,0\}$ in the second case,
while in the other cases it is either $\{2,2\}$, or
$\{1,1\}$ (only in the third case).

The first case is distinguished from the last one since in the first case
$\epsilon_1$ is an eigenvector for the action of $\sigma_1$ on
$\Lambda'$, while in the last case $\epsilon_1 = \frac{1 + \tau_1}{2}$, which
is not an eigenvector for the action of $\sigma_1$.

\bigskip

$S(\R) = 2K \sqcup T$: the first case is distinguished from the second one
since in the first case $\epsilon_1$ is an eigenvector for the action of
$\sigma_1$ on $\Lambda'$, while in the second case $\epsilon_1 =
\frac{1 + \tau_1}{2}$, which is not an eigenvector for the action of
$\sigma_1$.

\bigskip

If $\hat{G} = D_4$, the cases with $S(\R) = 2T$ are distinguished by the
parity of $\nu(\sigma_2)$, which is equal to $2$ in the first case, and equal
to $1$ in the second case.
\bigskip

We finally give the table for the non split case, where $G = {\Z}/2
\times {\Z}/2 $ is generated by elements $g,t$ such that $g(z_1, z_2) = (z_1
+ \eta_1, - z_2)$, $g_2(z_1, z_2) = (z_1 + \epsilon_1, z_2 +
\epsilon_2)$.  Recall the by now standard notation $\tau_j = x_j + i y_j$.

$$
\begin{tabular}{|c|c|c|c|c|c|}
\hline
$G$, $\hat{G}$ & $\tau_1$ & $\tau_2$ & $\tilde{\sigma_1}$, $\eta_1$,
$\epsilon_1$ &
$\tilde{\sigma_2}$, $\epsilon_2$ & $S({\R})$\\
\hline
\hline
$({\Z}/2)^2$ & $x_1 =0$, & $|\tau_2| = 1$ & $\bar{z_1} + \frac{1}{4}$, &
$
\tau_2
\bar{z_2} + \frac{1}{2}$, & $\emptyset$ \\
${\Z}/4 \times {\Z}/2$ & & or &$\epsilon_1 = \frac{1}{2}$    &
$\epsilon_2 =
\frac{1 +
\tau_2}{2}$    & \\   &   &$ x_2 = - \frac{1}{2}$ & $\eta_1 = \frac{1+
\tau_1}{2}$ & if
$|\tau_2| =1$& \\  &   & &or $- \bar{z_1} + \frac{\tau}{4}$   &$\bar{z_2} +
\frac{\tau}{2}$ &\\  & & & $\epsilon_1= \frac{\tau_1}{2}$ & $\epsilon_2 =
\frac{1}{2}$ &
\\  & & & $\eta_1 = \frac{1 + \tau_1}{2}$ & if $x_2 = - \frac{1}{2}$ & \\

\hline

\end{tabular}
$$

\bigskip

Address of the authors:\\

Fabrizio Catanese:\\ Mathematisches Institut der Georg - August -
Universit\"at G\"ottingen\\ Bunsentra\ss e 3-5  D-37073 G\"ottingen

e-mail: catanese@@uni-math.gwdg.de\\

Paola Frediani: \\ Dipartimento di Matematica della Universit\`a di Pisa\\
via Buonarroti, 2 I-56127 Pisa

e-mail: frediani@@dm.unipi.it\\

\end{document}